\documentstyle[11pt,leqno]{article}
\makeatletter
\ifx\@auxdefsloaded\relax\endinput
\else\let\@auxdefsloaded\relax\fi

\def\newenvironment{%
   \@ifnextchar *{\@@newenv{\global\@ignoretrue}}{\@@newenv{}*}}

\def\@@newenv#1*#2{%
   \@ifnextchar [{\@newenv{#1}{#2}}{\@newenv{#1}{#2}[0]}}

\long\def\@newenv#1#2[#3]#4#5{%
   \expandafter\newcommand\csname#2\endcsname[#3]{#4}%
   \expandafter\long\expandafter\def\csname end#2\endcsname{#5#1}}

\def\renewenvironment{%
   \@ifnextchar *{\@@renewenv{\global\@ignoretrue}}{\@@renewenv{}*}}

\def\@@renewenv#1*#2{%
   \@ifnextchar [{\@renewenv{#1}{#2}}{\@renewenv{#1}{#2}[0]}}

\long\def\@renewenv#1#2[#3]#4#5{%
   \expandafter\renewcommand\csname#2\endcsname[#3]{#4}%
   \expandafter\long\expandafter\def\csname end#2\endcsname{#5#1}}

\def\newoptcommand#1#2{%
   \@ifnextchar [{\@optargdef#1#2}{\@optargdef#1#2[1]}}

\def\renewoptcommand#1#2{%
   \edef\@tempa{\expandafter\@cdr\string#1\@nil}%
   \@ifundefined{\@tempa}{%
      \@latexerr{\string#1\space undefined}\@ehc}{}%
   \@ifnextchar [{\@reoptargdef#1#2}{\@reoptargdef#1#2[1]}}

\long\def\@optargdef#1#2[#3]#4{%
   \@ifdefinable #1{\@reoptargdef#1#2[#3]{#4}}}

\long\def\@reoptargdef#1#2[#3]#4{%
   \@tempcnta#3\relax \@tempcntb \@ne
   \let#1\relax \let\@tempa\relax
   \edef\@tempb{\long\def\csname\string#1\endcsname
      [\@tempa\the\@tempcntb]}%
   \advance\@tempcntb \@ne \advance\@tempcnta \m@ne
   \@whilenum\@tempcnta>0\do{%
      \edef\@tempb{\@tempb\@tempa\the\@tempcntb}%
      \advance\@tempcntb \@ne \advance\@tempcnta \m@ne}%
   \let\@tempa=##\@tempb{#4}%
   \def#1{\@ifnextchar [{\csname\string#1\endcsname}{%
      \csname\string#1\endcsname[#2]}}}

\def\newoptenvironment{%
   \@ifnextchar *{\@@newoptenv{\global\@ignoretrue}}{%
      \@@newoptenv{}*}}

\def\@@newoptenv#1*#2#3{%
   \@ifnextchar [{\@newoptenv{#1}{#2}{#3}}{%
      \@newoptenv{#1}{#2}{#3}[0]}}

\long\def\@newoptenv#1#2#3[#4]#5#6{%
   \expandafter\newoptcommand\csname#2\endcsname{#3}[#4]{#5}%
   \expandafter\long\expandafter\def\csname end#2\endcsname{#6#1}}

\def\renewoptenvironment{%
   \@ifnextchar *{\@@renewoptenv{\global\@ignoretrue}}{%
      \@@renewoptenv{}*}}

\def\@@renewoptenv#1*#2#3{%
   \@ifnextchar [{\@renewoptenv{#1}{#2}{#3}}{%
      \@renewoptenv{#1}{#2}{#3}[0]}}

\long\def\@renewoptenv#1#2#3[#4]#5#6{%
   \expandafter\renewoptcommand\csname#2\endcsname{#3}[#4]{#5}%
   \expandafter\long\expandafter\def\csname end#2\endcsname{#6#1}}

\newcounter{keepoptional}
\newcounter{optnestctr}

\newoptenvironment*{optional}{1}[1]{%
   \ifnum#1>\value{keepoptional}\relax
      \setcounter{optnestctr}{0}\@powerup\expandafter\@endopt\fi
   \ignorespaces}{}

\def\@powerup{\catcode`\{=12 \catcode`\}=12 \catcode`\\=12 \relax}
\def\@powerdown{\catcode`\{=1 \catcode`\}=2 \catcode`\\=0 \relax}
\begingroup \catcode`|=0 \catcode`[=1 \catcode`]=2 \@powerup
   |long|gdef|@endopt#1\end{optional}[%
      |@beginopt#1\begin{optional}*]%
   |long|gdef|@beginopt#1\begin{optional}[%
      |@ifstar[%
         |ifnum|value[optnestctr]>0|relax
            |addtocounter[optnestctr][-1]|let|@zzzap=|@endopt
         |else
            |@powerdown|end[optional]|let|@zzzap=|relax|fi
	 |@zzzap]%
      [
         |addtocounter[optnestctr][1]|@beginopt]]%
|endgroup

\makeatother
\makeatletter
\ifx\@auxdefsloaded\relax
\else \input{auxdefs.sty}\fi
\newskip\dgARROWLENGTH  \dgARROWLENGTH=2.5em\relax
\newskip\dgHORIZPAD     \dgHORIZPAD=1em\relax
\newskip\dgVERTPAD      \dgVERTPAD=2ex\relax
\newskip\dgLABELOFFSET  \dgLABELOFFSET=.7ex\relax
\newcount\dgARROWPARTS  \dgARROWPARTS=4\relax
\newcommand{\dgeverynode}{\displaystyle}
\newcommand{\dgeverylabel}{\scriptstyle}
\newskip\dgDOTSPACING   \dgDOTSPACING=0.35em
\newskip\dgDOTSIZE      \dgDOTSIZE=1.5\fontdimen8\tenln
\newcount\dgMAXSQUARE   \dgMAXSQUARE=4\relax
\newcount\dgMINDBLSQ    \dgMINDBLSQ=10\relax
\newskip\dgCOLUMNWIDTH  \dgCOLUMNWIDTH=2em\relax
\chardef\f@ur=4
\@namedef{dgo@}{\let\dg@VECTOR=\vector
   \dg@LBLPOS=\dgARROWPARTS \divide\dg@LBLPOS\tw@}%
\@namedef{dgo@-}{\let\dg@VECTOR=\line}%
\@namedef{dgo@!}{\let\dg@VECTOR=\dg@novector}%
\@namedef{dgo@1}{\dg@LBLPOS=\@ne\relax}%
\@namedef{dgo@2}{\dg@LBLPOS=\tw@\relax}%
\@namedef{dgo@3}{\dg@LBLPOS=\thr@@\relax}%
\@namedef{dgo@4}{\dg@LBLPOS=\f@ur\relax}%
\@namedef{dgo@5}{\dg@LBLPOS=5\relax}%
\@namedef{dgo@6}{\dg@LBLPOS=6\relax}%
\@namedef{dgo@7}{\dg@LBLPOS=7\relax}%
\@namedef{dgo@8}{\dg@LBLPOS=8\relax}%
\@namedef{dgo@9}{\dg@LBLPOS=9\relax}%
\def\dgt@e{\dg@DX=\@ne \dg@DY=\z@ \dg@SIZE=\@ne}%
\def\dgt@w{\dg@DX=\m@ne \dg@DY=\z@ \dg@SIZE=\@ne}%
\def\dgt@n{\dg@DX=\z@ \dg@DY=\@ne \dg@SIZE=\@ne}%
\def\dgt@s{\dg@DX=\z@ \dg@DY=\m@ne \dg@SIZE=\@ne}%
\def\dgt@ne{\dg@DX=\@ne \dg@DY=\@ne \dg@SIZE=\@ne}%
\def\dgt@se{\dg@DX=\@ne \dg@DY=\m@ne \dg@SIZE=\@ne}%
\def\dgt@nw{\dg@DX=\m@ne \dg@DY=\@ne \dg@SIZE=\@ne}%
\def\dgt@sw{\dg@DX=\m@ne \dg@DY=\m@ne \dg@SIZE=\@ne}%
\def\dgt@nne{\dg@DX=\@ne \dg@DY=\tw@ \dg@SIZE=\@ne}%
\def\dgt@nnw{\dg@DX=\m@ne \dg@DY=\tw@ \dg@SIZE=\@ne}%
\def\dgt@sse{\dg@DX=\@ne \dg@DY=-\tw@ \dg@SIZE=\@ne}%
\def\dgt@ssw{\dg@DX=\m@ne \dg@DY=-\tw@ \dg@SIZE=\@ne}%
\def\dgt@ene{\dg@DX=\tw@ \dg@DY=\@ne \dg@SIZE=\tw@}%
\def\dgt@ese{\dg@DX=\tw@ \dg@DY=\m@ne \dg@SIZE=\tw@}%
\def\dgt@wnw{\dg@DX=-\tw@ \dg@DY=\@ne \dg@SIZE=\tw@}%
\def\dgt@wsw{\dg@DX=-\tw@ \dg@DY=\m@ne \dg@SIZE=\tw@}%
\def\dggeometry{
   \dg@ZTEMP=\dg@XGRID \multiply\dg@ZTEMP\tw@
   \ifnum\dg@YGRID=\z@ \dg@ZTEMP=\tw@
   \else \divide\dg@ZTEMP\dg@YGRID \fi
   \ifnum\dg@ZTEMP>\f@ur \dg@ZTEMP=\f@ur \fi
   \ifnum\dg@ZTEMP<\@ne \dg@ZTEMP=\@ne \fi
   \unitlength=2sp\relax
   \ifnum\dg@ZTEMP<\tw@
      \advance\dg@ZTEMP\@ne
      \multiply\unitlength\dg@YGRID
   \else
      \multiply\unitlength\dg@XGRID \divide\unitlength\dg@ZTEMP
   \fi
   \dg@XGRID=\dg@ZTEMP \dg@YGRID=\tw@
   \dg@rmcommondiv\tw@\dg@XGRID\dg@YGRID
   \divide\unitlength\dg@YGRID \divide\unitlength\@m\relax}
\@namedef{dgo@..}{\let\dg@VECTOR=\dg@dotvector}%

\def\dg@dotvector(#1,#2)#3{%
   \begingroup
   \dg@XTEMP=#1\relax \dg@YTEMP=#2\relax
   \let\dg@NDOTS=\dg@XEND \let\dg@DOTDIAM=\dg@WEND
   \dg@NDOTS=\unitlength \multiply\dg@NDOTS #3\relax
   \dg@ZTEMP=\dg@YTEMP \dg@changesign\dg@YTEMP\dg@ZTEMP
   \ifnum\dg@XTEMP>\z@
      \ifnum\dg@YTEMP>\dg@XTEMP
         \multiply\dg@NDOTS\dg@YTEMP \divide\dg@NDOTS\dg@XTEMP \fi
   \else\ifnum\dg@XTEMP<\z@
      \ifnum\dg@YTEMP>-\dg@XTEMP
         \multiply\dg@NDOTS\dg@YTEMP \divide\dg@NDOTS-\dg@XTEMP \fi
   \fi\fi
   \dg@YTEMP=\dg@ZTEMP
   \divide\dg@NDOTS\dgDOTSPACING
   \ifnum\dg@NDOTS>\z@\else \dg@NDOTS=\@ne \fi
   \dg@ZTEMP=\unitlength \multiply\dg@ZTEMP #3\relax
   \divide\dg@ZTEMP\dg@NDOTS
   \ifnum\dg@XTEMP=\z@
      \dg@changesign\dg@ZTEMP\dg@YTEMP \dg@YTEMP=\dg@ZTEMP
   \else
      \dg@changesign\dg@ZTEMP\dg@XTEMP
      \multiply\dg@YTEMP\dg@ZTEMP \divide\dg@YTEMP\dg@XTEMP
      \dg@XTEMP=\dg@ZTEMP
   \fi
   \divide\dg@XTEMP\unitlength \divide\dg@YTEMP\unitlength
   \begin{picture}(0,0)
      \dg@DOTDIAM=\dgDOTSIZE \divide\dg@DOTDIAM\unitlength
      \multiput(0,0)(\dg@XTEMP,\dg@YTEMP){\dg@NDOTS}{%
         \circle*{\dg@DOTDIAM}}%
      \multiply\dg@XTEMP\dg@NDOTS \multiply\dg@YTEMP\dg@NDOTS
      \put(\dg@XTEMP,\dg@YTEMP){\vector(#1,#2){0}}%
   \end{picture}%
   \endgroup}%
\newif\ifdg@LATEXGEOM \dg@LATEXGEOMfalse

\@ifundefined{lamsvector}{}{%
   \@namedef{dgo@}{%
      \let\dg@VECTOR=\lamsvector
      \lamsshaft{?}\lamssource{?}\lamstarget{?}%
      \dg@LBLPOS=\dgARROWPARTS \divide\dg@LBLPOS\tw@\relax}%
   \@namedef{dgo@..}{\lamsshaft{.}}%
   \@namedef{dgo@=}{\lamsshaft{=}}%
   \@namedef{dgo@-}{\lamstarget{0}}%
   \@namedef{dgo@'}{\lamstarget{'}}%
   \@namedef{dgo@`}{\lamstarget{`}}%
   \@namedef{dgo@A}{\lamstarget{A}}%
   \@namedef{dgo@V}{\lamssource{V}}%
   \@namedef{dgo@J}{\lamssource{J}}%
   \@namedef{dgo@L}{\lamssource{L}}%
   \@namedef{dgo@S}{\lamssource{S}}%
   \def\dggeometry{
      \dg@ZTEMP=\dg@XGRID \multiply\dg@ZTEMP\tw@
      \ifnum\dg@YGRID=\z@ \dg@ZTEMP=\tw@
      \else \divide\dg@ZTEMP\dg@YGRID \fi
      \ifnum\dg@ZTEMP>6\relax \dg@ZTEMP=6\relax \fi
      \ifdg@LATEXGEOM\ifnum\dg@ZTEMP>\f@ur \dg@ZTEMP=\f@ur \fi\fi
      \ifnum\dg@ZTEMP<\@ne \dg@ZTEMP=\@ne \fi
      \unitlength=2sp\relax
      \ifnum\dg@ZTEMP<\tw@
         \advance\dg@ZTEMP\@ne
         \multiply\unitlength\dg@YGRID
      \else
         \multiply\unitlength\dg@XGRID \divide\unitlength\dg@ZTEMP
      \fi
      \dg@XGRID=\dg@ZTEMP \dg@YGRID=\tw@
      \dg@rmcommondiv\tw@\dg@XGRID\dg@YGRID
      \divide\unitlength\dg@YGRID \divide\unitlength\@m
      \dg@LATEXGEOMfalse}
   \def\dgt@nee{\dg@DX=\tw@ \dg@DY=\@ne \dg@SIZE=\tw@}%
   \def\dgt@see{\dg@DX=\tw@ \dg@DY=\m@ne \dg@SIZE=\tw@}%
   \def\dgt@nww{\dg@DX=-\tw@ \dg@DY=\@ne \dg@SIZE=\tw@}%
   \def\dgt@sww{\dg@DX=-\tw@ \dg@DY=\m@ne \dg@SIZE=\tw@}%
   \def\dgt@nnne{\dg@DX=\@ne \dg@DY=\thr@@ \dg@SIZE=\@ne}%
   \def\dgt@nnnw{\dg@DX=\m@ne \dg@DY=\thr@@ \dg@SIZE=\@ne}%
   \def\dgt@sssw{\dg@DX=\m@ne \dg@DY=-\thr@@ \dg@SIZE=\@ne}%
   \def\dgt@ssse{\dg@DX=\@ne \dg@DY=-\thr@@ \dg@SIZE=\@ne}%
   \def\dgt@nnnee{\dg@DX=\tw@ \dg@DY=\thr@@ \dg@SIZE=\tw@}%
   \def\dgt@nnnww{\dg@DX=-\tw@ \dg@DY=\thr@@ \dg@SIZE=\tw@}%
   \def\dgt@sssww{\dg@DX=-\tw@ \dg@DY=-\thr@@ \dg@SIZE=\tw@}%
   \def\dgt@sssee{\dg@DX=\tw@ \dg@DY=-\thr@@ \dg@SIZE=\tw@}%
   \def\dgt@nneee{\dg@DX=\thr@@ \dg@DY=\tw@ \dg@SIZE=\thr@@}%
   \def\dgt@nnwww{\dg@DX=-\thr@@ \dg@DY=\tw@ \dg@SIZE=\thr@@}%
   \def\dgt@sswww{\dg@DX=-\thr@@ \dg@DY=-\tw@ \dg@SIZE=\thr@@}%
   \def\dgt@sseee{\dg@DX=\thr@@ \dg@DY=-\tw@ \dg@SIZE=\thr@@}%
   \def\dgt@neee{\dg@DX=\thr@@ \dg@DY=\@ne \dg@SIZE=\thr@@
      \global\dg@LATEXGEOMtrue}%
   \def\dgt@nwww{\dg@DX=-\thr@@ \dg@DY=\@ne \dg@SIZE=\thr@@
      \global\dg@LATEXGEOMtrue}%
   \def\dgt@swww{\dg@DX=-\thr@@ \dg@DY=\m@ne \dg@SIZE=\thr@@
      \global\dg@LATEXGEOMtrue}%
   \def\dgt@seee{\dg@DX=\thr@@ \dg@DY=\m@ne \dg@SIZE=\thr@@
      \global\dg@LATEXGEOMtrue}%
}
\newcount\dg@HORIZ      \newcount\dg@VERT
\newcount\dg@XLPAD      \newcount\dg@YBPAD
\newcount\dg@XRPAD      \newcount\dg@YTPAD
\newcount\dg@X          \newcount\dg@Y
\newcount\dg@XGRID      \newcount\dg@YGRID
\newcount\dg@SIZE
\newcount\dg@USERSIZE
\newcount\dg@DX         \newcount\dg@DY
\newcount\dg@XLBL       \newcount\dg@YLBL
\newcount\dg@XOFFSET    \newcount\dg@YOFFSET
\newcount\dg@LBLOFF
\newcount\dg@XLBLOFF    \newcount\dg@YLBLOFF
\newcount\dg@LBLPOS
\newcount\dg@XTEMP      \newcount\dg@YTEMP
\newcount\dg@ZTEMP
\newcount\dg@XNODE      \newcount\dg@YNODE
\newcount\dg@XEND       \newcount\dg@YEND
\newcount\dg@WEND       \newcount\dg@HEND
\newcount\dg@COUNT
\newbox\dg@NODEBOX
\newoptenvironment*{diagram}{}[1]{%
   \global\advance\dg@COUNT\@ne \typeout{[diagram \the\dg@COUNT:}%
   \let\node=\dg@node \let\\=\dg@cr \let\arrow=\dg@arrow
   \def\dg@BIGNODE{#1}%
   \ifx\dg@BIGNODE\@empty\else
      \setbox\dg@NODEBOX=\hbox{$\dgeverynode{#1}$}%
      \dg@XTEMP=\wd\dg@NODEBOX
      \dg@YTEMP=\ht\dg@NODEBOX \advance\dg@YTEMP\dp\dg@NODEBOX
      \ifnum\dg@YTEMP=\z@ \dg@ZTEMP=\z@
      \else \dg@ZTEMP=\dg@XTEMP \divide\dg@ZTEMP\dg@YTEMP\fi
      \ifnum\dg@ZTEMP>\dgMAXSQUARE
         \ifnum\dg@ZTEMP<\dgMINDBLSQ \dg@XGRID=\thr@@ \dg@YGRID=\tw@
         \else \dg@XGRID=\tw@ \dg@YGRID=\@ne \fi
      \else \dg@XGRID=\@ne \dg@YGRID=\@ne \fi
      \advance\dg@XTEMP\dgHORIZPAD \advance\dg@YTEMP\dgVERTPAD
      \dg@XLPAD=\dg@XTEMP \divide\dg@XLPAD\tw@
      \dg@XRPAD=\dg@XTEMP \divide\dg@XRPAD\tw@
      \dg@YBPAD=\dg@YTEMP \divide\dg@YBPAD\tw@
      \dg@YTPAD=\dg@YTEMP \divide\dg@YTPAD\tw@
      \advance\dg@XTEMP\dgARROWLENGTH
      \divide\dg@XTEMP\dg@XGRID \divide\dg@XTEMP\@m
      \unitlength=1sp\relax \multiply\unitlength\dg@XTEMP
   \fi
   \typeout{saving...}%
   \dg@X=-\@ne \dg@Y=\z@ \dg@HORIZ=\z@\relax
   \let\dg@SLIST=\@empty
   \let\dg@NLIST=\@empty \let\dg@ALIST=\@empty
   \let\dg@PASS=\dg@savepass
}{
   \dg@VERT=-\dg@Y\relax
   \typeout{calculating...}%
   \ifx\dg@BIGNODE\@empty
      \dg@XGRID=\z@ \dg@YGRID=\z@
      \dg@XLPAD=\z@ \dg@XRPAD=\z@ \dg@YBPAD=\z@ \dg@YTPAD=\z@
      \let\dg@PASS=\dg@geompass
      \dg@NLIST \dg@ALIST
      \ifnum\dg@XGRID>\z@\else \dg@XGRID=\quad\relax \fi
      \ifnum\dg@YGRID>\z@\else \dg@YGRID=\quad\relax \fi
      \dggeometry
      \ifdim\unitlength>\z@\else \unitlength=\quad \fi
      \ifnum\dg@XGRID>\z@\else \dg@XGRID=\@ne \fi
      \ifnum\dg@YGRID>\z@\else \dg@YGRID=\@ne \fi
   \fi
   \dg@LBLOFF=\dgLABELOFFSET \divide\dg@LBLOFF\unitlength
   \multiply\dg@HORIZ\@m \multiply\dg@HORIZ\dg@XGRID
   \multiply\dg@VERT\@m \multiply\dg@VERT\dg@YGRID
   \divide\dg@XLPAD\unitlength \divide\dg@XRPAD\unitlength
   \divide\dg@YBPAD\unitlength \divide\dg@YTPAD\unitlength
   \typeout{drawing...}%
   \let\dg@PASS=\dg@drawpass
   \hspace*{\dg@XLPAD\unitlength}%
   \vcenter{%
      \hsize=0pt\relax\parindent=0pt\relax
      \parskip=0pt\relax\baselineskip=0pt\relax
      \vspace*{\dg@YTPAD\unitlength}%
      \begin{picture}(0,0)(0,0)\dg@NLIST\dg@ALIST\end{picture}%
      \raisebox{\z@}[\z@][\dg@VERT\unitlength]{}%
      \vspace*{\dg@YBPAD\unitlength}%
   }
   \hspace*{\dg@HORIZ\unitlength}%
   \hspace*{\dg@XRPAD\unitlength}%
   \typeout{done]}}%

\def\dg@savepass{s}
\def\dg@geompass{g}
\def\dg@drawpass{d}

\newoptcommand{\dg@node}{\@ne}[2]{%
   \ifx\dg@PASS\dg@savepass
      \dg@XTEMP=#1\relax
      \ifnum\dg@XTEMP<\@ne \dg@XTEMP=\@ne\fi
      \advance\dg@X\dg@XTEMP
      \ifnum\dg@HORIZ<\dg@X \dg@HORIZ=\dg@X \fi
      \setbox\dg@NODEBOX=\hbox{$\dgeverynode{#2}$}%
      \dg@XTEMP=\wd\dg@NODEBOX \advance\dg@XTEMP\dgHORIZPAD
      \dg@YTEMP=\ht\dg@NODEBOX \advance\dg@YTEMP\dp\dg@NODEBOX
      \advance\dg@YTEMP\dgVERTPAD
      \toks\z@=\expandafter{\dg@SLIST}%
      \edef\dg@SLIST{\the\toks\z@
         ,\noexpand\dg@XNODE=\number\dg@X\noexpand\relax
         \noexpand\dg@YNODE=\number\dg@Y\noexpand\relax
         \noexpand\dg@XTEMP=\number\dg@XTEMP\noexpand\relax
         \noexpand\dg@YTEMP=\number\dg@YTEMP\noexpand\relax}%
      \toks\z@=\expandafter{\dg@NLIST}%
      \toks\tw@={\dg@node{#2}}%
      \edef\dg@NLIST{\the\toks\z@
         \noexpand\dg@X=\number\dg@X\noexpand\relax
         \noexpand\dg@Y=\number\dg@Y\noexpand\relax
         \the\toks\tw@}%
   \else\ifx\dg@PASS\dg@geompass
      \ifnum\dg@X=\z@
         \dg@getnodesize
            {\dg@SLIST}{\dg@X}{\dg@Y}{\dg@WEND}{\dg@HEND}%
         \divide\dg@WEND\tw@
         \ifnum\dg@XLPAD<\dg@WEND \dg@XLPAD=\dg@WEND \fi\fi
      \ifnum\dg@X=\dg@HORIZ
         \dg@getnodesize
            {\dg@SLIST}{\dg@X}{\dg@Y}{\dg@WEND}{\dg@HEND}%
         \divide\dg@WEND\tw@
         \ifnum\dg@XRPAD<\dg@WEND \dg@XRPAD=\dg@WEND \fi\fi
      \ifnum\dg@Y=\z@
         \dg@getnodesize
            {\dg@SLIST}{\dg@X}{\dg@Y}{\dg@WEND}{\dg@HEND}%
         \divide\dg@HEND\tw@
         \ifnum\dg@YTPAD<\dg@HEND \dg@YTPAD=\dg@HEND \fi\fi
      \ifnum\dg@Y=-\dg@VERT\relax
         \dg@getnodesize
            {\dg@SLIST}{\dg@X}{\dg@Y}{\dg@WEND}{\dg@HEND}%
         \divide\dg@HEND\tw@
         \ifnum\dg@YBPAD<\dg@HEND \dg@YBPAD=\dg@HEND \fi\fi
   \else\ifx\dg@PASS\dg@drawpass
      \dg@XNODE=\dg@X \multiply\dg@XNODE\@m
      \multiply\dg@XNODE\dg@XGRID
      \dg@YNODE=\dg@Y \multiply\dg@YNODE\@m
      \multiply\dg@YNODE\dg@YGRID
      \setbox\dg@NODEBOX=\hbox{$\dgeverynode{#2}$}%
      \put(\dg@XNODE,\dg@YNODE){\dg@makebox{\box\dg@NODEBOX}}%
   \fi\fi\fi}%

\newoptcommand{\dg@cr}{\@ne}[1]{%
   \ifx\dg@PASS\dg@savepass
      \dg@YTEMP=#1\relax
      \ifnum\dg@YTEMP<\@ne \dg@YTEMP=\@ne \fi
      \advance\dg@Y -\dg@YTEMP\relax
      \dg@X=-\@ne\relax\fi}%
\newoptcommand{\dg@arrow}{\@ne}[2]{%
   \begingroup
   \dg@USERSIZE=#1\relax
   \ifnum\dg@USERSIZE<\@ne \dg@USERSIZE=\@ne \fi
   \dg@parse{#2}%
   \ifx\dg@PASS\dg@savepass
      \ifx\dg@label\dgl@b \let\dg@label=\dgl@t \fi
      \ifx\dg@label\dgl@r \let\dg@label=\dgl@l \fi
      \let\dg@process=\dg@save
   \else\ifx\dg@PASS\dg@geompass
      \let\dg@process=\dg@ignore
      \dg@geomcalc
   \else\ifx\dg@PASS\dg@drawpass
      \let\dg@process=\dg@draw
      \dg@drawcalc
   \fi\fi\fi
   \dg@label{\dg@process{#1}{#2}}}%

\newoptcommand{\arrow}{\@ne}[2]{%
   \dg@parse{#2}%
   \ifx\dg@label\dgl@b \let\dg@label=\dgl@t \fi
   \ifx\dg@label\dgl@r \let\dg@label=\dgl@l \fi
   \dg@label{\dg@textarrow{#1}{#2}}}%

\def\dg@textarrow#1#2#3#4{%
   \mathop{{\dgHORIZPAD=0pt\relax\dgVERTPAD=0pt\relax
      \begin{diagram}
         \node{}\arrow[#1]{#2}{#3}{#4}\node{}
      \end{diagram}}}}

\def\dg@parse#1{%
   \let\dg@label=\dgl@ \dgo@
   \let\dg@type=\@empty \let\dg@lbltype=\@empty
   \@for\dg@list:=#1\do{%
      \ifx\dg@type\@empty \let\dg@type=\dg@list
      \else\ifx\dg@lbltype\@empty \let\dg@lbltype=\dg@list
         \@ifundefined{dgo@\dg@list}{}{\@nameuse{dgo@\dg@list}}%
      \else
         \@ifundefined{dgo@\dg@list}{}{\@nameuse{dgo@\dg@list}}%
      \fi\fi}%
   \@ifundefined{dgt@\dg@type}{\dgt@e}{\@nameuse{dgt@\dg@type}}%
   \@ifundefined{dgl@\dg@lbltype}{}{%
      \dg@letname\dg@label{dgl@\dg@lbltype}}}

\def\dg@draw#1#2#3#4{%
   \put(\dg@X,\dg@Y){\dg@makebox{%
      \begin{picture}(0,0)%
         \thinlines
         \put(\dg@XOFFSET,\dg@YOFFSET){%
            \dg@VECTOR(\dg@DX,\dg@DY){\dg@SIZE}}%
         \put(\dg@XLBL,\dg@YLBL){\dg@makebox{%
            \begin{picture}(0,0)%
               \put(\dg@XLBLOFF,\dg@YLBLOFF){%
                  \dg@makebox[\dg@LBLONE]{$\dgeverylabel{#3}$}}%
               \put(-\dg@XLBLOFF,-\dg@YLBLOFF){%
                  \dg@makebox[\dg@LBLTWO]{$\dgeverylabel{#4}$}}%
            \end{picture}}}%
      \end{picture}}}%
   \endgroup}%

\def\dg@save#1#2#3#4{%
   \endgroup 
   \toks\z@=\expandafter{\dg@ALIST}%
   \toks\tw@={\dg@arrow[#1]{#2}{#3}{#4}}%
   \edef\dg@ALIST{\the\toks\z@%
      \noexpand\dg@X=\number\dg@X\noexpand\relax
      \noexpand\dg@Y=\number\dg@Y\noexpand\relax
      \the\toks\tw@}}%

\def\dg@ignore#1#2#3#4{\endgroup}

\def\dg@geomcalc{%
   \dg@XEND=\dg@SIZE \multiply\dg@XEND\dg@USERSIZE
   \ifnum\dg@DX=\z@
      \dg@YEND=\dg@XEND \dg@XEND=\z@
      \dg@changesign\dg@YEND\dg@DY
   \else
      \dg@changesign\dg@XEND\dg@DX \dg@YEND=\dg@XEND
      \multiply\dg@YEND\dg@DY \divide\dg@YEND\dg@DX
   \fi
   \advance\dg@XEND\dg@X \advance\dg@YEND\dg@Y
   \dg@getnodesize
      {\dg@SLIST}{\dg@XEND}{\dg@YEND}{\dg@WEND}{\dg@HEND}%
   \dg@XOFFSET=\dg@WEND \dg@YOFFSET=\dg@HEND
   \dg@getnodesize
      {\dg@SLIST}{\dg@X}{\dg@Y}{\dg@WEND}{\dg@HEND}%
   \advance\dg@XOFFSET\dg@WEND \divide\dg@XOFFSET\tw@
   \advance\dg@YOFFSET\dg@HEND \divide\dg@YOFFSET\tw@
   \dg@XTEMP=\dgARROWLENGTH \dg@YTEMP=\dgARROWLENGTH
   \ifnum\dg@DX>\z@
      \dg@ZTEMP=\dg@DX \multiply\dg@XTEMP\dg@DX
   \else \dg@ZTEMP=-\dg@DX \multiply\dg@XTEMP -\dg@DX \fi
   \ifnum\dg@DY>\z@
      \advance\dg@ZTEMP\dg@DY \multiply\dg@YTEMP\dg@DY
   \else \advance\dg@ZTEMP -\dg@DY \multiply\dg@YTEMP -\dg@DY\fi
   \ifnum\dg@ZTEMP=\z@\else
      \divide\dg@XTEMP\dg@ZTEMP \divide\dg@YTEMP\dg@ZTEMP
      \advance\dg@XOFFSET\dg@XTEMP \advance\dg@YOFFSET\dg@YTEMP
   \fi
   \divide\dg@XOFFSET\dg@SIZE \divide\dg@XOFFSET\dg@USERSIZE
   \divide\dg@YOFFSET\dg@SIZE \divide\dg@YOFFSET\dg@USERSIZE
   \ifnum\dg@DX=\z@ \dg@XOFFSET=\z@ \fi
   \ifnum\dg@DY=\z@ \dg@YOFFSET=\z@ \fi
   \ifnum\dg@XGRID<\dg@XOFFSET \global\dg@XGRID=\dg@XOFFSET\fi
   \ifnum\dg@YGRID<\dg@YOFFSET \global\dg@YGRID=\dg@YOFFSET\fi
   \relax}

\def\dg@drawcalc{%
   \dg@XEND=\dg@SIZE \multiply\dg@XEND\dg@USERSIZE
   \ifnum\dg@DX=\z@
      \dg@YEND=\dg@XEND \dg@XEND=\z@
      \dg@changesign\dg@YEND\dg@DY
   \else
      \dg@changesign\dg@XEND\dg@DX \dg@YEND=\dg@XEND
      \multiply\dg@YEND\dg@DY \divide\dg@YEND\dg@DX
   \fi
   \advance\dg@XEND\dg@X \advance\dg@YEND\dg@Y
   \dg@getnodesize
      {\dg@SLIST}{\dg@XEND}{\dg@YEND}{\dg@WEND}{\dg@HEND}%
   \divide\dg@WEND\unitlength \divide\dg@HEND\unitlength
   \multiply\dg@DX\dg@XGRID \multiply\dg@DY\dg@YGRID
   \dg@rmcommondiv\tw@\dg@DX\dg@DY
   \dg@rmcommondiv\tw@\dg@DX\dg@DY 
   \dg@rmcommondiv\thr@@\dg@DX\dg@DY
   \multiply\dg@SIZE\dg@USERSIZE \multiply\dg@SIZE\@m
   \ifnum\dg@DX=\z@
      \multiply\dg@SIZE\dg@YGRID
      \divide\dg@HEND\tw@ \advance\dg@SIZE -\dg@HEND
      \dg@getnodesize
         {\dg@SLIST}{\dg@X}{\dg@Y}{\dg@WEND}{\dg@YOFFSET}%
      \divide\dg@YOFFSET\unitlength \divide\dg@YOFFSET\tw@
      \advance\dg@SIZE -\dg@YOFFSET
      \dg@XOFFSET=\z@
      \def\dg@LBLONE{r}\def\dg@LBLTWO{l}%
      \dg@XLBL=\z@ \dg@YLBL=\dg@SIZE
      \multiply\dg@YLBL\dg@LBLPOS
      \divide\dg@YLBL\dgARROWPARTS\relax
      \advance\dg@YLBL\dg@YOFFSET
      \dg@changesign\dg@YLBL\dg@DY
      \dg@changesign\dg@YOFFSET\dg@DY
   \else
      \multiply\dg@SIZE\dg@XGRID
      \ifnum\dg@DY=\z@
         \divide\dg@WEND\tw@ \advance\dg@SIZE -\dg@WEND
         \dg@getnodesize
            {\dg@SLIST}{\dg@X}{\dg@Y}{\dg@XOFFSET}{\dg@HEND}%
         \divide\dg@XOFFSET\unitlength \divide\dg@XOFFSET\tw@
         \advance\dg@SIZE -\dg@XOFFSET
         \dg@YOFFSET=\z@
         \def\dg@LBLONE{b}\def\dg@LBLTWO{t}%
         \dg@YLBL=\z@ \dg@XLBL=\dg@SIZE
         \multiply\dg@XLBL\dg@LBLPOS
         \divide\dg@XLBL\dgARROWPARTS\relax
         \advance\dg@XLBL\dg@XOFFSET
         \dg@changesign\dg@XLBL\dg@DX
         \dg@changesign\dg@XOFFSET\dg@DX
      \else
         \divide\dg@WEND\tw@ \divide\dg@HEND\tw@
         \multiply\dg@HEND\dg@DX \divide\dg@HEND\dg@DY
         \ifnum\dg@HEND<\z@ \multiply\dg@HEND\m@ne \fi
         \ifnum\dg@WEND<\dg@HEND \advance\dg@SIZE -\dg@WEND
         \else \advance\dg@SIZE -\dg@HEND \fi
         \dg@getnodesize
            {\dg@SLIST}{\dg@X}{\dg@Y}{\dg@WEND}{\dg@HEND}%
         \divide\dg@WEND\unitlength \divide\dg@WEND\tw@
         \divide\dg@HEND\unitlength \divide\dg@HEND\tw@
         \multiply\dg@HEND\dg@DX \divide\dg@HEND\dg@DY
         \ifnum\dg@HEND<\z@ \multiply\dg@HEND\m@ne \fi
         \ifnum\dg@WEND<\dg@HEND \dg@XOFFSET=\dg@WEND
         \else \dg@XOFFSET=\dg@HEND \fi
         \advance\dg@SIZE -\dg@XOFFSET
         \dg@changesign\dg@XOFFSET\dg@DX
         \dg@YOFFSET=\dg@XOFFSET
         \multiply\dg@YOFFSET\dg@DY \divide\dg@YOFFSET\dg@DX
         \def\dg@LBLONE{br}\def\dg@LBLTWO{tl}%
         \ifnum\dg@DX<\z@ \ifnum\dg@DY>\z@
            \def\dg@LBLONE{bl}\def\dg@LBLTWO{tr}\fi\fi
         \ifnum\dg@DX>\z@ \ifnum\dg@DY<\z@
            \def\dg@LBLONE{bl}\def\dg@LBLTWO{tr}\fi\fi
         \dg@XLBL=\dg@SIZE
         \multiply\dg@XLBL\dg@LBLPOS
         \divide\dg@XLBL\dgARROWPARTS\relax
         \dg@changesign\dg@XLBL\dg@DX
         \dg@YLBL=\dg@XLBL
         \multiply\dg@YLBL\dg@DY \divide\dg@YLBL\dg@DX
         \advance\dg@XLBL\dg@XOFFSET
         \advance\dg@YLBL\dg@YOFFSET
      \fi
   \fi
   \dg@XLBLOFF=-\dg@DY \dg@changesign\dg@XLBLOFF\dg@DX
   \dg@YLBLOFF=\dg@DX \dg@changesign\dg@YLBLOFF\dg@DX
   \ifnum\dg@DX=\z@ \dg@XLBLOFF=\m@ne \fi
   \dg@XTEMP=\dg@DX \dg@changesign\dg@XTEMP\dg@DX
   \dg@YTEMP=\dg@DY \dg@changesign\dg@YTEMP\dg@DY
   \ifnum\dg@YTEMP>\dg@XTEMP \dg@XTEMP=\dg@YTEMP \fi
   \ifnum\dg@XTEMP=\z@ \dg@XTEMP=\@ne \fi
   \multiply\dg@XLBLOFF\dg@LBLOFF \divide\dg@XLBLOFF\dg@XTEMP
   \multiply\dg@YLBLOFF\dg@LBLOFF \divide\dg@YLBLOFF\dg@XTEMP
   \multiply\dg@X\@m \multiply\dg@X\dg@XGRID
   \multiply\dg@Y\@m \multiply\dg@Y\dg@YGRID
   \relax}%

\def\dg@rmcommondiv#1#2#3{%
   \dg@XTEMP=#2\relax
   \divide\dg@XTEMP #1\relax \multiply\dg@XTEMP #1\relax
   \dg@YTEMP=#3\relax
   \divide\dg@YTEMP #1\relax \multiply\dg@YTEMP #1\relax
   \ifnum\dg@XTEMP=#2\relax \ifnum\dg@YTEMP=#3\relax
      \divide#2#1\relax \divide#3#1\relax \fi\fi}%

\def\dg@changesign#1#2{%
   \ifnum #2<\z@ \multiply#1\m@ne
   \else\ifnum #2=\z@ #1=\z@ \fi\fi}%

\def\dg@getnodesize#1#2#3#4#5{%
   #4=\z@\relax #5=\z@\relax
   \expandafter\@for\expandafter\dg@trynode
   \expandafter:\expandafter=#1\do{%
      \dg@XNODE=\m@ne 
      \dg@trynode
      \ifnum #2=\dg@XNODE \ifnum #3=\dg@YNODE
         #4=\dg@XTEMP\relax #5=\dg@YTEMP\relax\fi\fi}}%

\newoptcommand{\dg@makebox}{}[2]{%
   \expandafter\makebox\expandafter(\expandafter
      0\expandafter,\expandafter0\expandafter)\expandafter
      [#1]{#2}}%

\def\dg@novector(#1,#2)#3{}%

\def\dg@letname#1#2{%
   \relax\expandafter
   \let\expandafter #1\csname #2\endcsname\relax}%

\def\dgl@#1{#1{}{}}%
\def\dgl@t#1#2{#1{#2}{}}%
\def\dgl@b#1#2{#1{}{#2}}%
\def\dgl@tb#1#2#3{#1{#2}{#3}}%
\def\dgl@l#1#2{#1{#2}{}}%
\def\dgl@r#1#2{#1{}{#2}}%
\def\dgl@lr#1#2#3{#1{#2}{#3}}%
\makeatother

\newcommand{\su  }{{\cal S}{\cal U}(r,L_0) }
\newcommand{\urd }{{\cal U}(r,d) }
\newcommand{\sustable  }{ {\cal S}{\cal U}^{s}(r,L_0) }
\newcommand{\cotangent  }{ {\cal X}(r,L_0) }
\newcommand{\wbase }{W^{\rm reg}_{r-1,r} }
\newcommand{\breg}{B^{\rm reg}}
\newcommand{\bregp}{B_p^{\rm reg}}
\newcommand{\phibase} {\phi _{r-1,r}}
\newcommand{\unspec }{\tilde {\cal C}  }
\newcommand{\spec  }{\tilde{C}_s  }
\newcommand{\unprymd  }{{\rm Prym}_{\tilde {d}}(\tilde {\cal C}, C) }
\newcommand{\prymd}{{\rm Prym}_{\tilde{d}}(\tilde{C}_s ,  C)}
\newcommand{\unprym }{{\rm Prym}(\tilde {\cal C}, C)  }
\newcommand{\prym}{{\rm Prym}( \tilde {C}_s, C)  }
\newcommand{\unjacd  }{J^{\tilde {d}} (\tilde {\cal C})  }
\newcommand{\unjac }{ J ^0 (\tilde { \cal C} ) }
\newcommand{\fibprod  }{\tilde {\cal C \;} _{ h}  \!\! \times _{\phi
\, h}
\tilde {\cal C} }
\newcommand{\higgs  }{{\cal M}(r, L_0)  }
\newcommand{\phitilded}{\tilde{\Phi} _{\tilde {d}} }

\newcommand{\phitilde}{\tilde{\Phi}}

\newcommand{\psid}{\psi _{\tilde {d}}^* }

\newcommand{\so}{S^{\circ}}
\newcommand{\pic}{{\rm Pic} \,}
\newcommand{\pr}{{\Bbb P} ^1}
\newcommand{\poincd}{{\cal P}^{\tilde d}}
\newcommand{\poinc}{{\cal P}}
\newcommand{\aff}{{\Bbb A}{\rm ff}}
\newcommand{\lra}{\longrightarrow}

\newtheorem{lem}{Lemma}[section]
\newtheorem{cor}{Corollary}[section]
\newtheorem{theor}{Theorem}
\newtheorem{prop}{Proposition}[section]
\newtheorem{defi}{Definition}[section]
\newtheorem{rem}{Remark}[section]

\begin{document}
 \addtolength{\baselineskip}{4pt}
\setlength{\unitlength}{1mm}
%
\expandafter\ifx\csname amssym.def\endcsname\relax \else\endinput\fi
%
\expandafter\edef\csname amssym.def\endcsname{%
       \catcode`\noexpand\@=\the\catcode`\@\space}
\catcode`\@=11
%

\def\undefine#1{\let#1\undefined}
\def\newsymbol#1#2#3#4#5{\let\next@\relax
 \ifnum#2=\@ne\let\next@\msafam@\else
 \ifnum#2=\tw@\let\next@\msbfam@\fi\fi
 \mathchardef#1="#3\next@#4#5}
\def\mathhexbox@#1#2#3{\relax
 \ifmmode\mathpalette{}{\m@th\mathchar"#1#2#3}%
 \else\leavevmode\hbox{$\m@th\mathchar"#1#2#3$}\fi}
\def\hexnumber@#1{\ifcase#1 0\or 1\or 2\or 3\or 4\or 5\or 6\or 7\or
8\or
 9\or A\or B\or C\or D\or E\or F\fi}

\font\tenmsa=msam10
\font\sevenmsa=msam7
\font\fivemsa=msam5
\newfam\msafam
\textfont\msafam=\tenmsa
\scriptfont\msafam=\sevenmsa
\scriptscriptfont\msafam=\fivemsa
\edef\msafam@{\hexnumber@\msafam}
\mathchardef\dabar@"0\msafam@39
\def\dashrightarrow{\mathrel{\dabar@\dabar@\mathchar"0\msafam@4B}}
\def\dashleftarrow{\mathrel{\mathchar"0\msafam@4C\dabar@\dabar@}}
\let\dasharrow\dashrightarrow
\def\ulcorner{\delimiter"4\msafam@70\msafam@70 }
\def\urcorner{\delimiter"5\msafam@71\msafam@71 }
\def\llcorner{\delimiter"4\msafam@78\msafam@78 }
\def\lrcorner{\delimiter"5\msafam@79\msafam@79 }
\def\yen{{\mathhexbox@\msafam@55 }}
\def\checkmark{{\mathhexbox@\msafam@58 }}
\def\circledR{{\mathhexbox@\msafam@72 }}
\def\maltese{{\mathhexbox@\msafam@7A }}

\font\tenmsb=msbm10
\font\sevenmsb=msbm7
\font\fivemsb=msbm5
\newfam\msbfam
\textfont\msbfam=\tenmsb
\scriptfont\msbfam=\sevenmsb
\scriptscriptfont\msbfam=\fivemsb
\edef\msbfam@{\hexnumber@\msbfam}
\def\Bbb#1{{\fam\msbfam\relax#1}}
\def\widehat#1{\setbox\z@\hbox{$\m@th#1$}%
 \ifdim\wd\z@>\tw@ em\mathaccent"0\msbfam@5B{#1}%
 \else\mathaccent"0362{#1}\fi}
\def\widetilde#1{\setbox\z@\hbox{$\m@th#1$}%
 \ifdim\wd\z@>\tw@ em\mathaccent"0\msbfam@5D{#1}%
 \else\mathaccent"0365{#1}\fi}
\font\teneufm=eufm10
\font\seveneufm=eufm7
\font\fiveeufm=eufm5
\newfam\eufmfam
\textfont\eufmfam=\teneufm
\scriptfont\eufmfam=\seveneufm
\scriptscriptfont\eufmfam=\fiveeufm
\def\frak#1{{\fam\eufmfam\relax#1}}
\let\goth\frak

\csname amssym.def\endcsname
\newsymbol\rtimes 226F
\newsymbol\nmid 232D
\newsymbol\varnothing 203F
\newread\epsffilein    
\newif\ifepsffileok    
\newif\ifepsfbbfound   
\newif\ifepsfverbose   
\newdimen\epsfxsize    
\newdimen\epsfysize    
\newdimen\epsftsize    
\newdimen\epsfrsize    
\newdimen\epsftmp      
\newdimen\pspoints     
\pspoints=1bp          
\epsfxsize=0pt         
\epsfysize=0pt         
\def\epsfbox#1{\global\def\epsfllx{72}\global\def\epsflly{72}%
   \global\def\epsfurx{540}\global\def\epsfury{720}%
   \def\lbracket{[}\def\testit{#1}\ifx\testit\lbracket
   \let\next=\epsfgetlitbb\else\let\next=\epsfnormal\fi\next{#1}}%
\def\epsfgetlitbb#1#2 #3 #4 #5]#6{\epsfgrab #2 #3 #4 #5 .\\%
   \epsfsetgraph{#6}}%
\def\epsfnormal#1{\epsfgetbb{#1}\epsfsetgraph{#1}}%
\def\epsfgetbb#1{%
%
%
\openin\epsffilein=#1
\ifeof\epsffilein\errmessage{I couldn't open #1, will ignore it}\else
%
doesn't
%
   {\epsffileoktrue \chardef\other=12
    \def\do##1{\catcode`##1=\other}\dospecials \catcode`\ =10
    \loop
       \read\epsffilein to \epsffileline
       \ifeof\epsffilein\epsffileokfalse\else
%
%
          \expandafter\epsfaux\epsffileline:. \\%
       \fi
   \ifepsffileok\repeat
   \ifepsfbbfound\else
    \ifepsfverbose\message{No bounding box comment in #1; using
defaults}\fi\fi
   }\closein\epsffilein\fi}%
%
%
\def\epsfclipstring{}
\def\epsfclipon{\def\epsfclipstring{ clip}}%
\def\epsfclipoff{\def\epsfclipstring{}}%
\def\epsfsetgraph#1{%
   \epsfrsize=\epsfury\pspoints
   \advance\epsfrsize by-\epsflly\pspoints
   \epsftsize=\epsfurx\pspoints
   \advance\epsftsize by-\epsfllx\pspoints
%
picture.
%
   \epsfxsize\epsfsize\epsftsize\epsfrsize
   \ifnum\epsfxsize=0 \ifnum\epsfysize=0
      \epsfxsize=\epsftsize \epsfysize=\epsfrsize
      \epsfrsize=0pt
%
arithmetic!
reasonably
%
     \else\epsftmp=\epsftsize \divide\epsftmp\epsfrsize
       \epsfxsize=\epsfysize \multiply\epsfxsize\epsftmp
       \multiply\epsftmp\epsfrsize \advance\epsftsize-\epsftmp
       \epsftmp=\epsfysize
       \loop \advance\epsftsize\epsftsize \divide\epsftmp 2
       \ifnum\epsftmp>0
          \ifnum\epsftsize<\epsfrsize\else
             \advance\epsftsize-\epsfrsize \advance\epsfxsize\epsftmp
\fi
       \repeat
       \epsfrsize=0pt
     \fi
   \else \ifnum\epsfysize=0
     \epsftmp=\epsfrsize \divide\epsftmp\epsftsize
     \epsfysize=\epsfxsize \multiply\epsfysize\epsftmp
     \multiply\epsftmp\epsftsize \advance\epsfrsize-\epsftmp
     \epsftmp=\epsfxsize
     \loop \advance\epsfrsize\epsfrsize \divide\epsftmp 2
     \ifnum\epsftmp>0
        \ifnum\epsfrsize<\epsftsize\else
           \advance\epsfrsize-\epsftsize \advance\epsfysize\epsftmp
\fi
     \repeat
     \epsfrsize=0pt
    \else
     \epsfrsize=\epsfysize
    \fi
   \fi
%
parse.
   \ifepsfverbose\message{#1: width=\the\epsfxsize,
height=\the\epsfysize}\fi
   \epsftmp=10\epsfxsize \divide\epsftmp\pspoints
   \vbox to\epsfysize{\vfil\hbox to\epsfxsize{%
      \ifnum\epsfrsize=0\relax
        \includegraphics{#1}%
      \else
        \epsfrsize=10\epsfysize \divide\epsfrsize\pspoints
        \includegraphics{#1}%
      \fi
      \hfil}}%
\global\epsfxsize=0pt\global\epsfysize=0pt}%
%
%
{\catcode`\%=12
\global\let\epsfpercent=
%
%
\long\def\epsfaux#1#2:#3\\{\ifx#1\epsfpercent
   \def\testit{#2}\ifx\testit\epsfbblit
      \epsfgrab #3 . . . \\%
      \epsffileokfalse
      \global\epsfbbfoundtrue
   \fi\else\ifx#1\par\else\epsffileokfalse\fi\fi}%
%
%
\def\epsfempty{}%
\def\epsfgrab #1 #2 #3 #4 #5\\{%
\global\def\epsfllx{#1}\ifx\epsfllx\epsfempty
      \epsfgrab #2 #3 #4 #5 .\\\else
   \global\def\epsflly{#2}%
   \global\def\epsfurx{#3}\global\def\epsfury{#4}\fi}%
%
%
\def\epsfsize#1#2{\epsfxsize}
%
%
\let\epsffile=\epsfbox

 \title {The automorphism group of the moduli space  \\ of semi
stable vector
bundles}
\author{ Alexis Kouvidakis\footnotemark[1]
$\;\;\;\;\;\;\;\;\;\;\;\;\;\;\;\;\;\;\;\;\;\;\;\;\;$
  Tony Pantev   \\  $ \;    $   \\
University of Pennsylvania \\
 Department of Mathematics, DRL \\
 Philadelphia, PA 19104-6395\\
 $\;\;$  alex@math.upenn.edu  $\;\;\;\;\;\;\;\;\;\;\;\;\;\;$
pantev@math.upenn.edu  }
\date{ }
\maketitle
\footnotetext{New address: Department of Mathematics, University of
Crete,
Iraklion 71409, Greece.}
\setcounter{section}{-1}
\section{Introduction}
To every smooth curve $C$ one can associate a natural variety $\urd$
- the
moduli space of semi stable vector bundles of rank $r$ and degree $d$
on $C$.
This canonical object has a rich geometrical structure which reflects
in a
beautiful way the geometry of the curve $C$ and its deformations.
In the
abelian case $r=1$ this relationship is classical and goes back to
the Riemann
inversion problem, the theory of  Jacobi theta-functions and the
Torelli
theorem. Understanding the  subtleties of the non-abelian case has
resulted
over the years in a multitude of new ideas and unexpected connections
with
other branches of
mathematics.  One of the remarkable similarities between the Jacobian
varieties
and the higher rank moduli spaces is the existence of theta line
bundles on
them, see \cite{dn}. These are ample determinantal line bundles,
naturally
arising from the interpretation of $\urd$ as a moduli space and
provide
valuable information about the projective geometry of $\urd$.

In this paper we describe the groups of automorphisms  and of
polarized
automorphisms of $\urd$ for $r > 1$. To understand the problem
better, consider
first the
case $r = 1$. It is well known that for any integer $d$, the group of
automorphisms of the Jacobian
 $J^{d}(C) = {\cal U}(1,d)$ is isomorphic to the semi-direct product
 \[ J^{0}(C) \rtimes {\rm Aut}^{\rm group}(J^{0}(C)) \cong {\rm
Aut}(J^{d}(C)),
\]
 where  ${\rm Aut}^{\rm group}(J^{0}(C))$ is the group of group
automorphisms
of $J^{0}$. The above isomorphism depends on the choice of a point
 $L_{0} \in
 J^{d}(C)$ and is given explicitly by $(\xi,\phi) \rightarrow
T_{\xi}\circ
T_{L_{0}}\circ \phi \circ T_{L_{0}}^{-1}$. Moreover, the subgroup of
automorphisms preserving the class of the theta bundle is just
 \[ {\rm Aut}^{\theta}(J^{d}(C)) \cong \left\{
 \begin{array}{ll}
 J^{0}(C) \rtimes ((\pm {\rm id})\times {\rm Aut}(C)) & C {\rm -not
\;
hyperelliptic,} \\
  J^{0}(C) \rtimes {\rm Aut}(C) & C{\rm -hyperelliptic}.
  \end{array} \right.
  \]

Some natural automorphisms of $\urd$ can be obtained by analogy with
this case.
For example, the Jacobian $J^{0}(C)$ acts by translations on the
moduli space
and, as
it turns out, it actually coincides with the identity component of
the
automorphism group. Furthermore, we have a natural morphism
\begin{equation}
 \det : \urd \longrightarrow J^{d}(C) \label{eq01}
\end{equation}
sending a bundle $E$ to its determinant $\det (E)$. The fibers of
this morphism
are moduli spaces on its own right, e.g. the fiber $\det^{-1}(L_{0})$
over a
point $L_{0}$ is just the moduli space $\su$ of semistable bundles of
rank $r$
and fixed determinant $L_{0}$.  One can use the fibration
(\ref{eq01}) to
construct some natural automorphisms of $\urd$, that is,  certain
automorphisms
of the fiber $\su$ can be glued with certain automorphisms of the
base
$J^{d}(C)$ to produce global automorphisms of  the moduli space.
After choosing
 a point $L_{0} \in
J^{d}(C)$ we can trivialize the bundle (\ref{eq01}) by pulling it
back to an
\'{e}tale cover
of $J^{d}(C)$:
\begin{equation}
{\divide\dgARROWLENGTH by 9
\begin{diagram}
\node{\su\times J^{0}(C)} \arrow[3]{e,t}{(E,L) \rightarrow E\otimes
L}
\arrow{s,l}{p_{2}} \node[3]{\urd} \arrow{s,r}{\det} \\
\node{J^{0}(C)} \arrow[3]{e,b}{L \rightarrow L^{\otimes r}\otimes
L_{0}}
\node[3]{J^{d}(C)}
\end{diagram}} \label{eq02}
\end{equation}
The top and the bottom rows of the diagram ({\ref{eq02}) are Galois
covers with
Galois group the group of $r$-torsion points $J^{0}(C)[r]$. Consider
the
subgroup ${\rm Aut}_{[r]}(J^{0}(C)) \subset {\rm Aut}(J^{0}(C))$
consisting of
group automorphisms of $J^{0}(C)$ which act trivially on the torsion
points
$J^{0}(C)[r]$.
Then given a translation, $T_{\mu}$, by a torsion point $\mu \in
J^{0}(C)[r]$
and
$\varphi \in
{\rm Aut}_{[r]}(J^{0}(C))$, the automorphism $T_{\mu} \times \varphi$
of
$\su\times J^{0}(C)$ commutes with the action of the Galois group and
hence
descends
to an automorphism of $\urd$.  In addition, every symmetry of the
curve $C$
will induce
an automorphism of $\urd$.

In the special case $r | 2d$, by choosing $\nu \in
J^{\frac{2d}{r}}(C)$ with
property
$\nu^{\otimes r} = L_{0}^{\otimes 2}$, we can construct an additional
automorphism
$\delta$ of order two:
\[
\begin{array}{cccc}
\delta : & \urd & \longrightarrow & \urd \\
&  E & \longrightarrow & E^{\vee}\otimes \nu
\end{array}
\]

To combine these, consider the subgroups
${\frak U}{\frak G}^{\circ}_{[r]} \subseteq {\frak U}{\frak G}_{[r]}
\subseteq
{\rm Aut}(J^{0}(C))\times
{\rm Aut}(C)$ defined by
\begin{defi} \label{def01}
\[ {\frak U}{\frak G}^{\circ}_{[r]} = \{ (\phi,\sigma) \; | \;
\phi^{-1}\circ
\sigma^{*} \in  {\rm Aut}_{[r]}(J^{0}(C))\},\]
and
\[ {\frak U}{\frak G}_{[r]} = \{ (\phi,\sigma) \; | \; \phi^{-1}\circ
(\pm {\bf
1})\circ
\sigma^{*} \in  {\rm Aut}_{[r]}(J^{0}(C)) \}. \]
\end{defi}
The group ${\frak U}{\frak G}^{\circ}_{[r]}$ is of index two in
${\frak
U}{\frak G}_{[r]}$ and we
have a natural homomorphism:
\begin{equation}
\begin{array}{ccl}
J^{0}(C) \rtimes  {\frak U}{\frak G}^{\circ}_{[r]} & \longrightarrow
& {\rm
Aut}(\urd ) \\
(\xi,\phi,\sigma) & \longrightarrow & (E \mapsto \sigma ^* E \otimes
\xi
\otimes \phi (\eta ) \otimes  \sigma ^*(\eta ^{-1} ))
\end{array} \label{eq03}
\end{equation}
where $\eta \in J^{0}(C)$ is an arbitrary line bundle with the
property $\eta
^{\otimes r} = \det E \otimes  L_0 ^{-1}$.  If in addition $r \, | \,
2d$,
then the homomorphism (\ref{eq03}) can
be extended to
\begin{equation} \label{eq04}
\begin{array}{ccl}
J^{0}(C) \rtimes  {\frak U}{\frak G}_{[r]} & \longrightarrow & {\rm
Aut}(\urd )
\\
(\xi,\phi,\sigma) & \longrightarrow & E \mapsto \left\{
\begin{array}{ll}
\sigma ^* E \otimes \xi \otimes \phi (\eta ) \otimes  \sigma ^*(\eta
^{-1} ) &
{\rm if} \;
\phi^{-1}\circ \sigma^{*} \in {\rm Aut}_{[r]}(J^{0}(C))  \\
 \sigma ^* E^{\vee} \otimes \xi \otimes \phi (\eta ) \otimes  \sigma
^*(\eta) &
{\rm if} \;
\phi^{-1}\circ ({\bf - 1}) \circ \sigma^{*} \in {\rm
Aut}_{[r]}(J^{0}(C))   \\
\end{array} \right.
\end{array}
\end{equation}

The commutative diagram (\ref{eq02}) suggests that the subgroups
(\ref{eq03})
(or (\ref{eq04}) in the case $r \, | \, 2d$) will not differ too much
from the
full automorphism group, if the moduli space $\su $ doesn't have
excess
automorphisms. The variety $\su $ has a lot of advantages. It is a
Fano variety
of  \linebreak Picard number one and by a theorem of Narasimhan and
Ramanan,
\cite{nr},  $H^{0}(\su ,T\su) = 0$ unless $g=2$, $r=2$ and
$d$ even. Therefore, in general, ${\rm Aut}(\su )$ is finite and
contains the
group $J^{0}(C)[r]$. In the case $r = 2$, odd degree and $g = 2$
the group ${\rm Aut}({\cal S}{\cal U}^{s}(2,L_0) )$ was described by
Newstead
\cite{ne}
as a consequence of his proof of the Torelli theorem for the variety
${\cal
S}{\cal U}^{s}(2,L_0)$. To generalize his result we adopt a different
viewpoint
- we use Hitchin's abelianization to prove our main theorem

\begin{theor} Let $C$ be a curve without automorphisms. Then the
automorphism
group of the moduli space $\su$ can be described as follows
\begin{enumerate}
\item If $r \nmid 2d$, then the natural map
\[ \begin{array}{ccl}
 J^{0}(C)[r] & \longrightarrow & {\rm Aut}(\su ) \\
 \mu & \longrightarrow & (E \mapsto E\otimes \mu)
 \end{array}
 \]
 is an isomorphism, and
\item  If $r \, | \, 2d$, then the natural map
\[
\begin{array}{ccl}
J^{0}(C)[r] \rtimes {\Bbb Z}/2{\Bbb Z} & \longrightarrow & {\rm
Aut}(\su ) \\
(\mu, \varepsilon) & \longrightarrow & (E \mapsto
\delta^{\varepsilon}(E)\otimes \mu)
\end{array}
\]
is isomorphism for $r \geq 3$ and has kernel ${\Bbb Z}/2{\Bbb Z}$ for
$r = 2$.
\end{enumerate}
\label{theor1}
\end{theor}

 \noindent
The strategy of the proof is to lift the automorphism $\Phi$ to a
birational
automorphism of  the moduli space of Higgs bundles and then study the
induced
automorphism $\phitilded$ of the family $\unprymd$ of smooth spectral
Pryms. To
determine the map $\phitilded$ explicitly we study the ring of
relative
correspondences over the universal spectral curve and show that the
family of
groups $\unprym$  does not have non-trivial
 global group automorphisms. Furthermore, we show that all the
sections of
$\unprymd$
 come from pull-back of line bundles on the base curve $C$ and
conclude that
 $\phitilded$  must be a multiplication by $\pm 1$ along the fibers
of
$\unprymd$
 followed by a translation by a pull-back bundle. To finish the
proof, we use
the
 rational map from a spectral Prym to the moduli space $\su$ to
recover from
 $\phitilded$ the original automorphism $\Phi$.

For curves with automorphisms the statement of the main theorem has
to be
modified appropriately.  To avoid some technical complications which
occur in
the case $g = 2$, we
assume that $g \geq 3$ throughout the paper. Consider the subgroups
${\frak  S}{\frak G}^{\circ}_{[r]} \subseteq {\frak  S}{\frak
G}_{[r]}
\subseteq
  J^{0}(C)\times {\rm Aut}(C)$
defined as
\begin{defi}\label{def02}
\[ {\frak  S}{\frak G}^{\circ}_{[r]} = \{ (\xi, \sigma) \; | \;
\xi^{r} =
L_{0}\otimes\sigma^{*}(L_{0}^{-1}) \}, \]
and
\[ {\frak  S}{\frak G}_{[r]} = \{ (\xi, \sigma) \; | \; \xi^{r} =
L_{0}\otimes\sigma^{*}(L_{0}^{\pm 1}) \}. \]
\end{defi}
Again ${\frak  S}{\frak G}^{\circ}_{[r]}$ is a normal subgroup of
index two in
${\frak  S}{\frak G}_{[r]}$
and we have

\begin{theor} Let $C$ be any smooth curve of genus $g \geq 3$. We
have
\begin{enumerate}
\item If $r \nmid 2d$, then the natural map
\[ \begin{array}{ccl}
 {\frak  S}{\frak G}^{\circ}_{[r]} & \longrightarrow & {\rm Aut}(\su
) \\
(\xi , \sigma) & \longrightarrow & (E \mapsto \sigma^{*}E\otimes \xi)
 \end{array}
 \]
 is surjective, and
\item  If $r \, | \, 2d$, then the natural map
\[
\begin{array}{ccl}
{\frak  S}{\frak G}_{[r]}& \longrightarrow & {\rm Aut}(\su ) \\
(\xi, \sigma) & \longrightarrow & E \mapsto \left\{
\begin{array}{ll}
\sigma^{*}E\otimes \xi & {\rm if} \; (\xi, \sigma) \in {\frak
S}{\frak
G}^{\circ}_{[r]} \\
\sigma^{*}E^{\vee}\otimes \xi & {\rm if} \; (\xi, \sigma) \in {\frak
S}{\frak
G}_{[r]} \setminus  {\frak  S}{\frak G}^{\circ}_{[r]}
\end{array} \right.
\end{array}
\]
is surjective.
\end{enumerate}
\label{theor2}
\end{theor}

\begin{rem}
{\rm  The maps in the above Theorem \ref{theor2} are actually
isomorphisms.
Our proof of the injectivity involves some standard arguments from
Prym theory
similar to those in \cite{nr2}; the proof is not included here
because of its
length. }
\end{rem}
\noindent
The first part of Theorem \ref{theor2} gives a positive answer to a
question
posed by
Tyurin in the survey paper \cite{ty}.

The main theorem, together with the observation that every
automorphism of
$\urd $
must lift to an automorphism of $\su \times J^{0}(C)$, see Lemma
\ref{lemurd3},
gives

 \begin{theor} Let $C$ be a smooth curve of genus $g \geq 3$. Then
the
automorphisms
 of the full moduli space $\urd $ can be described as follows
 \begin{enumerate}
 \item If $r \nmid 2d$, then the map {\rm (}\ref{eq03}{\rm )}
 \[ \begin{array}{ccl}
J^{0}(C) \rtimes  {\frak U}{\frak G}^{\circ}_{[r]} & \hookrightarrow
& {\rm
Aut}(\urd ) \\
(\xi,\phi,\sigma) & \rightarrow & (E \mapsto \sigma ^* E \otimes \xi
\otimes
\phi (\eta ) \otimes  \sigma ^*(\eta ^{-1} ))
\end{array} \]
is surjective, and
\item If $r | 2d$, then the map {\rm (}\ref{eq04}{\rm )}
\[\begin{array}{ccl}
J^{0}(C) \rtimes  {\frak U}{\frak G}_{[r]} & \hookrightarrow & {\rm
Aut}(\urd )
\\
(\xi,\phi,\sigma) & \rightarrow & E \mapsto \left\{
\begin{array}{ll}
\sigma ^* E \otimes \xi \otimes \phi (\eta ) \otimes  \sigma ^*(\eta
^{-1} ) &
{\rm if} \;
\phi^{-1}\circ \sigma^{*} \in {\rm Aut}_{[r]}(J^{0}(C))\\
 \sigma ^* E^{\vee} \otimes \xi \otimes \phi (\eta ) \otimes  \sigma
^*(\eta) &
{\rm if} \;
\phi\circ \sigma^{*} \in {\rm Aut}_{[r]}(J^{0}(C)) \\
\end{array} \right.
\end{array}
\]
is surjective.
\end{enumerate}
\label{theor3}
\end{theor}

The description of the polarized automorphisms of $\urd$ follows
easily from
Theorem \ref{theor3}:

\begin{theor} For any curve $C$ of genus $g \geq 3$ the
automorphisms of the
moduli space $\urd$ which preserve the class of the theta bundle are
those
belonging to the image of the subgroup \linebreak $J^{0}(C) \rtimes
((\pm {\rm
id})\times {\rm Aut}(C)) \subseteq  J^{0}(C) \rtimes  {\frak U}{\frak
G}_{[r]}$
under the map {\rm (}\ref{eq04}{\rm ).}
\label{theor4}
\end{theor}

As an application of the technics we use in this work, we prove in
the Appendix
at the end of the paper the Torelli Theorem for the moduli space of
vector
bundles:
\begin{theor}
Let $C_1$,  $C_2$ be  smooth curves of genus $g \geq 3$ and $L_1$,
$L_2$  line
bundles of degree $d$ on $C_1$,  $C_2$ respectively. If ${\cal
S}{\cal
U}_{C_1}(r, L_1) \simeq
 {\cal S}{\cal U}_{C_2}(r, L_2)$,  then  $C_1 \simeq C_2$.
\label{Torelli}
\end{theor}

The paper is organized as follows. In the first section we gather all
the facts
about
the moduli space of Higgs bundles, the linear system of spectral
curves and the
Hitchin map which we are going to use later on. We have included the
proofs of
some statements which seem to be well known to the experts in the
field but
which can not be found in the standard sources, as well as, the
proofs of some
facts about the discriminant locus in the Hitchin base which help us
simplify
our main argument. In the second section we construct the induced
automorphism
$\phitilded$ and study its first properties. The third section
contains a
discussion of the relative Picard of the universal spectral curve and
its
sections over a Zariski open set in the Hitchin base. In section four
we
analyze the ring of relative correspondences on the fibers of the
universal
spectral curve whose description is the key ingredient in the proof
of the main
theorem. In the fifth section we give the proof of the main theorem
and discuss
the case of curves with symmetries.  In the sixth section we derive
the
description of ${\rm Aut}(\urd )$.

\bigskip

\noindent
{\bf Acknowledgments:} We are very grateful to Ron Donagi for his
constant
support and for many valuable discussions and suggestions. We would
like to
thank
 Eyal Markman for explaining to us the proof of Lemma \ref{lem14} and
George
Pappas for helpful conversations.

 \bigskip

\noindent
{\bf Notation and Conventions}

\medskip

\noindent
 $\su \, :$   The moduli space of semi stable vector bundles of rank
$r$ and
fixed determinant $L_0$ of degree $d$. \\
${\cal S} {\cal U} ^s(r, L_0) \, :$   The moduli  space of  stable
vector
bundles of rank $r$ and fixed determinant $L_0$ of degree $d$. \\
${\cal X}(r, L_0) \, : $  The total space of the cotangent bundle
$T^*
\sustable$. \\
${\cal M}(r, L_0) \,  : $  The moduli  space of  semi stable  Higgs
pairs. \\
$H : {\cal M}(r, L_0) \lra W \, : $ The Hitchin map. \\
$ C \, : $  A smooth curve of genus $g \geq 3$.  \\
$  \omega _C  \, : $  The canonical bundle on $C$. \\
$ \so   \, : $  The total space of the canonical bundle $\omega _C$.
\\
$ S  \, : $   The space ${\Bbb P}({\cal  O} \oplus \omega _C)$. \\
$ \alpha : S  \lra C  \, : $   The canonical map.  \\
$ Y  \in H^0(S, {\cal O}_S(1) ) \, : $  The infinity section of $S$.
\\
$  X  \in H^0(S, {\cal O}_S(1)  \oplus \alpha ^*  \omega _C)  \, : $
The zero
section of $S$. \\
$  x= X/Y  \, : $  The tautological section of $\so $. \\
$  Y_{\infty } \, : $ The divisor at infinity  ${\rm div}(Y) \subset
S  $.  \\
$  X_0   \, : $ The zero divisor ${\rm div}(X) \subset \so $.  \\
$ \overline{W} = H^0(C, {\cal O}) \oplus H^0(C, \omega _C^{\otimes
2})  \oplus
\cdots \oplus  H^0(C, \omega _C^{\otimes r})  $.  \\
$  W_k= H^0(C, \omega _C^{\otimes k})  $.  \\
$  W \simeq W_2  \oplus \cdots  \oplus W_r$, embedded in
$\overline{W}$ as $\{
1  \} \oplus W_2 \oplus \cdots  \oplus W_r$. \\
$ W^{\infty }   \simeq W_2  \oplus \cdots  \oplus W_r$, embedded in
$\overline{W}$
 as $\{ 0  \} \oplus W_2 \oplus \cdots  \oplus W_r$. \\
$  {\cal D} \, : $ The discriminant divisor. \\
$ W^{{\rm reg}}  = W \setminus  {\cal D}$.  \\
$  \spec   \, : $  The spectral curve of genus $\tilde{g}
=r^2(g-1)+1$
   associated to $s \in W$. \\
$  \pi _s   :  \spec     \lra  C  \, : $  The projection  to the base
curve
$C$. \\
$\pi :  \unspec  \lra C  \, : $ The universal spectral curve and its
projection
to $C$. \\
$  \beta : \unspec  \lra \so \, : $ The map to $\so $. \\
$ \prym \, : $ The prymian of the map $\pi _s   :  \spec  \lra  C$.
\\
$ \unprym  \, : $ The universal prymian associated to the map $\pi :
\unspec
\lra  C $. \\
$  \prymd  \, : $   The   ``prymian''  of degree
$\tilde{d}=d+r(r-1)(g-1)$,
i.e. the set $\{ \tilde{L} \in J^{\tilde{d}} (\spec ) \; | \; \det
\pi _*
\tilde{L} =  L_0    \}$.   \\
$ \unprymd  \, : $ The universal spectral ``prymian''  of degree
$\tilde{d}$.

\section{Spectral curves}
\subsection{Linear systems of spectral curves}
Let $C$ be a smooth curve of genus $g$. Let $\so$ be the total space
of the
line bundle $\omega_{C} \rightarrow C$ and let $S = {\Bbb
P}(\omega_{C}\oplus
{\cal O}_{C})$ be the projective extension of $\omega_{C}$. Denote by
$\alpha :
S \rightarrow C$ the natural projection and let ${\cal O}_{S}(1)$ be
the
relative hyperplane bundle on $S$ corresponding to the vector bundle
$\omega_{C}\oplus {\cal O}_{C}$.

\begin{defi} An $r$-sheeeted spectral curve is an element
$\widetilde{C}$ of
the linear system $| \alpha^{*}\omega_{C}^{\otimes r}\otimes {\cal
O}_{S}(r)|$
having the property
$\widetilde{C} \subset \so$ \label{def11}
\end{defi}
The adjective spectral refers to another interpretation of the curve
$\widetilde{C}$ which we recall next. For any  vector bundle $E$ of
rank $r$
over $C$ and any $\omega_{C}$-twisted endomorphism $\theta \in
H^{0}(C,{\rm
End}E\otimes\omega_{C})$ of $E$ there is a suitable notion of  a
characteristic
polynomial. Define the $i$-th characteristic coefficient $s_{i} \in
H^{0}(C,\omega_{C}^{\otimes i})$ of  $\theta$ as $s_{i} =
(-1)^{i+1}{\rm
tr}(\wedge^{i}\theta)$. The characteristic polynomial $P(x)$ of
$\theta$ can be
written formally as $P(x) = x^{r} + s_{1}x^{r-1} + \ldots + s_{r}$.
Geometrically
the polynomial $P(x)$ corresponds to a subcheme $\widetilde{C}_{s}
\subset \so$
which via the projection $\alpha$ is an $r$-sheeted branch cover of
$C$.
Indeed, let
$X \in H^{0}(C,\alpha^{*}\omega_{C}\otimes{\cal O}_{S}(1))$ be the
zero section
of the ruled surface $S$ and let $Y \in H^{0}(C,{\cal O}_{S}(1))$ be
the
infinity section of $S$.
By setting $x := \frac{X}{Y}$, i.e. $x$ is the tautological section
of
$\alpha^{*}\omega_{C} \rightarrow \so$, we can reinterpret
$Y^{\otimes r}P(x) =
X^{r} + s_{1}X^{r-1}Y +
\ldots s_{r}Y^{r}$ as a section of $\alpha^{*}\omega_{C}^{\otimes
r}\otimes
{\cal O}_{S}(r)$ whose divisor is contained in $\so$, that is - a
spectral
curve, for more details see \cite{bnr}, \cite{h1}.

\begin{rem} \label{rem11} {\rm
We can view the space $\widetilde{W} := {\displaystyle
\oplus^{r}_{i=1}}H^{0}(C,\omega_{C}^{\otimes i})$ as the space of
characteristic polynomials of $\omega_{C}$-twisted endomorphisms of
rank $r$
vector bundles.  The vector space $\widetilde{W}$ can be identified
explicitly
as the locus of spectral curves in
$| \alpha^{*}\omega_{C}^{\otimes r}\otimes {\cal O}_{S}(r)|$ as
follows. The
push
forward map $\alpha_{*}$ induces an isomorphism
\[ \alpha_{*} : H^{0}(S,\alpha^{*}\omega_{C}^{\otimes r}\otimes {\cal
O}_{S}(r)) \longrightarrow
H^{0}(C,\alpha_{*}(\alpha^{*}\omega_{C}^{\otimes
r}\otimes {\cal O}_{S}(r))).\]
On the other hand,
\[
\begin{array}{lcl}
 \alpha_{*}(\alpha^{*}\omega_{C}^{\otimes r}\otimes {\cal O}_{S}(r))
& = &
\omega_{C}^{\otimes r} \otimes \alpha_{*}({\cal O}_{S}(r)) =
\omega_{C}^{\otimes r} \otimes ({\rm Sym}^{r}({\cal O}_{C} \oplus
\omega_{C}^{-1}) = \\
 & = &  {\cal O}_{C} \oplus \omega_{C}\oplus \ldots \oplus
\omega_{C}^{r}.
  \end{array}
  \]
  Therefore we get an isomorphism
  \[ \alpha_{*} : H^{0}(S,\alpha^{*}\omega_{C}^{\otimes r}\otimes
{\cal
O}_{S}(r)) \longrightarrow \oplus_{i=0}^{r}
H^{0}(C,\omega_{C}^{\otimes i}) ,\]
 whose inverse is
 \[
 \begin{array}{ccc}
  \oplus_{i=0}^{r} H^{0}(C,\omega_{C}^{\otimes i}) & \longrightarrow
&
  H^{0}(S,\alpha^{*}\omega_{C}^{\otimes r}\otimes {\cal O}_{S}(r)) \\
  (s_{0}, s_{1}, \ldots , s_{r}) & \longrightarrow & s_{0}X^{r} +
s_{1}X^{r-1}Y
+
  \ldots + s_{r}Y^{r}.
  \end{array}
  \]
 Here we have identified the spaces $H^{0}(C,\omega_{C}^{\otimes i})$
and
  $\alpha^{*}H^{0}(C,\omega_{C}^{\otimes i}) \subset H^{0}(C,
  \alpha^{*}\omega_{C}^{\otimes i})$ via $\alpha$.

  For any $s = (s_{0}, s_{1}, \ldots , s_{r})$ denote by
$\widetilde{C}_{s}$
the divisor
  ${\rm Div}(  s_{0}X^{r} + s_{1}X^{r-1}Y + \ldots + s_{r}Y^{r})$ on
$S$. Let
  ${\rm Div}(X) = X_{0}$ and ${\rm Div}(Y) = Y_{\infty}$. By
definition the
curve  $\widetilde{C}_{s}$ is spectral if $\widetilde{C}_{s} \subset
\so$, or
equivalently, if
  $\widetilde{C}_{s} \cap Y_{\infty} = \varnothing$. Since $X_{0}
\cap
Y_{\infty} = \varnothing$ it follows that $\widetilde{C}_{s} \cap
Y_{\infty}
\neq \varnothing$ if and only if $Y_{\infty} \subset
\widetilde{C}_{s}$, i.e.
if and only if $s_{0} = 0$.

  Thus the locus of all spectral curves in the projective space
  $| \alpha^{*}\omega_{C}^{\otimes r}\otimes {\cal O}_{S}(r)|$ is the
affine
open set
  $\widetilde{W} = \{ \widetilde{C}_{s} \, | \, s_{0} \neq 0 \}$.
}\end{rem}

  \bigskip

  For the purposes of this paper we will be mainly concerned with the
slightly
smaller
  space
  \[ W := H^{0}(C,\omega_{C}^{\otimes 2})\oplus
H^{0}(C,\omega_{C}^{\otimes
3})\oplus
  \ldots \oplus H^{0}(C,\omega_{C}^{\otimes r}),\]
  consisting of the characteristic polynomials of traceless twisted
endomorphisms of rank
  $r$ vector bundles.  Let $h : \unspec \rightarrow  W$ be the family
of
spectral
  curves parametrized by $W$. To construct the variety $\unspec $
consider the
line bundle
  $p_{S}^{*}(\alpha^{*}\omega_{C}^{\otimes  r}\otimes {\cal
O}_{S}(r))
\rightarrow
  W\times S$. There is a natural tautological section ${\frak c} \in
H^{0}(W\times S,
  p_{S}^{*}(\alpha^{*}\omega_{C}^{\otimes r}\otimes {\cal
O}_{S}(r)))$ given by
  \[ {\frak c}(((s_{2},\ldots ,s_{r});p)) =  X^{r}(p) +
s_{2}(p)X^{r-2}(p)Y^{2}(p) + \ldots +
  s_{r}(p)Y^{r}(p),\]
  and $\unspec = {\rm Div}({\frak c})$ is just the divisor of this
section.

  The universal family $\unspec$ is proper over $W$ and admits a
natural
compactification
  to a projective variety $\overline{\cal C}$ - the universal family
for the
linear system
  ${\Bbb P}(H^{0}(C,{\cal O}_{C})\oplus W)$.  Set $\overline{W} :=
H^{0}(C,{\cal
  O}_{C})\oplus W$.  Next we are going to study the linear system
${\Bbb
  P}(\overline{W})$ and the total spaces of the universal families
$\unspec$
and
  $\overline{\cal C}$.

  \begin{lem} \label{lem11} The linear system  ${\Bbb
P}(\overline{W})$ is
  base point free.
  Furthermore, if $r \geq 3$ the morphism $f_{{\Bbb P}(\overline{W})}
: S
\rightarrow
  {\Bbb P}(\overline{W})^{\vee}$ is an inclusion when restricted to
$\so$ and
contracts
  the infinity section $Y_{\infty}$.
  \end{lem}
  {\bf Proof.}  Let $p \in S$. If $p \in Y_{\infty}$, then the curve
${\rm
Div}(X^{r}) \in
  {\Bbb P}(\overline{W})$ does not pass through $p$. If $p \notin
Y_{\infty}$,
then choose
  a section $s_{r} \in H^{0}(C, \omega_{C}^{\otimes r})$ not
vanishing at
$\alpha(p)$. The
  curve ${\rm Div}(s_{r}Y^{r}) \in {\Bbb P}(\overline{W})$ does not
pass trough
$p$.

  Consider the morphism
  \[f_{|\overline{W}|} : S \longrightarrow {\Bbb
P}(\overline{W})^{\vee}.\]
  Since a curve $\widetilde{C} \in  | \alpha^{*}\omega_{C}^{\otimes
r}\otimes
{\cal O}_{S}(r)|$
  is either spectral or contains $Y_{\infty}$, see Remark
\ref{rem11}, it
follows that
  $f_{{\Bbb P}(\overline{W})}$ contracts $Y_{\infty}$.

  If $p\neq q \in S$ are points in the same fiber $F$ of $\alpha$,
then they
are separated by
  ${\Bbb P}(\overline{W})$ because by pulling ${\Bbb
P}(\overline{W})$ back on
$F$
  we get the linear system
  \[ {\Bbb P}({\rm Span}(x_{0}^{r},x_{0}^{r-2}x_{1}^{2}, \ldots
,x_{1}^{r}))
\subset
  |{\cal O}_{F}(r)|, \]
  which for $r \geq 3$ separates the points on $F \cong {\Bbb
P}^{1}$.

  If $p\neq q \in S$ are points in two different fibers of $\alpha$,
then by
choosing a section
  $s_{r} \in  H^{0}(C, \omega_{C}^{\otimes r})$ satisfying
$s_{r}(\alpha(p)) =
0$ and
  $s_{r}(\alpha(q)) \neq 0$, we get a curve ${\rm Div}(s_{r}Y^{r})$
passing
trough $p$ and
  not passing through $q$.

  The above considerations and the local triviality of $\so
\rightarrow C$ show
that
   ${\Bbb P}(\overline{W})$ separates also the tangent directions
because for
every point
   $p \in \so$ the morphism   $f_{{\Bbb P}(\overline{W})}$ maps the
fiber to
hrough $p$
   and a suitable translate of $X_{0}$ into two transversal curves in
   ${\Bbb P}(\overline{W})^{\vee}$.
   \begin{flushright}
   $\Box$
   \end{flushright}

   \bigskip

   Here are some immediate corollaries from the above lemma which are
going to
be
   usefull later on.

   \begin{cor} \label{cor11} The generic spectral curve $\spec$ is
smooth.
   \end{cor}
   {\bf Proof.} Follows from the fact that ${\Bbb P}(\overline{W})$
does not
have base points and from Bertini's theorem.
   \begin{flushright}
   $\Box$
   \end{flushright}

   \begin{cor} \label{cor12} The total spaces of the universal
families
$\unspec$ and
   $\overline{\cal C}$ are smooth.
   \end{cor}
   {\bf Proof.} Since $\unspec \subset  \overline{\cal C}$ is a
Zariski open
set it is enough to
   show that $\overline{\cal C}$ is smooth.

   Let $\pi : \overline{U} \rightarrow {\Bbb P}(\overline{W})^{\vee}$
be the
universal bundle of
   linear hypersurfaces in ${\Bbb P}(\overline{W})$. The total space
$\overline{\cal C}$
   fits in the fiber product diagram
   \[
   \begin{diagram}
   \node{\overline{\cal C}} \arrow{e} \arrow{s} \node{\overline{U}}
\arrow{s,r}{\pi} \\
   \node{S} \arrow{e,b}{f_{{\Bbb P}(\overline{W})}} \node{{\Bbb
P}(\overline{W})^{\vee}}
   \end{diagram}
   \]
   But both $S$ and $\overline{U}$ are smooth varieties and moreover
the
projection
   $\pi : \overline{U}  \rightarrow {\Bbb P}(\overline{W})^{\vee}$ is
a smooth
map. Therefore
   the fiber product  $\overline{\cal C} = S_{f_{|\overline{W}|}}
\!\!
\times_{\pi}{\cal
   P}$ is also smooth.
  \begin{flushright}
  $\Box$
  \end{flushright}

  \bigskip

  Let $\overline{W}_{r-1,r}$ (resp.  $W_{r-1,r}^{\infty}$)   be the
linear
subsystem
   of  $\overline{W}$ consisting of points
  of the form $(s_0,0,\ldots , s_{r-1}, s_{r})$ (resp. $(0,0,\ldots ,
s_{r-1},
s_{r})$).
   According to Remark \ref{rem11}, $\overline{W}_{r-1,r}$ (resp.
$W_{r-1,r}^{\infty}$)
   can be embedded naturally in
   the vector space $H^{0}(S, \alpha^{*}\omega_{C}^{\otimes r}\otimes
{\cal
   O}_{S}(r))$ and thus, $|(\overline{W}_{r-1,r}|$ (resp.
   $|(W^{\infty}_{r-1,r}|$)  can be viewed as a linear  subsystem
    of $|\alpha^{*}\omega_{C}^{\otimes r}\otimes {\cal O}_{S}(r)|$.
The
arguments
   in the proof of Lemma \ref{lem11} yield the following corollary.

   \begin{cor} \label{cor13}  The linear system
$\overline{W}_{r-1,r}$
 (resp.  $W_{r-1,r}^{\infty}$)  is  base point free.
  Furthermore if $r \geq 3$, then the morphism
$f_{|\overline{W}_{r-1,r}|} : S
\rightarrow
    {\Bbb P}(\overline{W}_{r-1,r})^{\vee}$ (resp.  $f_{|W^{\infty}
_{r-1,r}|} :
S \rightarrow
  {\Bbb P}(W^{\infty}_{r-1,r})^{\vee}$)        is an inclusion when
restricted
to $\so$ and
  contracts the infinity section $Y_{\infty}$.
  \end{cor}

  \bigskip

  In the next section we are going to prove the irreducibility of
various
discriminant loci.
  The following fact is an essential ingredient in the proof of
Corollary
\ref{cor16} and
  is also of independent interest.

  \begin{prop} \label{prop10} Let $\widetilde{C}$ be a spectral curve
of degree
$r$ and let
  $\Sigma \subset \widetilde{C}$ be an irreducible component of
$\widetilde{C}$. Then
   $\Sigma$ is a spectral curve of degree $l \leq r$.
 \end{prop}
 {\bf Proof.} Since $\Sigma \subset \widetilde{C} \subset \so$, we
only need to
show that
 \[ {\cal O}_{S}(\Sigma) = {\cal
O}_{S}(l)\otimes\alpha^{*}\omega_{C}^{\otimes
l}, \]
 for some $l \leq r$. But ${\rm Pic}(S) = {\Bbb Z}\cdot {\cal
O}_{S}(1) \oplus
 \alpha^{*}{\rm Pic}(C)$, and hence
 \[ {\cal O}_{S}(\Sigma) = {\cal O}_{S}(l)\otimes\alpha^{*}L, \]
 for some line bundle $L$ on $C$. Since $\Sigma \subset \so$ we have
$\Sigma\cdot
 Y_{\infty} = 0$ which yields
 \[ 0 = {\cal O}_{S}(\Sigma)\cdot{\cal O}_{S}(1) = {\cal
O}_{S}(l)\cdot {\cal
O}_{S}(1) +
 \alpha^{*}L\cdot {\cal O}_{S}(1) = - l(2g-2) + \deg L,\]
 i.e. $\deg L = l(2g - 2)$.

 On the other hand the line bundle ${\cal
O}_{S}(l)\otimes\alpha^{*}L$ has a
section
 $\Sigma$ which does not vanish on the infinity divisor $Y_{\infty}$.
But the
sections of
 ${\cal O}_{S}(l)\otimes\alpha^{*}L$ vanishing on $Y_{\infty}$ are
$H^{0}(S,{\cal
 O}_{S}(l)\otimes\alpha^{*}L\otimes {\cal O}_{S}(-Y_{\infty})) =
 H^{0}(S,{\cal O}_{S}(l-1)\otimes\alpha^{*}L)$.
 By pushing forward ${\cal O}_{S}(l-1)\otimes\alpha^{*}L$ on $C$ we
get
 \[
 \begin{array}{lcl}
 H^{0}(S,{\cal O}_{S}(l-1)\otimes\alpha^{*}L) & = &
 H^{0}(C,\alpha_{*}({\cal O}_{S}(l-1))\otimes L) = H^{0}(C,{\rm
Sym}^{l-1}({\cal O}_{C}\oplus
 \omega_{C}^{-1})\otimes L) = \\
  & = & H^{0}(C,L\oplus (L\otimes\omega_{C}^{-1})\oplus \ldots \oplus
  (L\otimes \omega_{C}^{-(l-1)})) = \oplus_{i=0}^{l-1} H^{0}(C,
L\otimes
  \omega_{C}^{-i}).
  \end{array}
  \]
  Similarly $H^{0}(S,{\cal O}_{S}(l)\otimes\alpha^{*}L) =
  \oplus_{i=0}^{l} H^{0}(C, L\otimes \omega_{C}^{-i})$ and hence
  \[ H^{0}(S,{\cal O}_{S}(l)\otimes\alpha^{*}L)/H^{0}(S,{\cal
  O}_{S}(l-1)\otimes\alpha^{*}L) \cong  H^{0}(C, L\otimes
\omega_{C}^{-l}).\]
  The section of the divisor $\Sigma$ gives a non-zero element in
  $H^{0}(S,{\cal O}_{S}(l)\otimes\alpha^{*}L)/H^{0}(S,{\cal
  O}_{S}(l-1)\otimes\alpha^{*}L)$ and therefore $H^{0}(C, L\otimes
\omega_{C}^{-l})
  \neq 0$.  Since $\deg (L\otimes \omega_{C}^{-l}) = 0$ we get that
$L\otimes
\omega_{C}^{-l} \cong {\cal O}_{C}$.
  \begin{flushright}
  $\Box$
  \end{flushright}

  \begin{rem} \label{rem10} {\rm The meaning of the above proposition
becomes
  transparent if we use the interpretation of the spectral curves as
characteristic
  polynomials of  $\omega_{C}$-twisted endomorphisms.

  Let $F \rightarrow \widetilde{C}$ be a line bundle. Then $F$ can be
viewed as
a sheaf
  on $S$ supported on $\widetilde{C} \subset \so$. Consider the rank
$r$ vector
bundle
  $E = \alpha_{*}F$ on $C$. The push-forward of the homomorphism
  \[ F \stackrel{\otimes x}{\longrightarrow}
F\otimes\alpha^{*}\omega_{C}\]
  is a $\omega_{C}$-twisted endomorphism of $E$:
  \[ \theta : E \longrightarrow E\otimes\omega_{C}.\]
  In this way the curve $\widetilde{C}$ then can be described as the
zero
scheme of the section
  \[ \det (\alpha^{*}\theta - x\cdot {\rm id}_{\alpha^{*}E}) \in
H^{0}(S,
  \alpha^{*}\omega_{C}^{\otimes r}\otimes {\cal O}_{S}(r)).\]

  Let now $\Sigma \subset \widetilde{C}$ be a component of
$\widetilde{C}$.
Consider the
  sheaf $F\otimes {\cal O}_{\Sigma}$ on $S$ supported on $\Sigma$.
Let $E' :=
  \alpha_{*}(F\otimes {\cal O}_{\Sigma})$ and let
  \[ \theta' :  E' \longrightarrow E'\otimes\omega_{C}\]
  be the push forward of the homomorphism
  \[ F\otimes {\cal O}_{\Sigma} \stackrel{\otimes x}{\longrightarrow}
(F\otimes
{\cal
  O}_{\Sigma})\otimes\alpha^{*}\omega_{C}.\]
  The commutative diagram of (torsion) sheaves on $\so$
  \[
  \begin{diagram}
  \node{F} \arrow{e,t}{\otimes x} \node{F\otimes\alpha^{*}\omega_{C}}
\\
  \node{F\otimes {\cal O}_{\Sigma}} \arrow{n} \arrow{e,t}{\otimes x}
  \node{(F\otimes {\cal O}_{\Sigma})\otimes\alpha^{*}\omega_{C}}
\arrow{n}
  \end{diagram}
  \]
  then pushes down to a commutative diagram of vector bundles on $C$:
  \[
  \begin{diagram}
  \node{E} \arrow{e,t}{\theta} \node{E\otimes\omega_{C}} \\
  \node{E'} \arrow{e,t}{\theta'} \arrow{n} \node{E'\otimes\omega_{C}}
  \arrow{n}
  \end{diagram}
  \]
  In particular $E'$ is $\theta$ invariant and $\theta' =
\theta_{|E'}$. Then
$\Sigma$ is
  the zero scheme of the section
  \[det (\alpha^{*}\theta' - x\cdot {\rm id}_{\alpha^{*}E'}) \in
H^{0}(S,
  \alpha^{*}\omega_{C}^{\otimes l}\otimes {\cal O}_{S}(l)),\]
  where $l = {\rm rk}E'$, i.e. $\Sigma$ is a spectral curve.}
  \end{rem}

  \begin{cor} \label{cor14} The locus $R \subset W$ consisting of
reducible
spectral curves
  has codimension \linebreak ${\rm codim}(R,W) \geq g-1$
  \end{cor}
  {\bf Proof.} Set $W_{i} := H^{0}(C,\omega_{C}^{\otimes i})$. Let
$R_{r}$ be
the
  image of $R$ under the natural projection $W \rightarrow W_{r}$. It
suffices
to show
  that ${\rm codim}(R_{r},W_{r}) \geq g-1$ or equivalently, since
$R_{r}$ is
invariant under
  dilations, that  ${\rm codim}({\Bbb P}(R_{r}),{\Bbb P}(W_{r})) \geq
g-1$.

  For every $i = 1, \ldots, r-1$ consider the variety $S_{i} \subset
  {\Bbb P}(W_{r})\times {\Bbb P}(W_{i})$ defined by
  \[ S_{i} = \{ (D,G) \; | \; D \geq G\}. \]
  Since a component of a spectral curve is a spectral curve we have
  \[  {\Bbb P}(R_{r}) \subset p_{{\Bbb P}(W_{r})}(S_{1}) \cup \ldots
\cup
  p_{{\Bbb P}(W_{r})}(S_{r-1}),\]
 and therefore ${\rm codim}({\Bbb P}(R_{r}),{\Bbb P}(W_{r})) \geq
\min_{1\leq i
\leq
   r-1}\{ {\rm codim}(p_{{\Bbb P}(W_{r})}(S_{i}),{\Bbb P}(W_{r}))\}$.
Furthermore
   \linebreak
   $\dim (p_{{\Bbb P}(W_{r})}(S_{i})) = \dim S_{i}$ because the map
   $p_{{\Bbb P}(W_{r})} : S_{i} \longrightarrow {\Bbb P}(W_{r})$ is
finite on
  its image. To compute $\dim S_{i}$  look at the second projection
  \begin{equation} \label{eq11}
  p_{{\Bbb P}(W_{i})} : S_{i} \longrightarrow {\Bbb P}(W_{i}).
  \end{equation}
  The map (\ref{eq11}) is obviously onto and for any $G \in {\Bbb
P}(W_{i})$ we
get by
  Riemann-Roch
  \[ \dim p_{{\Bbb P}(W_{i})}^{-1}(G) = h^{0}(C,\omega_{C}^{\otimes
r}(-G)) - 1
=
   h^{0}(C,\omega_{C}^{\otimes (r-i)}) - 1 \leq (2(r-i) - 1)(g-1). \]
   Consequently by the fiber-dimesion theorem we get
 $\dim S_{i} \leq   (2r - 2)(g-1) - 1$ and hence
 \[ {\rm codim}(p_{{\Bbb P}(W_{r})}(S_{i}),{\Bbb P}(W_{r})) = g-1, \]
 for any $i = 1, \ldots , r-1$.
 \begin{flushright}
 $\Box$
 \end{flushright}

  \subsection{Discriminant loci} \label{ss11}
  As we saw in Corollary \ref{cor11} the generic spectral curve is
smooth.
Therefore
  the space $W$ of spectral curves contains a natural divisor
parametrizing the
singular
  spectral curves:
  \[ {\cal D} := \{ s \in W \, | \, \spec{\rm - is \; not \; smooth}
\}, \]
  which due to Corollary \ref{cor12} is just the discriminant locus
for the map
$h : \unspec
  \rightarrow W$.

  The divisor ${\cal D}$ is irreducible. To see this consider the
full family
  $\widetilde{W}$ of spectral curves and its discriminant divisor
$\widetilde{\cal D}$.
  We have a natural inclusion $W \hookrightarrow \widetilde{W}$ and
clearly
  ${\cal D} = W \cap \widetilde{\cal D}$. The first step will be to
show that
  $\widetilde{\cal D}$ is irreducible and the second to deduce the
irreducibility of ${\cal
  D}$ from that.  To achieve this we need to introduce some
auxilliary objects.

  Let $\aff $ be the additive group of the vector space
$H^{0}(C,\omega_{C})$.
There is a
  natural affine action of $\aff$ on the bundle $\omega_{C}$ which
can be
considered as an
  action on its total space
  \[
  \begin{array}{lccl}
  \tau : & \aff & \longrightarrow & {\rm Aut}(\so) \\
     & \gamma & \longrightarrow & (p
\stackrel{\tau_{\gamma}}{\mapsto} p +
     \gamma(\alpha(p))).
  \end{array}
  \]
  The action $\tau$ on $\so$ has a natural lift to an action on the
bundle
  $\alpha^{*}\omega_{C} \rightarrow \so$ which gives an affine action
of $\aff$
on the
  space of its global sections: $\gamma \rightarrow (\mu \mapsto \mu
+
  \alpha^{*}\gamma)$.  Consequently we obtain a polynomial action on
the space
of
  spectral curves
  \[
  \begin{array}{lccl}
  \rho : & \aff & \longrightarrow & {\rm Aut}(\widetilde{W}) \\
    & \gamma & \longrightarrow & ( x^{r} + s_{1}x^{r-1} + \ldots +
s_{r}
    \stackrel{\rho_{\gamma}}{\mapsto}  (x+\gamma)^{r} +
s_{1}(x+\gamma)^{r-1} +
    \ldots + s_{r}).
  \end{array}
  \]
  For every $\gamma \in \aff$ we have a commutative diagram
  \begin{equation} \label{eq12}
  \begin{diagram}
  \node{\unspec} \arrow{e,t}{\tau_{\gamma}} \arrow{s,r}{h}
\node{\unspec}
  \arrow{s,l}{h} \\
  \node{\widetilde{W}} \arrow{e,t}{\rho_{\gamma}}
\node{\widetilde{W}}
  \end{diagram}
  \end{equation}
  where $\tau_{\gamma}$ is the automorphism of the universal spectral
curve
induced by
  $\tau_{\gamma} \in {\rm Aut}(\so)$.

  \begin{prop} \label{prop12} The discriminant divisor
$\widetilde{\cal D}
\subset
  \widetilde{W}$ is irrreducible.
  \end{prop}
  {\bf Proof.} Consider the incidence correspondence
$\widetilde{\Gamma}\subset
  \widetilde{W}\times \so$ defined by
  \[  \widetilde{\Gamma} := \{ (s,p) \in \widetilde{W}\times \so \; |
\; {\rm
Ord}_{p}\spec
  \geq 2 \}. \]
  Since $\widetilde{\cal D} = p_{\widetilde{W}}(\widetilde{\Gamma})$
it
suffices to show
  that $\widetilde{\Gamma}$ is irreducible. The commutativity of the
diagram
(\ref{eq12})
  implies that $\widetilde{\Gamma}$ is $\rho$-invariant. Furthermore,
if $s' =
  \rho_{\gamma}(s)$, then $s_{1}' = s_{1} + r\gamma$ and hence $\rho$
is a free
action.
  Therefore we can form the quotient variety
$\widetilde{\Gamma}/\aff$. Since
the
  fibration
  \[ \widetilde{\Gamma} \longrightarrow  \widetilde{\Gamma}/\aff \]
  is a locally trivial affine bundle, the irreducibility of
$\widetilde{\Gamma}$ is equivalent
  to the irreducibility of $\widetilde{\Gamma}/\aff$.

  The fibers of the projection
  \[ p_{C} := \alpha\circ p_{\so} : \widetilde{\Gamma}
\longrightarrow C \]
  are $\rho$-invariant and by taking the quotient we obtain a
morphism
  \[ p_{C} : \widetilde{\Gamma}/\aff  \longrightarrow C.\]
  Let $\Xi \subset \widetilde{W}\times C$ be the incidence
correspondence
  \[ \Xi  = \{ (s,a) \in \widetilde{W}\times C \; | \; \spec \; {\rm
is \;
singular \; at \; the \; point} \;
  \alpha^{-1}(a) \cap X_{0} \}. \]
  More explicitly by using the Jacobian criterion for smoothnes one
gets, see
\cite{bnr}
  Remark 3.5,
  \[  \Xi  = \{ (s,a) \in \widetilde{W}\times C \; | \; {\rm
Div}(s_{r}) \geq
2a, \; {\rm
  Div}(s_{r-1})
  \geq a \}.\]
  The fact that $\omega_{C}^{\otimes r}$ is very ample for any $r\geq
2$,
implies that
  the fibration $\Xi \rightarrow C$ is a vector bundle of rank
$\dim\widetilde{W} - 3$.
  On the other hand we have a natural morphism of fibrations
  \[
  \begin{diagram}
  \node{\Xi} \arrow[2]{e,t}{t} \arrow{se}
\node[2]{\widetilde{\Gamma}/\aff}
  \arrow{sw,r}{p_{C}} \\
  \node[2]{C}
  \end{diagram}
  \]
  where $t((s,a)) = {\rm Orb}_{\aff}((s,\alpha^{-1}(a)\cap X_{0}))$.
It is easy
to see that $t$ is onto.
  Indeed, if $(s,p) \in \widetilde{\Gamma}$, then ${\rm
Ord}_{p}(\spec) \geq
2$. Choose
  $\gamma \in H^{0}(C,\omega_{C})$ with the property
$\gamma(\alpha(p)) = p$.
Then
  \[(\rho_{-\gamma},\tau_{-\gamma})\cdot (s,p) = (\rho_{-\gamma}(s),
  \alpha^{-1}(\alpha(p))\cap X_{0})), \]
  and hence
  \[ {\rm Orb}_{\aff}((s,p)) = {\rm Orb}_{\aff}((\rho_{-\gamma}(s),
  \alpha^{-1}(\alpha(p))\cap X_{0}))). \]
  Consequently $t$ is onto and $\widetilde{\Gamma}/\aff$ is
irreducible.
  \begin{flushright}
  $\Box$
  \end{flushright}

  \begin{rem} \label{rem12} {\rm The variety
$\widetilde{\Gamma}/\aff$  (and
therefore
  $\widetilde{\Gamma}$) is actually smooth. Indeed, let $\Xi_{a}$ be
the fiber
of the
  bundle $\Xi$ over the point $a \in C$. Let $(s,a), (s',a) \in
\Xi_{a}$ be
such that
  $t((s,a)) = t((s',a))$. Then
  \[  {\rm Orb}_{\aff} ((s,\alpha^{-1}(a)\cap X_{0})) =  {\rm
Orb}_{\aff}
   ((s',\alpha^{-1}(a)\cap X_{0})),\]
  and therefore $s' = \rho_{\gamma}(s)$ for some $\gamma \in
H^{0}(C,\omega_{C})$
  satisfying $\gamma(a) = 0$.

  Let $\Xi_{0} \subset \Xi$ be the subbundle
  \[ \Xi_{0}  = \{ (s,a) \in \Xi \; | \; {\rm Div}(s_{1}) \geq a
\}.\]
  The morphism $t$ descends to an isomorphism of fibrations
  \[
  \begin{diagram}
  \node{\Xi /\Xi_{0}} \arrow[2]{e,t}{t} \arrow{se}
\node[2]{\widetilde{\Gamma}/\aff}
  \arrow{sw,r}{p_{C}} \\
  \node[2]{C}
  \end{diagram}
  \]
  and thus $\widetilde{\Gamma}/\aff \cong  \Xi /\Xi_{0}$ is smooth.
Furthermore,  observe
  that $p_{\widetilde{W}} : \widetilde{\Gamma} \rightarrow
\widetilde{\cal D}$
is
  birational because for the generic $s \in \widetilde{\cal D}$ the
curve
$\spec$ has a unique
ordinary double point. Therefore $\widetilde{\Gamma}$ can be viewed
as a
natural
  desingularization of $\widetilde{\cal D}$. }
  \end{rem}

   \begin{cor} \label{cor15} The discriminant divisor ${\cal D}
\subset W$ is
irreducible.
   \end{cor}
   {\bf Proof.} Consider the incidence correspondence
   \[ \Gamma = \{ (s,p) \; | \; {\rm Ord}_{p}(\spec ) \geq 2 \}. \]
   Let $q_{|\Gamma} :\Gamma \rightarrow  \widetilde{\Gamma}/\aff$ be
the
restriction to
   $\Gamma$ of the natural quotient morphism $q : \widetilde{\Gamma}
\rightarrow
   \widetilde{\Gamma}/\aff$. If $(s,p) \in \widetilde{\Gamma}$, then
   $(\rho_{-\frac{1}{r}\cdot s_{1}}(s), \tau_{-\frac{1}{r}\cdot
s_{1}}(s)) \in
\Gamma$.
   Combined with the fact that $\rho$ acts freely, this yields
   \[ \Gamma \cap   {\rm Orb}_{\aff}  ((s,p)) = \{
(\rho_{-\frac{1}{r}\cdot
s_{1}}(s),
    \tau_{-\frac{1}{r}\cdot
   s_{1}}(s)) \}. \]
   Consequently $q_{|\Gamma}$ is an isomorhism and in particular
$\Gamma$ and
   ${\cal D}$ are irreducible.
   \begin{flushright}
   $\Box$
   \end{flushright}

   There are two other discriminant loci whose irreducibility will be
usefull.
Define
   $W_{r} := H^{0}(C,\omega_{C}^{\otimes r})$ and $W_{r-1,r} :=
   H^{0}(C,\omega_{C}^{\otimes (r-1)})\oplus
H^{0}(C,\omega_{C}^{\otimes r})$.
   Let ${\cal D}_{r} = {\cal D}\cap W_{r}$ and ${\cal D}_{r-1,r}
={\cal D}\cap
W_{r-1,r}$
   be the discriminant divisors for the families of spectral curves
$W_{r}$ and
   $W_{r-1,r}$ respectively.

   \begin{prop} \label{prop13} The discriminant divisor ${\cal D}_{r}
\subset
W_{r}$
   is irreducible.
   \end{prop}
   {\bf Proof.} Consider the incidence correspondence $\Gamma_{r}
\subset
W_{r}\times
   C$ defined by
   \[ \Gamma_{r} = \{ (s_{r},a) \; |\; {\rm Div}(s_{r}) \geq 2a \}.\]
   Using the Jacobian criterion it is easy to check that a curve
$\spec$, $s
\in W_{r}$ is
   singular if and only if $s_{r}$ has a double zero.

   Therefore ${\cal D}_{r} = p_{W_{r}}(\Gamma_{r})$ and it suffices
to show
that
   $\Gamma_{r}$ is irreducible.  The divisor $\Gamma_{r}$ is
equivariant with
respect to
   the natural ${\Bbb C}^{\times}$ action on $W_{r}$ (extended
trivially to
   $W_{r}\times C$). Thus the problem reduces to showing that ${\Bbb
P}(\Gamma_{r})
   \subset {\Bbb P}(W_{r})\times C$ is irreducible.

   Consider the projection $p_{C} : {\Bbb P}(\Gamma_{r}) \rightarrow
C$. By
   Riemann-Roch $p_{C}$ is onto and ${\rm codim}(p_{C}^{-1}(a),{\Bbb
P}(\Gamma_{r}))
   \linebreak = 2$
   for every $a \in C$. Moreover, $p_{C}^{-1}(a) \subset {\Bbb
P}(W_{r})$ is a
linear
   subspace and hence all the fibers of $p_{C}$ are equidimensional
and
irreducible.
   Therefore ${\Bbb P}(\Gamma_{r})$ is irreducible.
   \begin{flushright}
   $\Box$
   \end{flushright}

   \begin{rem} \label{rem13} {\rm The divisor ${\Bbb P}({\cal
D}_{r})$ is a
divisor in the
   projective space of the complete linear system ${\Bbb
P}(H^{0}(C,\omega_{C}^{\otimes
   r}))$. This allows us to give more geometric description of
${\Bbb P}({\cal
D}_{r})$:
   \[
   {\Bbb P}({\cal D}_{r}) = {\Bbb P}(\{ H \subset {\Bbb
P}(W_{r}^{\vee}){\rm -
   hyperplane}\; | \; H \supset {\Bbb P}(T_{a}f_{|W_{r}|}(C)) \; {\rm
for \;
some \;} a \}), \]
   that is, ${\Bbb P}({\cal D}_{r}) = f_{|W_{r}|}(C)^{\vee}$ is the
dual
hypersurface of the
   $r$-canonical model of the curve $C$.}
   \end{rem}

   \begin{prop} \label{prop14} The divisor ${\cal D}_{r-1,r} \subset
W_{r-1,r}$
is
   irreducible.
   \end{prop}
   {\bf Proof.} Consider the incidence variety
   \[ \Gamma_{r-1,r} = \{ (a,b,p) \in W_{r-1,r}\times \so \; | \;
   {\rm Ord}_{p}(x^{r} + ax + b) \geq 2 \}. \]
   Again $p_{W_{r-1,r}}(\Gamma_{r-1,r}) = {\cal D}_{r-1,r}$ so it
suffices to
show that
   $\Gamma_{r-1,r}$ is irreducible.

   Consider the projection
   \[ p_{\so} : \Gamma_{r-1,r} \longrightarrow \so. \]
   Let $(U,z)$ be a coordinate chart on $C$. Then we can choose
natural
coordinates
   $(w,z) : \alpha^{-1}(U) \widetilde{\rightarrow} {\Bbb C}\times U$
on
$\alpha^{-1}(U)$,
   so that the tautological section $x$ is $x = w\alpha^{*}(dz)$.

   Let $(z_{0},w_{0}) \in \alpha^{-1}(U)$. Then for any $(a,b) \in
W_{r-1,r}$
one has in a
   neighborhood of $(z_{0},w_{0})$:
   \[
   \begin{array}{lclcl}
   a & = & (a_{0} + a_{1}(z-z_{0}) + \ldots)dz^{\otimes (r-1)} & = &
f(z)dz^{\otimes (r-1)} \\
   b & = & (b_{0} + b_{1}(z-z_{0}) + \ldots)dz^{\otimes r} & = &
g(z)dz^{\otimes r}.
   \end{array}
   \]
   Therefore in the local coordinates $(z,w)$ the equation of the
curve $x^{r}
+ ax + b$
   becomes
   \begin{equation} \label{eq13}
   w^{r} + f(z)w + g(z) = 0.
   \end{equation}
   The conditions for (\ref{eq13}) to have singularity at
$(z_{0},w_{0})$ are
   \[
   \begin{array}{rcl}
   b_{0} & = & 0 \\
   (rw^{r-1} + f(z))_{|(z_{0},w_{0})}  & = & 0 \\
   (f'(z)w + g'(z))_{|(z_{0},w_{0})}  & = & 0
   \end{array}
   \]
   or equivalently
   \begin{equation} \label{eq14}
   \begin{array}{rcl}
   b_{0} & = & 0 \\
   rw^{r-1}_{0} +a_{0} & = & 0 \\
   a_{1}w_{0} + b_{1} & = & 0.
   \end{array}
   \end{equation}
   The equations (\ref{eq14}) determine a codimension 3 affine
subspace of
$W_{r-1,r}$ and
   hence
   ${\Gamma_{r-1,r}}_{|\alpha^{-1}(U)} \longrightarrow
\alpha^{-1}(U)$
   is a trivial affine bundle. Therefore $\Gamma_{r-1,r}$ is an
affine bundle
over $\so$
   and hence is irreducible.
   \begin{flushright}
   $\Box$
   \end{flushright}

   \begin{cor} \label{cor16} Let $h : \unspec \rightarrow W$ be the
universal
spectral
   curve and let $h_{r-1,r} : \unspec_{|W_{r-1,r}} \rightarrow
W_{r-1,r}$ and
   $h_{r} : \unspec_{|W_{r}} \rightarrow W_{r}$ be the subfamilies
parametrized
by
   $W_{r-1,r}$ and $W_{r}$ respectively. Then the divisors
$h^{-1}({\cal D})$,
   $h_{r-1,r}^{-1}({\cal D}_{r-1,r})$ and $h^{-1}_{r}({\cal D}_{r})$
are
irreducible.
   \end{cor}
   {\bf Proof.}     Recall the following standard lemma, see
\cite{sh},
   \begin{lem} \label{lem12} Let $f : X \rightarrow Y$ be a proper
map between
algebraic
   varieties of pure dimension. If $Y$ is irreducible, the general
fiber of $f$
is irreducible
   and all the fibers are equidimensional then $X$ is irreducible.
   \end{lem}

   As we saw above each of  the discriminant loci ${\cal D}_{r-1,r}$,
   ${\cal D}_{r}$ and ${\cal D}$ is irreducible.  According to Lemma
\ref{lem12} it
   suffices to show that the generic fiber of each of the maps $h$,
$h_{r-1}$
and $h_{r-1,r}$
   is irreducible.  But, by Corollary \ref{cor14}, there exists a
point $s \in
{\cal D}_{r}
   \subset {\cal D}_{r-1,r} \subset {\cal D}$ such that $\spec$ is
irreducible
which
   finishes the proof.
   \begin{flushright}
   $\Box$
   \end{flushright}

   \subsection{The Hitchin map} \label{ss13}
   For a line bundle $L_{0} \in {\rm Pic}^{d}(C)$ denote by $\su$ the
moduli
space of
   semistable vector bundles of rank $r$ and determinant $L_{0}$. It
is well
known
   that $\su$ is a normal projective variety whose smooth locus
coincides with
the locus
   $\sustable$ of stable vector bundles, see \cite{se}. Since $C$ is
a curve,
the deformations
   of any vector bundle $E \rightarrow C$ are unobstructed and the
Kodaira-Spencer map
   \begin{equation} \label{eq15}
  T_{[E]}\sustable \longrightarrow H^{1}(C,{\rm End}^{o}E),
  \end{equation}
   is an isomorphism; here ${\rm End}^{o}E$ is the bundle of
traceless
endomorphisms
   of $E$.

   Using the isomorphism (\ref{eq15}) and Serre's duality, one gets a
canonical
   identification
   \[ T_{[E]}^{*}\sustable \cong H^{0}(C,{\rm
End}^{o}E\otimes\omega_{C}), \]
   of the fiber of the cotangent bundle to $\sustable$ at the point
$[E]$ and
the space of
   $\omega_{C}$-twisted traceless endomorphisms. Therefore the
collection of
   characteristic coefficients defined in Section \ref{ss11} gives
rise to a
morphism
   \begin{equation} \label{eq16}
   \begin{array}{lccl}
   H \; : \; & T^{*}\sustable & \longrightarrow & W \\
    & (E,\theta) & \longrightarrow & (s_{2}(\theta),
s_{3}(\theta),\ldots ,
s_{r}(\theta)).
    \end{array}
    \end{equation}
    Hitchin \cite{h1} was the first one to study the morphism
(\ref{eq16}) in
various
    situations.
    He discovered many of its remarkable properties and in particular
he showed
that
    $H$ endows the (holomorphic) symplectic manifold $\cotangent :=
T^{*}\sustable$
    with a structure of an algebraically completely integrable
hamiltonian
system.
    More specifically, he showed that there exists a partial
compactification
    \begin{equation} \label{eq17}
    \begin{diagram}
    \node{\cotangent} \arrow[2]{e} \arrow{se,r}{H} \node[2]{\higgs}
    \arrow{sw,r}{H} \\
    \node[2]{W}
    \end{diagram}
    \end{equation}
    with general fiber an abelian variety and such that the fibers of
(\ref{eq16}) embed as
    Zariski open sets in the fibers of $H : \higgs \rightarrow W$.
The variety
$\higgs$ is
    again a moduli space which parametrizes the equivalence classes
of
semistable Higgs
    pairs - that is, pairs $(E,\theta)$ of a bundle and an
$\omega_{C}$-twisted
    endomorphism such that for every $\theta$-invariant subbundle $U
\subset E$
the
     usual inequality $\mu(U) \leq \mu(E)$ for the slopes of $U$ and
$E$ holds.
The
     characteristic coefficients map (\ref{eq16}) is a well defined
    morphism on the whole variety $\higgs$ and is called the Hitchin
map of
$\higgs$.
    There is a natural ${\Bbb C}^{\times}$-action on the moduli space
$\higgs$
    \[ (E,\theta) \longrightarrow (E,t\theta), \]
    which induces via the Hitchin map a weighted ${\Bbb
C}^{\times}$-action on
the
    vector space $W$. A number $t \in {\Bbb C}^{\times}$ acts on $W$
by
multiplying
    the piece $H^{0}(C,\omega_{C}^{\otimes i})$ with $t^{i}$.
    The
    compactification diagram (\ref{eq17}) and its generalizations
have
    been studied extensively in the last years, see \cite{bnr},
\cite{h1},\cite{h2}, \cite{ni},
    \cite{si}. For our purposes, the most convenient is the approach
in
\cite{bnr} which we
    proceed to describe.

    \begin{prop}[Beauville-Narasimhan-Ramanan]\label{propbnr} Assume
that
$\pi_{s} :
    \spec \rightarrow C$ is an integral spectral curve. Then the
push-forward
map
    $\pi_{s*}$ induces a canonical bijection between
    \begin{itemize}
    \item Isomorphism classes of rank one torsion-free sheaves $F$ on
$\spec$
of Euler
    characteristic $d - r(g-1)$ satisfying $\det (\pi_{s*}F) =
L_{0}$.
    \item Isomorphism classes of Higgs pairs $(E,\theta)$ satisfying
$\det E =
L_{0}$,
    ${\rm tr}\theta = 0$, $H(\theta) = s$.
    \end{itemize}
    \end{prop}

    \begin{rem} \label{rem14} {\rm Let $s \in W$ be such that $\spec$
is
integral. Then
    every
    Higgs pair satisfying $H(\theta) = s$ is stable. Indeed, if we
assume that
there exists a
    proper $\theta$-invariant subbundle $U \subset E$,  then the
spectral curve
of
    $(U,\theta_{|U})$ will be contained in $\spec$ which contadicts
the
integrality of $\spec$.
    Consequently $E$ does not have $\theta$-invariant subbundles and
in
particular
    $(E,\theta)$ is stable.

   Under the bijection in Proposition \ref{propbnr} this translates
into the
well known fact
   that the moduli functor of rank one torsion free sheaves on an
integral
curve possesing
   only planar  singularities is representable by an irreducible fine
moduli
space - the
   compactified Jacobian of the curve, see \cite{aik}. Therefore the
push-forward map
   $\pi_{s*}$ induces an isomorphism
   \[  \pi_{s*} : {\rm Prym}_{\tilde{d}}(\spec,  C)  \longrightarrow
H^{-1}(s)
, \]
   where ${\rm Prym}_{\tilde{d}}(\spec,  C) = \{ F \in
\overline{J}_{d -
r(g-1)}(\spec) \; | \;
   \det (\pi_{s*}F) = L_{0} \}$, and $\overline{J}_{d -
r(g-1)}(\spec)$ is the
generalized
   Jacobian parametrizing rank one torsion free sheaves of Euler
characteristic
$d - r(g-1)$
   on $\spec$.}
   \end{rem}

   \bigskip

   \noindent
   \begin{rem} \label{rem15} {\rm If the curve $\spec$ is smooth,
then the
Prymian
   ${\rm Prym}_{\tilde{d}}(\spec,  C)$ is in a natural way torsor
over an
abelian variety.
   Indeed, it is easy to see that
   \[ \det (\pi_{s*}F) = {\rm Nm}_{\pi_{s*}}(F)\otimes \det
(\pi_{s*}{\cal
   O}_{\spec}),\]
   for any line bundle $F$ on $\spec$. Since $\pi_{s*}{\cal
O}_{\spec} =
   {\cal O}_{C}\oplus \omega_{C}^{-1} \oplus \ldots \oplus
\omega_{C}^{-(r-1)}$, see
   \cite{bnr}, we get
      \[ {\rm Prym}_{\tilde{d}}(\spec,  C) = \{ F \in
J^{\tilde{d}}(\spec) \; |
\;
   {\rm Nm}_{\pi_{s*}}(F) = L_{0}\otimes \omega_{C}^{\otimes
\frac{r(r-1)}{2}}
\},\]
   with $\tilde{d} = d + r(r-1)(g-1)$.}
   \end{rem}

   \bigskip

   \noindent
   \begin{rem} \label{rem16} {\rm
   We will need more explicit geometric description for the fiber
$H^{-1}(s)$
of the Hitchin
   map over a generic point $s \in {\cal D}$ of the discriminant.
Observe first
that there
   exists a point $s \in {\cal D}$ for which $\spec$ is irreducible
and has a
unique ordinary
   double point as singularity. Indeed, let $s_{r} \in
H^{0}(C,\omega_{C}^{\otimes r})$ be a
   section having a unique double zero $a \in C$. It follows then by
the
Jacobian criterion
   for smoothness that the curve
   \[ \spec : X^{r} + s_{r}Y^{r} = 0, \]
   has a unique node lying over $a$. Consequently, by upper
semicontinuity
there is a
   non-empty Zariski open set  ${\cal D}^{o} \subset {\cal D}$
satisfying
   \[ {\cal D}^{o} = \{s \in {\cal D} \; | \; \spec \; {\rm is \;
irreducible
\; and \; has \; a \; unique \; ordinary \; double \;
   point} \}. \]

   Let now $s \in {\cal D}^{o}$ and say $p \in \spec$ is its ordinary
double
point. Set
   $a = \pi_{s}(p)$. To study the geometric properties of  ${\rm
Prym}_{\tilde{d}}(\spec,
    C)$, we recall the construction of the compactified Jacobian
    $\overline{J}_{d - r(g-1)}(\spec)$ in this simple case. Let
    \[
    \begin{diagram}
    \node{\spec^{\nu}} \arrow{e,t}{\nu} \arrow{se,r}{\pi_{s}^{\nu}}
\node{\spec}
    \arrow{s,r}{\pi_{s}}  \\
    \node[2]{C}
    \end{diagram}
    \]
    be the normalization of the curve $\spec$. Let $\nu^{-1}(p) = \{
p^{+},
p^{-} \} \subset
    \spec^{\nu}$. For any rank one torsion free sheaf $F \rightarrow
\spec$ we
have the
    exact sequence
    \begin{equation} \label{eq18}
     0 \longrightarrow F \longrightarrow \nu_{*}\nu^{*}F
\longrightarrow {\Bbb
C}_{p}
    \longrightarrow 0.
    \end{equation}
    Therefore $\chi(F) + 1 = \chi(\nu_{*}\nu^{*}F) =  \chi(\nu^{*}F)$
and hence
the
    pull-back induces a morphism
    \begin{equation} \label{eq19}
    \nu^{*} :  \overline{J}_{d - r(g-1)}(\spec) \longrightarrow J_{d
+ 1 -
    r(g-1)}(\spec^{\nu}),
    \end{equation}
    where $J_{d + 1 - r(g-1)}(\spec^{\nu})$  is the component of the
Picard
group of
    $\spec^{\nu}$ consisting of line bundles of Euler characteristic
$d + 1 -
r(g -1)$. One can
    show that (\ref{eq19}) is a bundle with structure group ${\Bbb
C}^{\times}$
and with
    fibers isomorphic
    to a ${\Bbb P}^{1}$ glued at two points. To identify the bundle
(\ref{eq19})
    explicitly, consider a Poincare bundle
    \[
    \begin{diagram}[J]
    \node{{\cal P}} \arrow{s} \\
     \node{J_{d + 1 - r(g-1)}(\spec^{\nu})\times\spec^{\nu}}
     \end{diagram}
    \]
    and set ${\cal P}_{+} = {\cal P}_{|J_{d + 1 -
r(g-1)}(\spec^{\nu})\times
\{p^{+}\} }$ and
    ${\cal P}_{-} = {\cal P}_{|J_{d + 1 - r(g-1)}(\spec^{\nu})\times
\{p^{-}\}
}$. The ${\Bbb
    P}^1$-bundle
    \[ {\Bbb P}({\cal P}_{+}\oplus {\cal P}_{-}) \longrightarrow J_{d
+ 1 -
    r(g-1)}(\spec^{\nu}) \]
    does not depend on the choice of the Poincare bundle ${\cal P}$
and is
   furnished with two sections $X_{+}$ and $X_{-}$ corresponding to
${\cal
P}_{+}$ and
   ${\cal P}_{-}$ respectively. Furthermore, there is a bundle
isomorphism
   \[
   \begin{diagram}
   \node{\overline{J}_{d - r(g-1)}(\spec)} \arrow[2]{e,t}{\cong}
\arrow{se,r}{\nu^{*}}
   \node[2]{{\Bbb P}({\cal P}_{+}\oplus {\cal P}_{-})/X_{+}\sim
X_{-}}
   \arrow{sw} \\
   \node[2]{J_{d + 1 - r(g-1)}(\spec^{\nu})}
   \end{diagram}
   \]
    To describe the subvariety ${\rm Prym}_{\tilde{d}}(\spec, C)
\subset
\overline{J}_{d +1 -
    r(g-1)}(\spec)$ in these terms, we need the following lemma
    \begin{lem} \label{lem13} If $F_{1}, F_{2} \in \overline{J}_{d +
1-
r(g-1)}(\spec)$ are
    such that
    \[ \det (\pi_{s*}F_{1}) =  \det (\pi_{s*}F_{2}), \]
     then
     \[ \det (\pi_{s*}^{\nu} \nu^{*}F_{1}) =  \det (\pi_{s*}^{\nu}
\nu^{*}F_{2}). \]
     \end{lem}
     {\bf Proof.}  For any $F \in \overline{J}_{d + 1-
r(g-1)}(\spec)$ we have
the short exact
     sequence (\ref{eq18}). After taking direct images, we get
     \[
     {\divide\dgARROWLENGTH by 4
     \begin{diagram}
     \node{0} \arrow[2]{e} \node[2]{\pi_{s*}F} \arrow[2]{e}
      \node[2]{\pi_{s*}\nu_{*}\nu^{*}F}
     \arrow[2]{e} \arrow{s,=,-} \node[2]{\pi_{s*}{\Bbb C}_{p}}
\arrow[2]{e}
\arrow{s,=,-}
     \node[2]{R^{1}\pi_{s*}F} \arrow{s,=,-} \\
     \node[5]{\pi_{s*}^{\nu} \nu^{*}F} \node[2]{{\Bbb C}_{a}}
\node[2]{0}
     \end{diagram}}
     \]
     i.e. we have a short exact sequence of sheaves on $C$:
     \[ 0 \longrightarrow \pi_{s*}F \longrightarrow \pi_{s*}^{\nu}
\nu^{*}F
     \longrightarrow  {\Bbb C}_{a} \longrightarrow 0.\]
     Therefore the bundle $\pi_{s*}F$ is a Hecke transform of the
bundle
     $\pi_{s*}^{\nu} \nu^{*}F$ with center at an $(r-1)$-dimensional
subspace
of
     the fiber $(\pi_{s*}^{\nu} \nu^{*}F)_{a}$. Thus
     \[ \det(\pi_{s*}F) = \det(\pi_{s*}^{\nu} \nu^{*}F)\otimes {\cal
O}_{C}(-x), \]
     which yields the statement of the lemma.
     \begin{flushright}
     $\Box$
     \end{flushright}

     According to the above lemma, if a point $F \in \overline{J}_{d
+ 1-
r(g-1)}(\spec)$
     belongs to the subvariety $\prymd$, then every point in the
fiber
     $(\nu^{*})^{-1}(\nu^{*}(F))$ also belongs to this subvariety.

     Define ${\rm Prym}(\spec^{\nu},C) := \{ F \in J_{d + 1 -
r(g-1)}(\spec^{\nu}) \; | \;
     \det (\pi_{s*}^{\nu}F) = L_{0}(-a) \}$. Then $\prymd$ fits in
the fiber
square
     \[
     \begin{diagram}
     \node{\prymd} \arrow{e} \arrow{s,l}{\nu^{*}}
\node{\overline{J}_{d + 1-
     r(g-1)}(\spec)} \arrow{s,r}{\nu^{*}} \\
     \node{{\rm Prym}(\spec^{\nu},C)} \arrow{e} \node{J_{d + 1 -
r(g-1)}(\spec^{\nu})}
     \end{diagram}
     \]
     In particular, $\prymd$ is a bundle over the abelian variety
${\rm
     Prym}(\spec^{\nu},C)$ whose total space is singular along a
divisor and
has a smooth
     normalization which is a ${\Bbb P}^{1}$-bundle over ${\rm
Prym}(\spec^{\nu},C)$.}
     \end{rem}

     \bigskip

     The total space $\cotangent$ of the cotangent bundle of
$\sustable$ is a
Zariski
     open set in the moduli space of Higgs bundles $\higgs$. It has
the
advantage that we
     have morphisms both to $\sustable$ and $W$:
     \[
     \begin{diagram}
     \node[2]{\cotangent} \arrow{sw,l}{\Pi} \arrow{se,l}{H} \\
     \node{\sustable} \node[2]{W}
     \end{diagram}
     \]
      whose fibers are generically transversal. The natural
projection $\Pi :
\cotangent
      \rightarrow \sustable$ is onto by definition and the Hitchin
map $H$ is
dominant by
      an argument of Beauville-Narasimhan-Ramanan, see \cite{bnr}.
In
particular,
      this implies that for the generic element $s \in W$ the
push-forward map
      $\pi_{s*} : H^{-1}(s)\cap \cotangent \rightarrow \sustable$ is
dominant.
      In specific geometric situations it is important to know what
is the
codimension of the
      locus of  those $s \in W$ for which $\pi_{s*}$ is not dominant.
A general
result
      to this extend can be deduced easily from the following Lemma
\ref{lem14}
      communicated to us by E. Markman.

      First we introduce some notation.
      \begin{defi} A vector bundle $E \in \urd$ is called very stable
if it
does not have
      non-zero nilpotent $\omega_{C}$-twisted endomorphism.
      \end{defi}

      By a result of Drinfeld and Laumon, see \cite{la}, \cite{bnr},
very
stable vector bundles
      always exist. Let $U \subset \sustable$ be the non-empty
Zariski open set
       consisting of very stable bundles.
      For an $E \in \cotangent$ denote by ${\cal X}_{E} :=
T^{*}_{[E]}\sustable
=
      H^{0}(C,{\rm End}^{o}E\otimes \omega_{C})  \subset \cotangent$
the fiber
of the
       cotangent bundle at $[E]$.

      \begin{lem} \label{lem14}  For any $E \in U$ the Hitchin map
      \begin{equation} \label{eq111}
       H_{E} = H_{|{\cal X}_{E}}  :  {\cal X}_{E} \longrightarrow W
       \end{equation}
      is surjective.
      \end{lem}
      {\bf Proof.} Since $E$ is very stable we have $H_{E}^{-1}(0) =
0$ and
hence the
      dimension of the generic fiber of $H_{E}$ is zero.  On the
other hand,
$\dim {\cal
      X}_{E} = \dim W$ and therefore $H_{E} : {\cal X}_{E}
\rightarrow W$ is
      dominant.

      Denote by $W^{o} \subset W$ the Zariski open subset
      \[ W^{o} = \{ s \in W \; | \; \dim H_{E}^{-1}(s) = 0 \}. \]
      Clearly the subvariety ${\cal X}_{E}$ is preserved by the
${\Bbb
C}^{\times}$ action
      on $\higgs$. Furthermore the ${\Bbb C}^{\times}$-equivariance
of $H$
implies  that
      $W^{o}$ is a ${\Bbb C}^{\times}$-invariant subset of $W$. Since
$H_{E}^{-1}(0) = \{
      0 \}$,
      we have that
      \begin{equation}\label{eq110}
       H_{E} :  {\cal X}_{E}\setminus \{ 0 \} \longrightarrow
W\setminus \{ 0
\}
       \end{equation}
      is a ${\Bbb C}^{\times}$-equivariant morphism and thus descends
to a
morphism
      \[ \overline{H}_{E} : {\Bbb P}({\cal X}_{E}) \longrightarrow
{\Bbb
P}_{{\rm weight}} ,
      \]
      where ${\Bbb P}_{{\rm weight}} = ( W\setminus \{ 0 \})/{\Bbb
C}^{\times}$
is the
      corresponding  weighted projective space.

      The morphism $\overline{H}_{E}$ is dominant because its range
contains
the Zariski
      open set $( W^{o}\setminus \{ 0 \})/{\Bbb C}^{\times} \subset
{\Bbb
P}_{{\rm
      weight}}$. But ${\Bbb P}({\cal X}_{E})$ is a projective variety
and
therefore
      $\overline{H}_{E}({\Bbb P}({\cal X}_{E}))$ is projective and
thus
$\overline{H}_{E}$ is
      onto.

      We obtain a commutative diagram of  ${\Bbb C}^{\times}$-bundles
over
      ${\Bbb P}({\cal X}_{E})$ and ${\Bbb P}_{{\rm weight}}$
respectively
      \[
      \begin{diagram}
      \node{{\cal X}_{E}\setminus \{ 0 \}} \arrow{e,t}{H_{E}}
\arrow{s}
      \node{W\setminus \{ 0 \}} \arrow{s} \\
      \node{{\Bbb P}({\cal X}_{E})} \arrow{e,t}{\overline{H}_{E}}
      \node{{\Bbb P}_{{\rm weight}}}
      \end{diagram}
      \]
      Since the ${\Bbb C}^{\times}$ on the fibers of ${\cal
X}_{E}\setminus \{
0 \}
      \rightarrow  {\Bbb P}({\cal X}_{E})$ and $W\setminus \{ 0 \}
\rightarrow
      {\Bbb P}_{{\rm weight}}$ is simply transitive, we obtain that
the map
(\ref{eq110})
      (and consequently the map (\ref{eq111}))  is surjective.
      \begin{flushright}
      $\Box$
      \end{flushright}

      \begin{cor} \label{cor17} For any $s \in W$ the rational map
      \[
      \begin{diagram}[\pi_{s*} = \Pi_{|H^{-1}(s)} : \; \;  H^{-1}(s)]
      \node{\pi_{s*} = \Pi_{|H^{-1}(s)} : \; \;  H^{-1}(s)}
\arrow{e,..}
      \node{\su}
      \end{diagram}
      \]
      is dominant.
      \end{cor}
      {\bf Proof.} By Lemma \ref{lem14} the range $\Pi(H^{-1}(s))$
contains the
nonempty
      Zariski open set $U$.
      \begin{flushright}
      $\Box$
      \end{flushright}

       \section{The abelianization of an automorphism} \label{s2}
\subsection{The map $\varphi_{W}$} \label{ss21}
Start with an automorphism $\Phi$ of $\su$. We would like to
abelianize it,
i.e. to lift $\Phi$ to the moduli space of Higgs bundles and induce
from it an
automorphism of the family of spectral Pryms. The hope is that in
this way we
will get a simpler picture because of the fact that the automorphisms
of the
abelian varieties are pretty well
understood. As we will see in the subsequent sections, this hope
turns out to
be unsubstantiated in that simple form since the spectral Pryms are
not generic
abelian
varieties and can have non-trivial group automorphisms. This is the
reason why
we have to
work with the whole family of Pryms rather then the individual Prym
varieties
and to make essential use of its ``rigidity'' as a group scheme over
$W$.

To construct the abelianized version of $\Phi$, observe first that
the
codifferential of $\Phi$
provides a natural lift of $\Phi_{|\sustable}$ to the total space of
the
cotangent bundle
\[ d\Phi^{*} : \cotangent \longrightarrow \cotangent.\]
The automorphism $d\Phi^{*}$ commutes with the projection $\Pi :
\cotangent
\rightarrow \sustable$ by deffinition. The behaviour of $d\Phi^{*}$
with
respect to the Hitchin map is described by the following proposition.

\begin{prop} \label{prop21} There exists an automorphism $\varphi_{W}
\in
{\rm Aut}(W)$ making the following diagram commutative
\begin{equation} \label{eq21}
\begin{diagram}
\node{\cotangent} \arrow{e,t}{d\Phi^{*}} \arrow{s,l}{H}
\node{\cotangent}
\arrow{s,r}{H} \\
\node{W} \arrow{e,b}{\varphi_{W}} \node{W}
\end{diagram}
\end{equation}
\end{prop}
{\bf Proof.} Let $s \in W\setminus {\cal D}$. Then the spectral curve
$\spec$
is smooth and
(see Remark \ref{rem14}) $H^{-1}(s) = \prymd$. Consider the fiber
\[ H^{-1}_{\cal X}(s) :=  H^{-1}(s)\cap \cotangent = \{ F \in \prymd
\; | \;
\alpha_{*}F {\rm -
stable} \}. \]
The argument in Remark 5.2 of \cite{bnr} gives ${\rm
codim}(H^{-1}_{\cal X}(s),
H^{-1}(s)) \geq 2$. Therefore since $\prymd = H^{-1}(s)$ is smooth,
the regular
map $H\circ d\Phi^{*} :  H^{-1}_{\cal X}(s) \rightarrow W$ extends by
Hartogs
theorem
to a morphism $\overline{H\circ d\Phi^{*}} :  H^{-1}(s) \rightarrow
W$. But
$\prymd$
is projective and $W$ is an affine space; hence
\[ \overline{H\circ d\Phi^{*}}(\prymd ) = {\rm point \; in \;} W. \]
Denote this point by $\varphi_{W}(s)$. In this way we obtain a
rational map
$\varphi_{W}$ on $W$,  with possible singularities along ${\cal D}$
which makes
the diagram
(\ref{eq21}) commutative.

This implies that $\varphi_{W}$ is ${\Bbb C}^{\times}$-equivariant
since
$d\Phi^{*}$ is linear on the fibers of the cotangent bundle and $H$
is ${\Bbb
C}^{\times}$-equivariant.
Thus $\varphi_{W}$ descends to a rational map
\[
\begin{diagram}[\overline{\varphi}_{W} : \; \; {\Bbb P}_{{\rm
weight}}]
\node{\overline{\varphi}_{W} : \; \; {\Bbb P}_{{\rm weight}}}
\arrow{e,..}
\node{{\Bbb P}_{{\rm weight}}}
\end{diagram}.
\]
But ${\Bbb P}_{{\rm weight}}$ is a normal projective variety and
therefore
$\overline{\varphi}_{W}$ is not defined in a locus of codimension
$\geq 2$.
Finally the ${\Bbb C}^{\times}$-equivariance of $\varphi_{W}$
together with the
fact that ${\Bbb C}^{\times}$ acts transitively on the fibers of
$W\setminus \{
0 \}
\rightarrow  {\Bbb P}_{{\rm weight}}$ yield that $\varphi_{W}$ itself
is not
defined on a locus
of codimension $\geq 2$. Since $\varphi_{W}$ is a map from a vector
space to a
vector space, it extends by Hartogs to the whole $W$.
\begin{flushright}
$\Box$
\end{flushright}

Since we want to reduce the study of the automorphism $\Phi$ to the
study of a
suitable fiberwise automorphism of the family of spectral Pryms, the
best
situation for us would have been if we could say that $\varphi_{W} =
{\rm id}$.
Unfortunately, at this stage the only information we can extract from
the above
proposition is that $\varphi_{W}$ commutes with the ${\Bbb
C}^{\times}$ action
on $W$. In the next lemma we use this to determine some further
properties of
$\varphi_{W}$.

\begin{lem} \label{lem21} The map $\varphi_{W}$ preserves the
subspace
$W_{r-1,r}
\subset W$ and moreover $\phibase := {\varphi_{W}}_{| W_{r-1,r} }$ is
of the
form
$\phibase = (\phi_{r-1},\phi_{r})$ where
\[ \phi_{r-1} : H^{0}(C,\omega_{C}^{\otimes (r-1)}) \longrightarrow
H^{0}(C,\omega_{C}^{\otimes (r-1)}) \]
and
\[ \phi_{r} : H^{0}(C,\omega_{C}^{\otimes r}) \longrightarrow
H^{0}(C,\omega_{C}^{\otimes r}) \]
are linear maps.
\end{lem}
{\bf Proof.} Consider the algebra of regular functions $S := {\Bbb
C}[W]$ on
$W$.
Decomposing $S$ into isotypical components with respect to the ${\Bbb
C}^{\times}$
action, we get a new grading on it:
\[ S = \otimes_{n \geq 0} S_{n}, \]
where $S_{n}$ consists of all functions on $W$ on which a number $t
\in
{\Bbb C}^{\times}$ acts by a multiplication by $t^{n}$. More
explicitly, if we
choose
coordinates  $x_{i\scriptsize{1}}, \ldots , x_{in_{i}}$ in $W_{i} =
H^{0}(C,\omega_{C}^{\otimes i})$,
then the elements in $S_{n}$ are just polynomials in $x_{ij}$'s with
admissible
 multidegrees $\{ a_{ij}\}$ satisfying
\[ \sum_{i,j} i\cdot a_{ij} = n .\]
Any automorphism of $W$ commuting with the ${\Bbb C}^{\times}$ action
will
preserve this grading and in particular will preserve the subspaces
of the form
$W_{i}\oplus W_{i+1}\oplus \ldots \oplus W_{r}$ for any $i$.  In
paticular,
since $\varphi_{W}$ commutes with the ${\Bbb C}^{\times}$ action we
have
\[ \phibase : W_{r-1,r} \longrightarrow W_{r-1,r}.\]
 Moreover, since ${\rm g.c.d.}(r-1,r) = 1$, it follows that
$\phibase$ has the
required form.
 \begin{flushright}
 $\Box$
 \end{flushright}

 \subsection{Extension to the family of Pryms} \label{ss22}
 As we have seen in Section \ref{ss13} the fiber of the Hitchin map
$H : \higgs
  \longrightarrow W$ over every point $s \in W^{\rm reg} :=
W\setminus {\cal
D}$ is
  isomorphic to an abelian variety. Moreover, one can show, see
\cite{ni}, that
${\cal D}$ is
  exactly the locus of critical values of $H$. Define $\unprymd$ to
be the
preimage
  of $W^{\rm reg}$ under the Hitchin map. Then $\unprymd$ is a smooth
variety,
see
  \cite{ni}, and we have a smooth fibration $H : \unprymd \rightarrow
W^{\rm
reg}$.
  Our next goal is to extend the automorhism $d\Phi^{*}$ to an
automorphism of
  the family $\unprymd$. The first result in this direction is the
following
proposition

  \begin{prop} \label{prop22} The map $\varphi_{W}$ preserves the
discriminant
divisor
  ${\cal D} \subset W$.
  \end{prop}
  {\bf Proof.} Let $s \in W$. Corollary \ref{cor17} guarantees that
  $\prymd \cap \cotangent \neq \varnothing$ and hence, due to
Proposition
\ref{prop21},
  $d\Phi^{*}$ induces a birational automorphism between $H^{-1}(s)$
and
  $H^{-1}(\varphi_{W}(s))$.

  Let $s \in {\cal D}^{o}$. If we assume that $\varphi_{W}(s) \not\in
{\cal
D}$, then
  $\widetilde{C}_{\varphi_{W}(s)}$ is smooth and ${\rm
  Prym}(\widetilde{C}_{\varphi_{W}(s)},C)$ is isomorphic to an
abelian variety.
  Therefore $\prymd$ is birational to an abelian variety. On the
other hand
$\prymd$
  is birationally isomorphic to a ${\Bbb P}^{1}$ bundle over an
abelian
variety, see
  Remark \ref{rem16},  which is a contradiction because such a bundle
has
Kodaira
  dimension $- \infty$.

  Consequently $\varphi_{W}({\cal D}^{o}) \subset {\cal D}$. But
$\varphi_{W}$
is an
  automorphism and hence is a closed map which yields
$\varphi_{W}({\cal D}) =
{\cal D}$
  since the discriminant is irreducible, see Corollary \ref{cor15}.
  \begin{flushright}
  $\Box$
  \end{flushright}

  \begin{cor} \label{cor21} Let  $\phi_{r}$ be the map defined in
Lemma
\ref{lem21}.  Then
  there exists a number $\lambda \in {\Bbb C}^{\times}$ and an
automorphism
  $\sigma$ of the curve $C$ such that
  \[ \phi_{r} = m_{\lambda}\circ \sigma^{*} , \]
  where $m_{\lambda}$ is the dilation by $\lambda$.
  \end{cor}
  {\bf Proof.} By Lemma \ref{lem21} the map $\phi_{r}$ is a linear
map. On the
other
  hand, Proposition \ref{prop22} implies in particular that the map
$\phi_{r}$
  preserves the discriminant divisor ${\cal D}_{r}$. Therefore the
induced map
  $\overline{\phi}_{r}$ on the projective space ${\Bbb P}(W)$
preserves the
divisor
  ${\Bbb P}({\cal D}_{r})$. But, by Remark \ref{rem14}, the divisor
${\Bbb
P}({\cal
  D}_{r})$ is the dual variety of the $r$-canonical model of $C$ and
therefore
the dual
  map $\overline{\phi}_{r}^{\vee}$ preserves the curve $f_{{\Bbb
P}(W_{r})}(C)
\subset
  {\Bbb P}(W_{r})^{\vee}$. Let $\sigma$ be the induced automorphism
of $C$.
  Then,  since $f_{{\Bbb P}(W_{r})}(C)$ spans ${\Bbb
P}(W_{r})^{\vee}$ and
  $\overline{\phi}_{r}^{\vee}$ is linear, it follows that
$\overline{\phi}_{r}
= \sigma^{*}$.
  \begin{flushright}
  $\Box$
  \end{flushright}
  \begin{rem} \label{rem21} {\rm If the curve $C$ does not have
  automorphisms, then the above corollary implies that $\phi_{r}$ is
just a
dilation.}
  \end{rem}
  The proposition above yields the commutative diagram
  \[
  \begin{diagram}
  \node{\unprymd} \arrow{e,t,..}{d\Phi^{*}} \arrow{s,l}{H}
\node{\unprymd}
  \arrow{s,r}{H} \\
   \node{W^{\rm reg}} \arrow{e,b}{\varphi_{W}} \node{W^{\rm reg}}
\end{diagram}
\]
Furthermore, as we saw in the proof of Proposition \ref{prop22} the
map
$d\Phi^{*}$
gives a birational isomorphism between  $H^{-1}(s)$ and
$H^{-1}(\varphi_{W}(s))$
for any $s$. In particular, in the case $s \in  W^{\rm reg}$ we get a
birational automorphism between abelian varieties and hence
$d\Phi^{*}_{|H^{-1}(s)}$ extends to a
biregular isomorphism between  $H^{-1}(s)$ and
$H^{-1}(\varphi_{W}(s))$.
Therefore
$d\Phi^{*}$ extends as a continuous automorphism $\phitilded$ of the
whole
variety $\unprymd$ and since $\unprymd$ is smooth, we get:

\begin{prop} \label{prop23} There exists a biregular automorphism
$\phitilded$ extending $d\Phi^{*}$ which fits in the commutative
diagram
 \[
  \begin{diagram}
  \node{\unprymd} \arrow{e,t}{\phitilded } \arrow{s,l}{H}
\node{\unprymd}
  \arrow{s,r}{H} \\
   \node{W^{\rm reg}} \arrow{e,b}{\varphi_{W}} \node{W^{\rm reg}}
\end{diagram}
\]
 \end{prop}
 \section{The fibration $\unspec \lra \wbase$ }
 \label{section3 }
  \subsection{The Picard group of $\unspec $ }
 \label{picunspec}
 We start with some notation: For a vector space $V$, the points of
the
Grassmanian $Gr(k,V)$  correspond to  codimension $k$ linear
subspaces of $V$.
Also, for a linear space of sections $W$ on $S$ and a given point $p$
on $S$,
$W_p$ denotes the linear subspace of sections of $W$ vanishing at
$p$. We now
make the assumption that $r \geq 3$.
 We are going  first to prove our main Theorem \ref{theor1}   in the
case $r
\geq 3$ and next, in Section \ref{rank2}, we outline the
modifications needed
for the proof of the rank $2$ case.
   For simplicity  denote the space $W_{r-1,r} $ by $B$ and  the
spaces
$\overline{W}_{r-1,r}$  and  $\wbase = W_{r-1,r} \setminus {\cal D}$
by
$\overline{B}$ and $B^{reg } $ respectively.  In the following
we are going to use repeatedly  Proposition 6.5 in \cite{ha}, which
compares
the
 Picard group  of a variety $X$ with that of a Zariski open subset
$U$ of $X$.
\\
\\
{\bf Notation.}  From now on  we will  denote by $\unspec $ the part
of the
universal spectral curve sitting over $\breg$,  unless otherwise
stated. \\
\\
We have the diagram:
\vskip.1in
\[
\begin{diagram}
 \node{ \unspec }  \arrow{se,l}{\pi } \arrow{s,l}{\beta}
\arrow{e,t}{h}
\node{\breg }  \\
   \node{S^{\circ  } }  \arrow{e,t}{\alpha } \node{C}
\end{diagram}
\]
 \vskip.1in

\begin{prop}
 The Picard group of the variety $\unspec $ is isomorphic to the pull
back of
the  Picard group of  the  base  curve $C$  by the map $\pi $,  i.e.
$\pic
\unspec \simeq  \pi ^*  \pic C $.
\label{proppicunspec}
 \end{prop}

\noindent
{\bf Proof.}
The variety $\unspec $ can be constructed as follows. Consider the
diagram:

 \vskip.1in
\begin{equation}
\begin{diagram}[aaaaaaaa]
 \node{f^{*}_{| \overline B |} U^{\rm reg}|_{S^{\circ}} \simeq
\unspec   }
\arrow{e,t}{ } \arrow{s,l} {\beta  }
 \node{ U^{\rm reg} \subset U   } \arrow{e}    \arrow{s,r} {  }
\node{Gr(1,\overline{B}) \times \overline {B}}  \arrow{sw,l}{p_1}
\arrow{s,r}{p_2}   \\
\node{\;\;\;\;\;\;\;  S^{\circ} \subset S } \arrow{e,b}{ f_{|
\overline {B} |}}
\node{Gr(1,\overline{B}) }
\node{\;\;\;\;\;\;\;\;\;\;\;\;\;\;\; \overline{B} \supseteq  B
\supseteq  {\cal
D} }
\end{diagram}
\label{diagunspec}
\end{equation}
 \vskip.1in

\noindent
The map  $ f_{| \overline {B} |}$ is the map associated to the base
point free
linear system $| \overline {B} |$ i.e.  it sends a point $p$ to the
point
$[\overline {B}_p] $ representing the linear subspace $\overline
{B}_p$. The
bundle $U$ is the universal bundle over the Grassmanian.   The fiber
of the
bundle $U \cap  p_2^{-1}(B) $ over a point  $[H]$  in the Grassmanian
is
isomorphic to $H \cap B$. The variety  $U \cap p_2^{-1}(B) $ is an
affine
bundle over
  $ f_{| \overline {B} |}(S^{\circ} )$  with fiber over the point
$f_{|
\overline {B} |}(p)$ isomorphic to $B _p$.
Define $U^{\rm reg}\stackrel{\rm df}{=}  U \cap p_2^{-1}(B \setminus
{\cal D} )
=U \cap p_2^{-1}
 (B ^{\rm reg})$. Then $    \unspec   \simeq    f^{*}_{| \overline B
|}
    U^{\rm reg}|_{S^{\circ}}  $.  Note that by Corollary \ref{cor16},
the
divisor
  $ f_{| \overline {B} |}^* ( p_2 ^{-1} (\Delta ) \cap U)$ is
irreducible.
Also it is linearly equivalent  to zero  since $\pic B $ is trivial.

\noindent
To complete the proof of the proposition, we have:
 \[ \begin{array}{rlll}
 \pic \unspec  & \simeq &  \pic  f^{*}_{|\overline{B }|} U^{\rm
reg}|_{\so }  &
 \\
     & \simeq     &  \pic  f^{*}_{|\overline {B }|}  (  U \cap
p_2^{-1}(B
\setminus
   {\cal D} )  ) |_{S^{\circ}}    &      \\
       & \simeq     &  \pic  f^{*}_{|\overline {B }|}  (  U \cap
p_2^{-1}(B  )
) |_{\so }    & \mbox{since}\;\; U \cap  p_2^{-1} {\cal D} \;\; \mbox
{is  an
irreducible divisor lin. equiv. to} \; \; 0,
           \\
            & \simeq  &    \beta ^* \pic  S^{\circ}  &  \mbox {since}
\;\;
  f^{*}_{|\overline {B }|}  (  U \cap p_2^{-1}(B  ) ) |_{\so }\;\;
    \mbox{is an affine bundle over} \;\; \so  , \\
           &  \simeq   &  \pi ^* \pic C  & \mbox{since} \;\; \so
\;\;
\mbox{is an affine bundle over} \;\; C  .
\end{array} \]
\begin{flushright}
   $\Box$
   \end{flushright}

\subsection{The Picard group of $\unspec _p  $}
\label{picunspecp}
For a fixed point $p$ in $S^{\circ }$ we denote by $\unspec _p$ the
subvariety
of $\unspec $ consisting of those curves whose image in $S^{\circ}$
passes
through $p$.
  $\unspec _p$ sits over the subspace $\breg _p $ of $\breg $,
consisting of
those  sections vanishing at $p$.  Let  $h: \unspec _p  \lra   \bregp
$ denote
again the restriction of the  map  $h$  on $\unspec _p$.  We are
going to
calculate the Picard group of $ \unspec _p  $. Observe  that the
above
fibration $h$  has a section corresponding to the point $p$. Take the
restriction of the map $\beta $ to $\unspec  _p$ i.e. $\beta: \unspec
_p
\lra S^{\circ}$. Let $H_{p}$ be the preimage of $p$ via the map
$\beta $.   By
a construction  similar to that for the variety  $\unspec $ in the
above
Section \ref{picunspec}, one can see that the fibration  $\beta:
\unspec _p
\setminus H_p  \lra \so \setminus \{ p \}  $ is  the complement of an
irreducible divisor linearly equivalent to zero in an affine
fibration.
 We thus have  $\pic
  (\unspec  _p\setminus H_p)   \simeq \beta ^*  \pic ( \so  \setminus
\{ p \})
\simeq \beta ^*  \pic \so  \simeq   \pi ^* \pic C$.  Since  $H_p$ is
an
irreducible divisor in $\unspec _p$,  we conclude that  $\pic
\unspec _ p $
is generated by $\pi ^* \pic C$ and $H_p$.  We proceed by showing
that  we
actually  have:

\begin{prop}
    $  \pic \unspec _p   \simeq  \pi ^* \pic C   \oplus {\Bbb Z}[H_p]
$.
\label{proppicunspecp}
\end{prop}

\noindent
{\bf Proof.}    Pick a generic pencil $\pr $ in $\overline {B}_p$
such that all
the fibers of the restriction of the compactified universal  spectral
curve
over the pencil  are irreducible,  except the one over the infinity
point $b$
which has exactly two irreducible components, see Corollary
\ref{cor14}. We may
also assume that none of the singular points  of the fibers lies  on
the $N=
2r^2(g-1)$ base points of the system. The restriction of the
compactified
universal spectral curve  over  $\pr $,  is a smooth surface $X$
which is the
blow up of the surface $S$ over the $N$ base points $p_1=p, \ldots
,p_N$  of
the pencil. We have the following picture:

\vspace{4in}

  \noindent
  In the above picture, $b$ is the point at infinity,  $b_1, \ldots
,b_k$
correspond to the
 singular fibers and  $E_1, \ldots  ,E_N$ are the exceptional
divisors. The
curve over $b$ consist of two components $\tilde {Y}_{\infty }$ and
$\tilde{C}_b^{'}$,  one of which, $\tilde {Y}_{\infty }$, is the
proper
transform of the divisor at infinity on $S$. The intersection of $X$
with the
variety $\unspec _p $ is $X^{\circ } = X \setminus h^{-1}(b,b_1,
\ldots ,b_k)$.
We define $E_i^{\circ} =  E_i \cap X^{\circ}$. Note that $E_1^{\circ}
= H_p
\cap X$.    To prove the claim, it is enough to show that on
$X^{\circ }$ the
line bundle $E^{\circ} $ is independent from $\pi ^* \pic C$.
We have
\begin{equation}
  \pic X \simeq  \beta ^* \pic S \oplus {\Bbb Z}[E_i]_{i=1, \ldots ,
N}
\simeq  \pi ^* \pic C \oplus {\Bbb Z}[\tilde {Y}_{\infty }]
\oplus_{i=1} ^N
(\oplus  {\Bbb Z}[E_i]).
\label{eqpicunspecp1}
\end{equation}
 Now $\beta ^* (rY_{\infty} + r \alpha ^* \omega _C ) = h
^{-1}(\mbox{point}) +
\sum _{i=1}^NE_i $, where the equality stands for the  linear
equivalence of
divisors.  Therefore, $r \tilde {Y}_{\infty } + r \pi ^* \omega _C =
h
^{-1}(\mbox{point}) + \sum _{i=1}^NE_i = \tilde {C} ^{'}_b  + \tilde
{Y}_{\infty }+ \sum _{i=1}^NE_i$.   In other words,
\begin{equation}
\tilde {C} ^{'}_b= (r-1)  \tilde{Y}_{\infty} + r \pi ^* \omega _C -
\sum _{i=1}^NE_i  .
\label{eqpicunspecp2}
\end{equation}
  All the fibers of $h$ define linear equivalent divisors on $X$ and
the
restriction
 on $X \setminus  h ^{-1} (b)$  of the line bundles corresponding to
the
divisors
 $\tilde {C} ^{'}_b$ and $\tilde {Y}_{\infty }$ are trivial.  We thus
get
\[  \begin{array}{rlll}
\pic  X ^{\circ }  & \simeq &  \pic X \setminus  h ^{-1}(b)  \simeq
\pic X
\left/  (\tilde {Y}_{\infty } =0, \tilde {C} ^{'}_b =0 )     \right.
& \\
    &  \simeq   &   \pi ^* \pic C \oplus_{i=1} ^N (\oplus  {\Bbb
Z}[E^{\circ
}_i])  \left/ (r\pi ^* \omega _C = \sum _{i=1}^N E^{\circ} _i )
\right.  &
\mbox{by}  \;\;
   (\ref{eqpicunspecp1})  \;\; \mbox{and} \;\;
(\ref{eqpicunspecp2}).
\end{array} \]
Assume now that $\pi ^* L + m H_p = 0 $ on $\unspec _p$. Then $\pi ^*
L +
 m E^{\circ }_1 = 0$ on $X^{\circ}$. By the description of the $\pic
X ^{\circ
} $ we conclude that $\pi ^* L =0$ and $m=0$  and this completes the
proof of
the proposition.
  \begin{flushright}
   $\Box$
   \end{flushright}

 \subsection{ The sections of  $ H: \unjacd   \lra \wbase$ }
\label{sections}

\begin{prop}
 The only sections of the map  $ H: \unjacd   \lra \wbase =  \breg$
are
those coming from a pull back of a fixed line bundle on $C$.  In
other words,
if $\sigma : \breg \lra \unjacd $ is a section of $H$, then $\sigma
(s) = [\pi
_s ^* M] $ where $M$ is a fixed line bundle on $C$.
In particular,  if $r$ does not divide $\tilde {d}$, then the map $H$
has no
sections.
\label{propsections}
\end{prop}
\noindent
We start with two lemmas:
  \begin{lem} On $\unjacd \times _{\breg} \unspec $  there exists a
line bundle
$\poincd   $
such that $\poincd |_{[L] \times \spec }  \simeq L^{\otimes n}$ for
some
integer $n$.
\label{lempoincare}
\end{lem}
{\bf Proof.}  See \cite{mr}.
 \begin{flushright}
   $\Box$
   \end{flushright}

\begin{rem}
{ \rm One can actually prove that the minimum such positive integer
$n$ is
equal to ${\rm g.c.d.(r,d)}$.  }
\label{remsections1}
\end{rem}

\begin{lem}
  If $\sigma $ is a section of the map $H$ whose image lies in the
locus $\pi
^* \pic C$ i.e. $\sigma (s)=[\pi ^* M_s] $, where $M_s$ is a line
bundle on
$C$, then $M_s=M$ for all $s \in \breg $.
\label{lemsectionimage}
\end{lem}
\noindent
{\bf Proof.}   Consider the map $\gamma:  \unspec \lra  \unjacd
\times _{\breg
} \unspec $ which sends a point $q$ sitting over $s \in \breg $ to $(
\sigma
(s), q )$. Then $
[\gamma ^* \poincd |_{[L]} ]  =  n \sigma (s)$. By the description of
$\pic
\unspec $, see
  Proposition   \ref{proppicunspec},  we get that
  $\gamma ^* \poincd \simeq \pi ^* M_1$ for a fixed  $M_1$ in $\pic
C$. Hence,
  $n [\pi ^* M_s] = n \sigma (s) = [\pi ^* M_1]$ for all $s\in \breg
$. Since
the map $\pi ^* $ is one to one, see  Remark 3.10 in  \cite{bnr}, we
conclude
that the map $\breg  \lra \pic C$ that sends $s$ to $[M_s]$ has
finite image.
Since $\breg $ is connected,   the map is constant.
  \begin{flushright}
   $\Box$
   \end{flushright}

   \bigskip

\noindent
{\bf Proof of Proposition \ref{propsections}.} Say that $\sigma $ is
not coming
from a pull back. By the above Lemma \ref{lemsectionimage}, we may
assume that
there exists an $s_0 \in \breg $ such that $\sigma (s_0) $ is not of
the form
$\pi ^* A$ for some $A$ in $\pic C$. Take a point $p$ on $\so $ that
lies on
$\tilde {C}_{s_0} \subseteq \so $. The family of curves $\unspec _p $
has a
section and therefore on  $J^{\tilde {d}} (\unspec _p) \times _{\breg
_p}
\unspec _p $ there exists a Poincare bundle $\poincd _p $, see
\cite{mr}. We
have
 \[
[\gamma ^* \poincd  |_{\spec }] =n \sigma (s) \;\; \mbox{and}\;\;
[\gamma ^*
\poincd _p|_{\spec }] = \sigma (s) \;\; \mbox{for all} \;\; s \in
\breg _p.
\]
Therefore,
  $n   \gamma ^* \poincd _p|_{\spec } \simeq \gamma ^* \poincd
|_{\spec }$ for
all
   $s\in \breg _p$.  Since $\pic \breg _p $ is trivial, we conclude
by the
see-saw principle, see \cite{mu},
 that
 \begin{equation}
     n  \gamma ^* \poincd  _p  \simeq  \gamma ^* \poincd|_{\unspec _p
}   \;\;
\mbox{on}\;\;
     \unspec _p.
\label{eqpicunspecp3}
\end{equation}

 On the other hand,  by  the description of the Picard groups of
$\unspec _p $
   and   $\unspec $, see Propositions \ref{proppicunspec}  and
\ref{proppicunspecp},
 we have
 \begin{equation}
   \gamma ^* \poincd  _p \simeq \pi ^* M + m H_p  \;\; \mbox{and}
\;\;  \gamma
^* \poincd
 \simeq \pi ^* L.
\label{eqpicunspecp4}
\end{equation}
 Hence,  by (\ref{eqpicunspecp3}) and (\ref{eqpicunspecp4}), we have
on
$\unspec _p $
  that $n \pi ^* M + n m H_p \simeq \pi ^* L $. Proposition
\ref{proppicunspecp} implies
 that $m=0$ i.e. $\gamma ^* \poincd  _p \simeq \pi ^* M $. But then,
$\sigma
(s_0) = [\gamma ^* \poincd  _p |_{\tilde {C}_{s_0}}] = [ \pi ^* M]$
which
contradicts the assumption on $\sigma (s_0) $.
  \begin{flushright}
   $\Box$
   \end{flushright}

 \section{The ring of correspondences}
\label{section4}
\subsection{The map  $ \beta ^{'} : \fibprod \lra \so \times \so $ }
\label{soso}
 Throughout the Sections \ref{section4} and \ref{section5}, we will
denote by
$\phi $
  the map  $\phi _{r-1,r} : B \lra B$.  According to Lemma
\ref{lem21}, we have that  $\phi _{r-1,r}$ has the form $\phi
_{r-1,r}=(\phi
_{r-1}, \phi _r)$ where $\phi _i : H^0(C, \omega _C^i) \lra  H^0(C,
\omega
_C^i)$ are linear automorphisms  for $i=r-1, r$.  Furthermore, by
Proposition
\ref{prop22},  we have that $\phi ({\cal D}) = {\cal D}$. Therefore,
the
restriction of the map $\phi $ on $\breg $,
 which we
  will denote again by $\phi $, induces an automorphism $\phi :
\breg \lra
\breg $.
 In this section we study   the fiber  product $\fibprod $   defined
by the
diagram

\vskip.1in
\[
 {\divide\dgARROWLENGTH by 8
\begin{diagram}
 \node{\so \times \so} \node[3]{\fibprod }\arrow[3]{w,t}{\beta ^{'}}
\arrow[2]{e,t}{ }
\arrow[2]{s,l}{ }
 \node[2]{\unspec   }    \arrow{s,r} {h }\\
\node{  }  \node[3]{ } \node[2]{ \breg  } \arrow{s,r}{\phi} \\
\node{ } \node[3]{\unspec} \arrow[2]{e,t}{h}  \node[2]{\breg}
    \end{diagram}}
\]
\vskip.1in

\noindent
  We define $\overline {\phi }  : \overline{B} \lra \overline{B}$ to
be the map
which extends $\phi$  to $\overline{B}$
  as  $ \overline {\phi }=   (1,\phi _{r-1}, \phi _r) $.
 To investigate the map $\beta ^{'}$ , we define  in the product $S
\times S$
the  following
 loci $A$ and  $\Gamma $.

\begin{defi}  {\rm
 \[ \begin{array}{rll}
A   &  \stackrel{\rm df }{=} & \so \times \so \setminus {\rm
Im}\beta^{'}   \\
  \Gamma &  \stackrel{\rm df }{=} &  \{ (p,q) \in  S \times S  \;\;
\mbox{such
that} \;\;  \overline{B}_p =\overline {\phi } (\overline{B}_q)   \}
   \end{array} \]  }
\end{defi}

\noindent
  On $S \times S \setminus (A \cup \Gamma) $ we can define a map $f$
to the
 Grassmanian  $Gr(2,\overline{B})$ of codimension $2$ linear
subspaces of
 $\overline{B}$, by sending the point $(p,q)$ to the class of the
plane
 $\overline{B}_p  \cap \overline {\phi }(\overline{B}_q) $. We have
the
following diagram:
 \vskip.1in

\begin{equation}
\begin{diagram}[aaaaaaaaaaaaaaaaaaaaa]
 \node{\;\;\;\;\;\; \fibprod |_{\so \times \so \setminus \Gamma}
\simeq  f^{*}
U^{\rm reg}|_{\so \times \so
\setminus \Gamma } }  \arrow{e,t}{ } \arrow{s,l} {\beta ^{'} }
 \node{ U^{\rm reg} \subset U   } \arrow{e}   \arrow{s,r} {  }
\node{Gr(2,\overline{B}) \times \overline {B}}  \arrow{sw,l}{p_1}
\arrow{s,r}{p_2}   \\
\node{\;\;\;\;\;\;\;\;\;\;\;\;\;\;\;\;\;\;  \so \times \so \setminus
(A \cup
\Gamma)    \subseteq S \times S \setminus (A \cup \Gamma)   }
\arrow{e,t}{ f}
\node{Gr(2,\overline{B}) }
\node{\;\;\;\;\;\;\;\;\;\;\;\;\;\;\; \overline{B} \supseteq  B
\supseteq  {\cal
D} }
    \end{diagram}
\label{diagfibprod2}
\end{equation}
 \vskip.1in

\noindent
 The notation is in complete analogy with that of the Diagram
(\ref{diagunspec}).
 As in the  case of $\unspec $,  the variety $ \fibprod |_{\so \times
\so
\setminus (A \cup \Gamma)} $ is an affine bundle over the image of
$\so \times
\so \setminus (A \cup \Gamma) $ under the map $f$.
\begin{lem}
 ${\rm dim}A \leq 2$
\label{lemdimA}
\end{lem}
\noindent
{\bf Proof.}  Let $B^{\infty }$ denote the space of sections
$W_{r-1,r}^{\infty}$, i.e. the linear subspace of $\overline{W}$
consisting of
points of the form $(0,0, \ldots , 0, s_{r-1}, s_r)$.  We
 have
\[   \begin{array}{ll}
   A  & = \{ (p,q) \in  \so \times \so   \;\; \mbox{such that} \;\;
U|_{f(p,q)}
\cap p_2^{-1} B =\varnothing  \} \\
    &   = \{ (p,q) \in  \so  \times \so   \;\; \mbox{such that} \;\;
B_p
  \cap  \phi (B _q)  = \varnothing  \}  \\
   &    = \{ (p,q) \in  \so  \times \so   \;\; \mbox{such that} \;\;
 B^{\infty}_p   =  \overline {\phi}  (B^{\infty} _q)   \}
\end{array}  \]
\noindent
By  Corollary \ref{cor13}, the linear system defined by $B^{\infty }
$
separates points on $\so $. Therefore,  given a point  $p \in \so $
there
exists at most one $q \in \so $ with $(p,q) \in A$ and this proves
the lemma.
 \begin{flushright}
   $\Box$
   \end{flushright}

   \bigskip

We  give now a better description of the locus $\Gamma $. Consider
the diagram

\vskip.1in
\[
\begin{diagram}
 \node{S  } \arrow{e,t}{^{f_{| \overline{B}| }}}
\arrow{se,r}{\overline{\phi}
 ^{\, *} \, f_{| \overline{B}| } }
 \node{{\Bbb P}(  \overline{B}^{\, \vee}) }    \arrow{s,r} {\overline
{\phi}
^{\, *}  }\\
\node { }\node { {\Bbb P}(  \overline{B}^{\, \vee}) }\end{diagram}
\]
\vskip.1in

\noindent
In the above diagram the map $\overline{\phi}  ^{\, *} $ is the
induced
automorphism of
${\Bbb P}(  \overline{B}^{\, \vee})$ by the linear automorphism
$\overline{\phi} $.  We denote by $\pi _1 $ and $\pi _2 $ the two
projections
of $\Gamma $ to the surface $S$. It is easy to see that
\[   \begin{array}{rll}
 &   \pi _1 (\Gamma)  & = \{ p  \in S   \;\; \mbox{such that} \;\;
f_{|
\overline{B}| } (p) \in
  \overline{\phi}  ^{\, *}   \, f_{| \overline{B}| } (S) \}  \\
    &   &  =  f_{| \overline{B}| }^{-1}  \left( f_{| \overline{B}|
}(S) \cap
 \overline{\phi}  ^{\, *}  \, f_{| \overline{B}| } (S) \right)   \\
  \mbox{and similar} \;\;\;\;\;\;\;\;\;\;   \\
  &    \pi _2(\Gamma ) & =  (\overline{\phi}  ^{\, *}   \, f_{|
\overline{B}|
})^{-1}  \left( f_{| \overline{B}| }(S) \cap  \overline{\phi}  ^{\,
*}   \,
f_{| \overline{B}| } (S)  \right)
  \end{array}  \]

We describe the fibers of the map $\pi _1$.  If $p \in Y_{\infty }$,
then it is
easy to see that $\pi _1 ^{-1} (p)$ consists of the whole divisor  $
Y_{\infty}
$. If $p \in \so $, then
  $\pi _1 ^{-1} (p)$ consists of at most one point since the system
$\overline{B}$
 separates points on $\so $.  By the above discussion we conclude
that

\begin{lem}
 ${\rm dim} \left( \Gamma \cap (\so \times \so )  \right) \leq 2$.
\label{lemdimB}
\end{lem}

\subsection{The Picard group of  $\fibprod  $ }
\label{picfibprod}
 We are going to use the following diagram

\vskip.1in
\[
\begin{diagram}
 \node{\fibprod  } \arrow{se,l}{\pi ^{'}}   \arrow{s,l}{\beta {'}} \\
 \node{\so \times \so  }    \arrow{e,t} {\alpha ^{'}  }  \node{C
\times C}
\end{diagram}
\]
\vskip.1in
\noindent
to calculate the Picard group of  $\fibprod  $.
\begin{lem}
  $\pic  ( \fibprod |_{\so \times \so \setminus \Gamma} )   \simeq
\pi ^{' *}
\pic
 (C \times C )$.
\label{lempicfibprod1}
\end{lem}

\noindent
{\bf Proof.}  We have
 $\fibprod |_{\so \times \so \setminus \Gamma} \simeq
f^{*} U^{\rm reg}|_{\so \times \so   \setminus \Gamma }  \simeq f^{*}
U^{\rm
reg}|_{\so \times \so   \setminus (\Gamma  \cup A )} $ since the
fiber over $A$
is empty. Therefore,

\noindent
$  \pic  ( \fibprod |_{\so \times \so \setminus \Gamma}  )
\simeq    $
\nopagebreak
\vskip-.1in
\[ \begin{array}{lll}
      \;\;\; \simeq  &    \pic  f^{*} U^{\rm reg}|_{\so \times \so
\setminus
(\Gamma  \cup A )}  &  \\
     \;\;\; \simeq     &  \pic  f^* (U \cap p_2^{-1}(B))|_{\so \times
\so
  \setminus (\Gamma  \cup A)}     &  \mbox{since the preimage of}
\;\, {\cal D}
 \;\, \mbox{is irred. divisor lin. equiv. to}\;\; 0,      \\
      \;\;\; \simeq     &  \beta ^{' *} \pic (  \so \times \so
\setminus
(\Gamma \cup A))   & \mbox{since it is an affine bundle over} \;\;
\so \times
\so \setminus (\Gamma \cup A), \\
         \;\;\; \simeq  &  \beta ^{' *} \pic ( \so \times \so  )
&  \mbox
{since} \;\; {\rm codim}(\Gamma \cup A) \geq 2 ,   \\
          \;\;\;  \simeq   & \pic  (C \times C )   & \mbox{since}
\;\; \so
\times  \so \;\; \mbox{is a rank two  affine bundle over} \;\; C
\times C.
\end{array} \]
 \begin{flushright}
   $\Box$
   \end{flushright}

   \bigskip

 The variety $\fibprod $ splits as
$$\fibprod = \fibprod |_{\so \times \so \setminus \Gamma}  \, \amalg
\,
 \fibprod |_{\Gamma}. $$
\noindent
 To calculate its  Picard group, we  have  to consider the following
two
cases.\\
\\
{\bf  Case A:}   dim$\, \Gamma \cap (\so \times \so ) \leq 1$:   Then
one can
easily see that dim$\, \fibprod |_{\Gamma}  \leq  {\rm dim}(\fibprod
) -2$.  We
thus have:
\begin{lem}
 In case A, the
$\pic ( \fibprod )  \simeq \pic ( \fibprod|_{\so \times \so \setminus
\Gamma} )
 \simeq
 \pi ^{' *}  \pic(C \times C) $.
\label{lemcasea}
\end{lem}
\noindent
{\bf Case B:}   dim$\, \Gamma \cap (\so \times \so ) =2 $:  Following
the
notation we used in the description of $\Gamma $ in the above Section
\ref{soso}, we have  $\pi _1(\Gamma ) = S$.  This implies that
$\overline
{\phi} ^{\, *} $  induces an automorphism of  $f_{|\overline{B}|}
(S)$.  The
following summarizes the basic properties of the induced
 automorphism $\overline {\phi} ^{\, *} $.
\begin{enumerate}
\item   The point $p_{\infty}=f_{|\overline{B}|}(Y_{\infty})$
remains fixed.
\item   $\overline {\phi} ^{\, *}$ sends a fiber of the map $\alpha $
to a
fiber: otherwise we get a $\pr $ cover of the curve $C$, which
contradicts the
assumption  that genus $g(C) >0$.
\item   Since $C$ has no automorphisms,  $\overline {\phi}  ^{\, *}$
sends a
fiber of the map $\alpha $  to the same fiber.
\end{enumerate}
Therefore $\overline {\phi} ^{\, *} $  induces an automorphism  $\chi
: \so
\lra \so $
 which  preserves the fibers of the map $\alpha: \so \lra C$.

\begin{lem}
Let $s\in B$ and let $\spec $ denote also  the image of the
corresponding
spectral curve on
$\so $.   Then $\chi (\spec ) = \tilde{C}_{\phi ^{-1}(s)}$.
\label{lemsurfautom}
\end{lem}

\noindent
{\bf Proof.}  Let $H_s$ denote the hyperplane in $ {\Bbb P}(
\overline{B}^{\,
\vee})$  corresponding  to the section $s$.  Then
$$    \spec   \simeq  f_{|\overline{B}|}(\spec )  = H_s \cap
f_{|\overline{B}|}(\so ) $$
and
$$  \tilde{C}_{\phi ^{-1}(s)} \simeq
f_{|\overline{B}|}(\tilde{C}_{\phi
^{-1}(s)})  =
H_{\phi  ^{-1}(s)}  \cap  f_{|\overline{B}|}(\so ). $$
It is $\overline{\phi }^{\, *}  (H_s)  = H_{\phi ^{-1}(s)}$ and so,

\[ \begin{array}{llll}
    f_{|\overline{B}| }(\tilde{C}_{\phi ^{-1}(s)})    &    =   &
\overline{\phi }^{\, *}  (H_s)  \cap    f_{|\overline{B}| }( \so ) &
\\
       & =  &   \overline{\phi }^{\, *}  (H_s)  \cap  \overline{\phi
}^{\, *}
  f_{|\overline{B}| }
   ( \so )      & \mbox{since}  \;\;    \overline{\phi }^{\, *}
f_{|\overline{B}| } ( \so )  =
     f_{|\overline{B}| } ( \so )  \;\;  \mbox{in} \;\; {\Bbb
P}(\overline{B}^*),   \\
          & =    &    \overline{\phi }^{\, *} (H_s \cap
f_{|\overline{B}| } (
\so )  )   &  \mbox{since}  \;\;  \overline{\phi }^{\, *}  \;\;
\mbox{is an
automorphism,}           \\
         &  =   &   \overline{\phi }^{\, *}    f_{|\overline{B}| }
(\spec ).
  &      \end{array} \]
 \begin{flushright}
   $\Box$
   \end{flushright}
   \bigskip

By the above Lemma  \ref{lemsurfautom},  the automorphism $\chi $ on
$\so $
induces an automorphism $\psi $ of $\unspec $ over $B$ which  makes
the
following diagram commutative
\vskip.1in
\[
 \begin{diagram}
 \node{\unspec } \arrow{s,l}{h}  \arrow{e,t}{\psi }
 \node{\unspec }   \arrow{s,r} {h} \\
 \node{B  } \arrow{e,t}{^{\phi ^{-1}}} \node{B }
 \end{diagram}
\]
 \vskip.1in

\noindent
Since $\so $ is the total space of the line bundle $\omega _C $ on
$C$,  the
map $\chi $  acts by   a  dilation  and  a translation by a section
of $\omega
_C$ on the fibers i.e. it has the form $\chi = m_{\lambda } \circ
T_s$ where
$\lambda \in {\Bbb C}^*$ and $s \in  H^0(C, \omega _C) $.   We
actually claim
that $s=0$ and that $\phi ^{-1} $ has the form
$\phi ^{-1}(s_{r-1},  s_r)= (\lambda ^{-(r-1)}  s_{r-1}, \lambda
^{-r} s_r)$.
Indeed,
 \begin{equation}
\chi  \circ h^{-1}(s_{r-1},  s_r)= h^{-1}( \phi ^{-1}(s_{r-1},
s_r)).
\label{eqpicfibprod1}
\end{equation}
  To prove the first claim, we apply  (\ref{eqpicfibprod1}) to
$(s_{r-1},
s_r)=(0,0)$. Then,
  $\chi  \circ h^{-1}(0,0)= h^{-1}(0,0)$ since  $\phi ^{-1}$ is a
linear map.
But $h^{-1}(0,0)$ is the curve $x^r=0$, and so, $\chi \circ
h^{-1}(0,0)$  is
the curve $(\lambda x +s)^r=0$, see beginning of Section \ref{ss11}.
This
implies that $s=0$.

For the second claim:
$h^{-1}(s_{r-1},  s_r)$ is the curve $x^r +s_{r-1}x +s_r=0$ and so,
$\chi \circ  h^{-1}(s_{r-1},  s_r)$  is the curve $\lambda ^r x^r
+\lambda
s_{r-1}x +s_r=0$
i.e. the curve $x^r +\lambda ^{-(r-1)} s_{r-1}x + \lambda^{-r} s_r=0$
i.e. the
curve $h^{-1}(\lambda ^{-(r-1)} s_{r-1}, \lambda^{-r} s_r)$ which
completes the
proof.

\begin{cor}
$\phi  (s_{r-1},  s_r)=(\lambda ^{r-1} s_{r-1}, \lambda^{r} s_r). $
\label{corformphi}
\end{cor}

  We now proceed with our discussion about the Picard group of
$\fibprod $
 in the case of  dim$\Gamma \cap (\so \times \so ) =2 $.  We have
that
 $\fibprod |_{\Gamma}= \Delta _{\lambda}= \{   (p, \psi ^{-1}(p)), p
\in
\unspec   \}$.  Note that if $\lambda=1$,  then  $\Delta _{\lambda} $
is
exactly the diagonal in the fiber
 product
$\tilde {\cal C \;} _{ h}  \!\!\times _{h}  \tilde {\cal C}$.
Combining the
latter with Lemma
\ref{lempicfibprod1} we get:

\begin{prop}
In case B, the  $\pic ( \fibprod ) $ is generated by  $\pi ^{' *}
\pic (C
\times C)$ and $ \Delta _{\lambda}$.
\label{proppicfibprodb}
\end{prop}

\begin{rem}
 { \rm  It is easy to see that   $\pic (\fibprod ) \simeq  \pi ^{' *}
\pic (C
\times C) \oplus
{\Bbb Z}[ \Delta _{\lambda}] $. }
\label{rempicfibprod1}
\end{rem}
\section{Proof of the main Theorem}
\label{section5}
\subsection{The form of $\tilde{\Phi }$}
\label{formphi}
We start with a definition.  Let $C$ and $C_1$ are two smooth curves.
Given a
line bundle ${\cal L}$ on the product $C \times C_1 $, then ${\cal
L}$ induces
a map
$$  \psi _{\cal L}:  \pic C \lra \pic C_1  $$
\noindent  defined by  $ \psi _{\cal L} ({\cal O}(p))  =  \alpha
^*({\cal
L}|_{C_1^p})$, where $C_1^p$  is the fiber in the product over the
point $p \in
C$ and $\alpha : C_1 \lra C_1^p $ the natural isomorphism.    The
extension of
the definition to a point $[L] \in \pic C$ is given by taking a
meromorphic
section of the bundle $L$.  In the same way, whenever we
 have a line bundle on  a fiber product of two families of curves we
get a map
between their
  relative Picard groups

Consider  the diagram

\[
\begin{diagram}
\node{} \node{\fibprod}  \arrow{se,l}{\pi _1} \arrow{s,l}{\pi ^{'}}
\arrow{sw,l}{\pi _2}  \arrow{e,t}{h^{'}}  \node{\breg}   \\
\node{C}  \node{C \times C} \node {C}
\end{diagram}
\]
\vskip.1in

\noindent
Let ${\cal L}$ be a line bundle  on $C \times C$. Then $\pi ^{' *}
{\cal L}$ is
a line bundle on the fiber product $\fibprod $. Following the above
notation,
it is easy to see  that
\begin{lem}
 $\psi _{\pi ^{'*} {\cal L}} = \pi _2 ^* \,  \psi _{\cal L} \,  {\rm
Nm}_1 $,
where ${\rm Nm}_1$ is the norm map of $\pi _1$.
\label{lemblabla1}
\end{lem}

Note that on the level of fibers we have
\[ \begin{array}{cccc}
  \psi _ {\pi ^{'*} {\cal L}}   \,: & \pic \spec   & \lra         &
\pic
\tilde{C}_{\phi  (s) }\\
               & {\cal O}(p)        & \longmapsto  &   \pi  ^*_ {\phi
(s) }
\,  \psi _{\cal L} \, {\rm Nm}_s( {\cal O}(p))
\end{array} \]

\noindent
where the subscripts refer to the restriction on the corresponding
fiber over a
point of $\breg $.

\begin{rem}  {\rm  We  recall the following fact about maps of
abelian torsors
and induced maps of abelian schemes.  Let   $p: {\cal G} \lra Z $ be
an abelian
scheme and
  $\pi : {\cal T} \lra Z$  a ${\cal G}$-torsor.  Given two elements
$t_1$ and
$t_2$   in the same fiber of ${\cal T}$ over $z \in Z$, we denote by
$t _1 - t
_2 $ the unique element $g \in {\cal G}$ over $z \in Z$,  with the
property:
$T_g t_2 =t_1$.  Let now $\phi _{\cal T} : {\cal T} \lra  {\cal T}$
be  an
automorphism of the abelian torsor ${\cal T}$ i.e. an automorphism
that sends a
fiber of the map $\pi $ to another fiber and preserves the action of
${\cal
G}$. To that, one can associate a {\em group}  automorphism  $\phi
_{\cal G}$
of  the abelian  scheme ${\cal G}$ as follows.  Given $g \in {\cal
G}_z$ i.e.
an element of ${\cal G}$ sitting over $z \in Z$, choose an element
$t_z \in
{\cal T}_z$ and define
\[  \phi _{\cal G} (g_z)= \phi _{\cal T} (T_{g_z} t_z) -   \phi
_{\cal T}( t_z)
.      \]
 It is easy to check that this  is independent from the choice of
$t_z$ in the
fiber ${\cal T}_z$ and that it defines a group automorphism.  Note
that given
$t_1$ and $t_2$ two elements in the same fiber  ${\cal T}_z $,  then
 \[   \phi _{\cal T} (t_1) -     \phi _{\cal T}(t_2)    =   \phi
_{\cal G} (t_1
-t_2) . \]
 \noindent }
\label{remtorsor}
\end{rem}

Following the notation of Proposition \ref{prop23}, let $\phitilde :
\unprym
\lra
 \unprym $ be the group automorphism  associated to the automorphism
  $\phitilded :  \unprymd \lra  \unprymd $ of   $\unprym$-torsors.
We determine now the form of the map $\phitilde $.    Consider the
map
\[ \begin{array}{cccc}
   \mu   \,: &    \fibprod  & \lra         &  \unprym  _{\; H}  \!\!
\times
_{h} \unspec \\
               & (p_s, q_{\phi  (s)})     & \longmapsto  &    \left(
\phitilde (rp_s - \pi^*_s {\rm Nm}_s (p_s)),  q_{\phi  (s)}
\right)
\end{array} \]
where the notation for a point on a curve, stands also for the line
bundle
which  the point defines.
 According to  Lemma \ref{lempoincare}, on the product  $\unprym
_{\; H}  \!\!
\times _{h} \unspec$ we have a line bundle $\poinc $ with the
property $\poinc
|_{ [L] \times  \tilde{C}_{\phi (s)}}  \simeq  n \, L $ for some
integer $n$.
Hence,
\begin{equation}
  \mu ^* \poinc |_{ [p_s ] \times  \tilde{C}_{\phi  (s)}}  \simeq  n
\phitilde (rp_s - \pi^*_s {\rm Nm}_s (p_s)).
\label{eqformphi1}
\end{equation}
  By using (\ref{eqformphi1}) and the knowledge of the Picard group
of
$\fibprod$ we will
  derive the form of the map $\phitilde $.    \\
\\
{\bf Case A:} dim$\, \Gamma \cap (\so \times \so ) \leq 2$. Then, by
Lemma
\ref{lemcasea},  we have that  $\pic (\fibprod)  \simeq \pi ^{'*}
\pic (C
\times C)$.
 We thus get
\[  \mu ^* \poinc  \simeq \    \pi ^{'*} {\cal L}   \;\; \mbox{and
so,} \;\;
\mu ^* \poinc|_{[p_s ] \times  \tilde{C}_{\phi (s)}  }  \simeq \
\pi ^{'*}
{\cal L} |_{[p_s ] \times
 \tilde{C}_{\phi (s)}  }  , \]
i.e.
 \begin{equation}
  \phitilde (nrp_s -n  \pi^*_s {\rm Nm}_s (p_s)   )  \simeq   \pi
^*_ {\phi  }
  \psi _{\cal L} {\rm Nm}_s(p_s) .
\label{eqformphi2}
 \end{equation}
\noindent
   Take the map $\phitilde _s :
 \prym  \lra  {\rm Prym}(\tilde{C}_{\phi  (s)} , C)$
  Now given a line bundle  $\tilde{L}_s \in \prym  $, choose  a line
bundle
$\tilde{M}_s$
   in $\prym$ such that
$\tilde{L}_s = nr \tilde{M}_s$.  By choosing a meromorphic section,
we can
write $\tilde{M}_s
= {\cal O}( \sum _i (p_i^1 -p_i^2) )$ where ${\rm Nm}(\sum _i (p_i^1
-p_i^2))=0$. We have

 \[ \begin{array}{rlll}
    \phitilde _s({\cal L}_s)  & =  &  \phitilde _s (rn {\cal M}_s)  =
\phitilde
_s \left( \sum _i rn (p_i^1 -p_i^2) \right)=  &  \\
     & =    &   \sum _i  \left(   \phitilde _s (rn p_i^1- n\pi _s ^*
{\rm
Nm}_s (p_i^1))   -
 \phitilde _s (rn p_i^2- n\pi _s ^*   {\rm Nm}_s (p_i^2)) \right) \\
&   &  + \sum _i \phitilde _s ( n\pi _s ^*   {\rm Nm}_s (p_i^1-
p_i^2))   \\
& =  &  \sum _i  \pi ^*_{\phi  (s)} \psi _{\cal L} {\rm Nm}_s (p_i^1)
- \sum
_i  \pi ^*_{\phi  (s)} \psi _{\cal L} {\rm Nm}_s (p_i^2)  &
\mbox{by} \;\;
(\ref{eqformphi2}),    \\
& =  &  \pi ^*_{\phi  (s)} \psi _{\cal L} {\rm Nm}_s (\sum_i  (p_i^1
-p_i^2)) &
\\
   &  =   &  0 &  \mbox{since} \;\;   {\rm Nm}(\sum _i (p_i^1
-p_i^2))=0.
\end{array} \]
\noindent
The later contradicts the fact that $\phitilde $ is an isomorphism,
which means
that case A cannot occur.  Therefore  dim$\Gamma \cap (\so \times \so
) =  2$,
i.e. we are in case B which we  examine bellow:\\
\\
{\bf Case B:} dim$\Gamma \cap (\so \times \so ) =  2$.   By the
discussion in
Section
 \ref{picfibprod},    the  $\pic (\fibprod)  $ is generated by $\pi
^{'*} \pic
(C \times C)$ and the divisor $\Delta _{\lambda}$. Hence,
$$ \mu ^* \poinc \simeq \pi ^{'*} {\cal L}+ n \Delta _{\lambda}. $$
 Working as before and using the definition of $\Delta _{\lambda}$,
we
conclude that
 $\phitilde  ({\cal L}) = n \psi ^*({\cal L}) $.  Since $\phitilde  $
and
$\psi ^*$  are
 isomorphisms,   we get that  $n=\pm 1$. To summarize,

\begin{prop}
$\phitilde  = \pm  \psi ^*$,  where $\psi $ is the map defined in
Section
\ref{picfibprod}.
\label{propplusminus}
\end{prop}

\subsection{The conclusion of the proof}
\label{conclusionproof}
We now conclude the proof of the main  Theorem \ref{theor1}  in the
case of $r
\geq 3$   by examining the two cases $\phitilde  = \pm \psi ^*$.   We
start
with some notation. We denote by $\psid :  \unprymd \lra \unprymd $
the  pull
back of the map $\psi :  \unspec \lra \unspec $. Note that $\psi  ^*
= \psi
_0^*$. We  have the diagram:
\vskip.1in

 \[
 \begin{diagram}
  \node{\unprymd } \arrow{s,r}{H}  \arrow{e,t} {  ^{ \psi ^*_{
\tilde{d} }} }
 \node{\unprymd  }   \arrow{s,r} {H} \\
 \node{\breg   } \arrow{e,t}{\phi  } \node{\breg }
 \end{diagram}
 \]
 \vskip.1in

\begin{lem}  Let $x$ denote an element  in $\unprymd$. We have that
\begin{enumerate}
\item If $\phitilde = \psi ^* $,  then  $\psi _{\tilde{d}}^{*-1}
\phitilded
(x)-x   $ is independent from $x$ on the fibers of the map $H$.
\item If  $\phitilde = -\psi ^* $,  then  $\psi
_{\tilde{d}}^{*-1}\phitilded
(x)+x   $ is independent from $x$ on the fibers of the map $H$.
\end{enumerate}
\label{lempsi}
\end{lem}

\noindent
{\bf Proof.}  For the first:  Let $y$ be apoint in the  fiber of $H$
through
   $x$.     Since $\psi^{*-1} \phitilde =1$, it is enough to show
that $\psi
_{\tilde{d}}^{*-1} \phitilded (x)  - \psi
_{\tilde{d}}^{*-1}\phitilded (y)
=\psi ^{*-1} \phitilde (x-y)$.  Since $\psi ^ {*-1} \phitilde$ is the
group
homomorphism associated to the map   $\psi _{\tilde{d}}^{*-1}
\phitilded$ of
abelian torsors, the later  is true by  Remark \ref{remtorsor}
For the second: It is enough to show that $\psi _{\tilde{d}}^{*-1}
\phitilded
(x)  - \psi _{\tilde{d}}^{*-1}\phitilded (y) =y-x $. But since
$\phitilde =
-\psi ^* $ we have
  $ \psi  ^{*-1} \phitilde (x-y)=y-x $  and this case follows as
well.
 \begin{flushright}
   $\Box$
   \end{flushright}

\begin{prop} Let ${\bf  -1}$ denote the inversion along the fibers of
$\unprym
$. Then
\begin{enumerate}
\item If $\phitilde = \psi ^* $,  then $\phitilded = T_{\pi ^* \mu}
\circ \psid
$, where $\mu$ is an $r$-torsion line bundle  on the base curve $C$.
\item If  $\phitilde = -\psi ^* $,  then $ r|2d$ and $\phitilded =
T_{\pi ^*
\nu _1} \circ(\psid)
  \circ ({ \bf -1})$, where $ \nu $ is a line bundle on the base
curve $C$
which satisfies
 $ \nu _1 ^{\otimes r}=L_0^{\otimes 2} \otimes  \omega _C  ^{\otimes
r(r-1)} $.
\end{enumerate}
\label{propaaa}
\end{prop}

\noindent
{\bf Proof.}  For the first:  According to  Lemma \ref{lempsi},  on
each fiber
$H^{-1}(s) = {\rm Prym}_{\tilde{d}}(\spec, C) $ of the map $H$,  the
maps
$\phitilded$ and $\psid$ differ by a translation by a unique  element
$a(s) \in
\prym$.  This defines a section of the map $H :  \unjac \lra \breg $
and
therefore by  Proposition \ref{propsections}, it must have
  the form $\pi ^* M$ for some fixed line bundle $\mu $ on $C$. Since
the image
of the section is in the $\unprym  $ we get that $\mu $ is an
$r$-torsion line
bundle.

For the second:  According to  Lemma \ref{lempsi},   we can construct
a section
of the map $H : J^{2 \tilde{d}} (\unspec ) \lra \breg $  by assigning
to the
point $s \in \breg$ the line bundle
$\psi _{\tilde{d}}^{*-1}\phitilded (x)+x   $ for some   point $x \in
{\rm
Prym}_{\tilde{d}}(\spec, C) $.   By Proposition  \ref{propsections},
we get
that $r|2 \tilde{d}$. Since $\tilde{d}= d +r(r-1)(g-1)$ this is
equivalent to
$r|2d$. The same proposition  implies that the above  section must
have the
form $\pi ^* \nu _1$ for some  fixed line bundle $\nu _1$ on $C$.  To
complete
the proof of the second part of the proposition,  observe  that the
map $\psi
^*$ commutes with the translations by an element of the form  $\pi
^*\nu _1$
and that ${\rm Nm} \, \pi ^* \nu _1 =  L_0  ^{\otimes 2}  \otimes
\omega _C
^{\otimes r(r-1)}$.

    \begin{flushright}
   $\Box$
   \end{flushright}

   \bigskip

We are ready now to complete the proof of the Theorem \ref {theor1}
in the case
  $r \geq 3$.   Pick an element  element $s \in \breg $. Then, by
Corollary
\ref{cor17},
 the prymian   $\prymd  $ maps dominantly to  $\su $. Let ${\cal V}$
be its
image.
 It is enough to prove the theorem for the restriction of the map
$\Phi $ on
${\cal V}$.
 We have the following commutative diagram:

\vskip.1in
\[
 \begin{diagram}
 \node{{\rm Prym }_{\tilde {d}}(\spec ,C)} \arrow{s,l}{\pi _{s *}}
\arrow{e,t}
  {\phitilde _{\tilde{d} _{ }  }}
 \node{ {\rm Prym }_{\tilde {d}}(\tilde{C}_{\phi  (s)} ,C) }
\arrow{s,r}
 {\pi _{\phi (s) *} } \\
 \node{{\cal V}  } \arrow{e,t}{\Phi  } \node{{\cal V} }
 \end{diagram}
 \]
 \vskip.1in
\noindent where the maps $\pi _{s *}$ and $\pi _{\phi (s) *}$ are
rational
maps.

Assume first that $\phitilde = \psi ^* $. Let $E$ a vector bundle in
${\cal
V}$.  Then, there exists a line bundle $\tilde {L}_s$ in ${\rm Prym
}_{\tilde
{d}}(\spec ,C)$ such that
 $\pi _{s *} (\tilde{L}_s) = E$. By Proposition \ref{propaaa}, the
above
diagram and
  the fact that the map $\psi $ commutes with the projections    $\pi
_{s *}$
and
 $\pi _{\phi (s) *}$  to the base curve $C$,  we conclude that
$$ \Phi (E) = \Phi  \pi _{s *}  (\tilde{L}_s) = \pi _{\phi (s) *}
\phitilded
(\tilde{L}_s)
 = \pi _{\phi (s) *} (\psid (\tilde{L}_s) \otimes \pi _{\phi (s) *}
\mu ) =\pi
_{\phi (s) *}
 \psid (\tilde{L}_s) \otimes \mu  = E \otimes \mu . $$

Assume next that  $\phitilde =-\psi ^* $. Then,  following the  above
notation,
we have
$$ \Phi (E) = \Phi \pi _{s*} (\tilde{L}_s) = \pi _{\phi (s) *}
\phitilded
(\tilde{L}_s) =
  \pi _{\phi (s) *}  (\psid (\tilde{L}_s^{-1} ) \otimes \pi _{\phi
(s) *}  \nu
_1) =
 \pi _{\phi (s) *} \psid (\tilde{L}_s^{-1}) \otimes \nu _1 = \pi _{s
*}
  (\tilde{L}_s^{-1}) \otimes \nu _1  . $$
We claim that $ \pi _{s *}  (\tilde{L}_s^{-1}) = E^{\vee } \otimes
\omega
_C^{-(r-1)}$. Indeed, consider the map $\pi _s  : \spec \lra C$. By
relative
duality we have
\[ {\cal R}^0 \pi _{s *}  (\tilde{L}_s^{-1})  \simeq {\cal R}^0 \pi
_{s *}
 (\omega _{\pi _s} \otimes \tilde{L}_s  )^{\vee}  .  \]
By the adjunction formula we have, see  e.g. \cite{h1}, that
$\omega_{\pi  _s
}  \simeq
 \pi _s ^* \omega _C ^{r-1}$. We thus  get
\[  \pi _{s *}  (\tilde{L}_s^{-1})  \simeq E^{\vee} \otimes
\omega_C^{-(r-1)}
.  \]
 By choosing $\nu = \nu _1 \otimes  \omega_C ^{-(r-1)}$,  we get that
$  \phi
(E)  = E^{\vee } \otimes  \nu $
 where  $\nu ^{\otimes r} =L_0 ^{\otimes 2}$.

Thus, we have shown the surjectivity of the maps (1.) and (2.) in the
statement
of Theorem
\ref{theor1}. To show the injectivity of the map (1.), pick up a
point $\mu
\neq 0 \in
J^{0}[r]$ and consider the set of fixed points $\su^{\langle T_{\mu}
\rangle}$.
A vector
bundle $E$ is fixed under the action of $T_{\mu}$, if we have an
isomorphism
\[ \kappa : E \longrightarrow E\otimes\mu .\]
 If $p \, | \, r$ is the order of the torsion point $\mu$, then the
$\mu$-twisted Higgs bundle
$(E,\kappa)$ gives a $p$-sheeted unramified spectral cover $\pi_{\mu}
: C_{\mu}
\rightarrow C$ and a semistable vector bundle $F_{\kappa}$ of rank
$r/p$ on it
with the property $\pi_{\mu *}(F_{\kappa}) = E$, see \cite{nr2} for
details.
Therefore the locus
$\su^{\langle T_{\mu} \rangle}$ can be identified with the image of
the moduli
space
${\cal S}{\cal U}_{C_{\mu}}(r/p)$ under the pushforward map $\pi_{\mu
*}$ and
hence is a proper subvariety.

We will sketch the proof for the injectivity of the map (2.) in the
case $r
\geq 3$ and when $L_{0} =
{\cal O}$. The modifiations of the argument for general $L_{0}$ are
minor and
are left
to the reader. If $L_{0} = {\cal O}$, then it suffices to check that
the map $E
\rightarrow
E^{\vee}$ is not the identity on ${\cal S}{\cal U}(r,{\cal O})$. But
if $E$ is
a stable vector
bundle satisfying $E \cong E^{\vee}$, then  the bundle $E^{\otimes
2}$ has a
unique (up to scaling) non-zero section $t$. But $H^{0}(C,E^{\otimes
2}) =
H^{0}(C,{\rm Sym}^{2}E)\oplus H^{0}(C,\wedge^{2}E)$ and therefore
either
 $H^{0}(C,E^{\otimes 2}) = 0$ or $H^{0}(C,\wedge^{2}E) = 0$. This
implies that
the
isomorphism $t$ between $E$ and $E^{\vee}$ is either symmetric or
skew-symmetric. Thus, $E$ is either orthogonal or symplectic and
hence it lies
either in the moduli space  of stable $SO(r)$-bundles or in the
moduli space of
stable $Sp(r)$-bundles. But for $r \geq 3$, those are proper
subvarieties of
${\cal S}{\cal U}(r, {\cal O})$ and this yields  that the map (2.) is
injective
for $r \geq 3$.

\begin{rem} {\rm For $r=2$,  the  moduli space ${\cal S}{\cal U}(2,
{\cal O})$
coincides with the moduli space of $Sp(2)$-bundles; if $E \in  {\cal
S}{\cal
U}(2, {\cal O})$, then $E \simeq E ^{\vee }$}.
\label{remrank2}
\end{rem}

\subsection{The rank $2$ case}
\label{rank2}
The proof in the rank $2$  case is a modification of the proof of the
rank
$\geq 3$ case.
 The main difference is that the linear system defined by
   $B= W_2=H^0(C, \omega _C^2)$   does not separate the
  points on the surface $\so $: A section $s \in B$ corresponds to
the curve
$x^2 + s=0$ on
  $\so $. If $m_{-1}$ is the dilation by $-1$ on $\so $, then a curve
in the
linear system
 $|B|$ that passes through a point $p$,  passes  also through the
point
$m_{-1}(p)$.  Therefore the map $f_{|B|}$ sends $\so $ in a $2:1$ way
to the
projective space. We now show briefly the adjustments  for the rank
$2$ case of
the argument we used  in the rank
  $\geq 3$ case.

 At first one can see e.g.  by  Corollary \ref{cor21},  that the map
$\phi   =
\phi _2$ is a multiplication by a $\lambda \in {\Bbb C}^*$.

   Following a similar argument with that of Section \ref{picunspec},
we can
prove that the Picard group of $\unspec \lra \breg $ is again
$$  \pic \unspec = \pi ^* \pic C  . $$
 We also have
\begin{equation}
\pic \unspec _p =  \pic  \pi ^* C  \oplus
{\Bbb Z} [H_p]
\label{eqrank21}
\end{equation}
but in this case the technicalities of the proof are slightly
different: Let
$p^{'} = m_{-1}(p)$. On the fiber of $\so $ over $c= \alpha(p)
=\alpha
(p^{'})$, only the points $p $ and $p^{'}$ belong in the image of the
map
$\beta $. We write $\unspec _p = \unspec _p\setminus (H_p \cup
H_{p^{'}})
\amalg  (H_p \cup H_{p^{'}})$.   Then $\pic \unspec _p$ is generated
by
$\pic(C\setminus  \{ c \})$ and $H_p, H_{p^{'}}$. Observe that $\pi
^*(c) =
H_p+ H_{p^{'}}$. We  thus  have that $\pic \unspec _p $ is generated
by $\pic
\pi ^* C $ and $H_p$.  To prove that this is a direct sum,
  we again take  a pencil as in the proof of Proposition \ref
{proppicunspecp},
but now the curve at infinity embedded on $S$,
 consists of the  $Y_{\infty }$ and a bunch of $N_1$ fibers  over the
points
$c_1, \ldots , c_{N_1}$ on C,  which pass through the $N= 2N_1$ base
points of
the pencil.  Therefore the fiber over the infinity point $b$ on $X$
consists of
the infinity divisor $\tilde{Y}_{\infty}$ and a bunch of $N_1$
divisors $A_1,
\ldots , A_{N_1}$ which satisfy $\pi  ^* (c_i) = A_i + E_i +
E^{'}_i$.  Using
that we get
$$ \pic  X^{\circ } = \pi ^* \pic C \oplus _{i=1}^{N_1} ( \oplus
{\Bbb Z}[E_i])
  \oplus _{i=1}^{N_1} (\oplus {\Bbb Z}[E^{'}_i] )  \left/ (\pi
^*(c_i)=E_i
+E^{'}_i)_{i=1, \ldots N_1}
  \right.     .  $$
 Working as in the rank $\geq 3$ case,  this proves relation
(\ref{eqrank21}).
To prove  the analogue of Proposition \ref{propsections},  we proceed
in the
 same way  as in the rank $\geq 3$ case.

For the Picard group of $\fibprod $: Let again $\psi :\unspec \lra
\unspec $ be
the analogue of the map introduced in Section \ref{picfibprod}. Then
the locus
$\Gamma $ consists of two components namely $\Gamma _1 = \{ (p, \psi
(p)) \;\;
\mbox{for} \;\; p \in \unspec  \}$
and  $\Gamma _2 = \{ (p, \psi (m_{-1}(p))) \;\; \mbox{for} \;\; p \in
\unspec
\}$. By
 observing that $\pi ^{'*} (\mbox{Diagonal in} \;\; C \times C ) =
\Gamma _1 +
 \Gamma _2 $ we conclude that $\pic (\fibprod )$ is generated by
$\pic (C
\times C)$
  and $\Gamma _1$. The rest of the argument, combined with Remark
\ref{remrank2},
  proceeds  as in the rank $r \geq 3 $ case.

\subsection{Curves with automorphisms}
\label{automorphisms}
For the rank $r \geq 3$ case, the only modification that has to be
done,  is in
the argument of Section \ref{picfibprod}, about the induced
automorphism of the
embedded surface $f_{|\overline{B}  |}$.   For the rank $2$ case, the
only
modification is  in the use of
Corollary  \ref{cor21}. We leave to the reader to fill up the details
for the
proof  of Theorem \ref{theor2}.
\section{The automorphisms of $\urd $}
\label{urd}
\begin{prop}
Let $\Phi $ be an automorphism of $\urd $. Then $\Phi $ factors
through the
determinant map i.e. we have the following commutative diagram
\begin{equation}
\begin{diagram}[aaaaaaa]
 \node{\urd   }   \arrow{e,t}{\Phi } \arrow{s,l} {\rm det}
 \node{ \urd  }    \arrow{s,r} {\rm det }\\
\node{J^d(C) } \arrow{e,t}{ \phi _d} \node{ J^d(C) }
\end{diagram}
\label{diagurd1}
\end{equation}
\label{propurd1}
\end{prop}
\noindent
{\bf Proof.}  We would like  to show that the map ${\rm det} \circ
\phi : {\rm
det}^{-1}(L_0) \lra J^d(C)$ is constant for all $L_0 \in J^d(C)$. By
definition,  ${\rm det}^{-1}(L_0) = \su $ and ${\rm dim} \, \su =
(r^2-1)(g-1)
> g, \;\; g \geq 2 $. Therefore the map  ${\rm det} \circ \phi :
 \su  \lra J^d(C)$ has positive dimentional fibers and so, no  pull
back of a
line bundle
from $J^d(C)$ can be ample.  According to  \cite{dn},  $\pic \su =
{\Bbb
Z}[\Theta ]$. Since
 $\su $ is a projective variety, $\Theta $ is an ample divisor. Let
${\cal L}$
be a line bundle on $J^d(C)$. We thus have  $({\rm det} \circ \phi)
^*{\cal L}
\simeq  n \Theta$ for some integer $n$. We claim that $n=0$: If not,
then  say
first that $n >0$. Then $  n \Theta$ is ample and so $({\rm det}
\circ \phi)
^*{\cal L} $ is ample, which is  a contradiction.  Next, if
 $n<0$,  then the same argument applied to the line bundle ${\cal
L}^{-1}$
leads to a
 contradiction. Therefore $({\rm det} \circ \phi) ^*{\cal L} ={\cal
O}$ for all
${\cal L} \in \pic J^d(C)$. Since $\su $ is irreducible, this implies
that the
map is constant. Indeed, if the image has  positive dimension, then
there
exists a positive divisor on that - the hyperplane section of its
embedding
in the projective space. But then the pull back of that divisor on
$\su $ by
the map ${\rm det} \circ \phi $ defines a non - trivial line bundle,
which is a
contradiction.
 \begin{flushright}
   $\Box$
   \end{flushright}

   \bigskip

As in the case of the automorphisms of $\su $, an automorphism  $\Phi
$ of
$\urd $ induces an automorphism $\phitilded$ of the Jacobian
fibration $H:
\unjacd \lra W^{\rm reg}$ which
 sends a fiber of the Hitchin map $H$ to  a fiber of $H$.  Let
  $N_1 : \unjacd \lra J^d(C)$ be the map defined by $N_1 = {\rm det}
\circ \pi
_*$.
  Note that ${\rm Nm} = T_{\frac{r(r-1)}{2} \omega _C } \circ N_1$.

\begin{lem}
Let $s $ be a  point in $W^{\rm reg}$. Then the following diagram
commutes
\[
\begin{diagram}
 \node{J^{\tilde {d}}(\spec ) }   \arrow{e,t}{\phitilde _{\tilde{d}
_{   }}  }
 \arrow{s,l} {\rm N_1}
 \node{ J^{\tilde {d}}( \tilde{C}_{\phi _W (s)})  }    \arrow{s,r}
{\rm N_1 }\\
\node{J^d(C) } \arrow{e,t}{ \phi _d} \node{ J^d(C) }
\end{diagram}
\]
\label{lemurd1}
\end{lem}
\noindent
{\bf Proof.}  For the proof we are going to use the following
commutative
diagram

\begin{equation}
\begin{diagram}
 \node{{\cal X}(r,d) }   \arrow{e,t}{d\Phi ^*  } \arrow{s,l} {\pi _*}
 \node{{\cal X}(r,d)   }    \arrow{s,r} {\pi _*}\\
\node{\urd }  \arrow{e,t}{\Phi }  \arrow{s,l} {{\rm det}}
\node{\urd}
\arrow{s,r}
{{\rm det}}\\
\node{J^d(C) } \arrow{e,t}{ \phi _d} \node{ J^d(C) }
\end{diagram}
\label{diagurd3}
\end{equation}

\noindent
It is enough to show that $N_1 \phitilded = \phi _d N_1 $ on a
Zariski open $U$
 in  $ J^{\tilde {d}}(\spec )  $. Choose $U$ to be the intersection
of the
cotangent bundle ${\cal X}(r,d) $ to $\urd $ with the  Jacobian
$J^{\tilde
{d}}(\spec )$.
 According to Corollary \ref{cor17},  this is a non empty Zariski
open in
$\urd$.
The proof of the Lemma is now a consequence of the commutativity of
the
 above diagram.

 \begin{flushright}
   $\Box$
   \end{flushright}

   \bigskip

Let $\phitilde $ and $\phi $ be the group maps associated  to
$\phitilded $ and
$\phi _d$,
 see Remark \ref{remtorsor}.   By using the above Lemma
\ref{lemurd1},  it is
easy to see that
\begin{cor}
The following diagram is commutative
\begin{equation}
\begin{diagram}
 \node{J^0 (\spec ) }   \arrow{e,t}{\phitilde } \arrow{s,l} {\rm Nm}
 \node{ J^0( \tilde{C}_{\phi _W (s)} ) }    \arrow{s,r} {\rm Nm }\\
\node{J^0(C) } \arrow{e,t}{ \phi  } \node{ J^0(C) }
\end{diagram}
\label{diagurd4}
\end{equation}
where ${\rm Nm}$ is the norm map.
\label{corurd1}
\end{cor}

\begin{lem}
Following the notation of Corollary \ref{corurd1} , the  diagram
bellow  is
commutative
 \[
\begin{diagram}
 \node{J^0 (\spec ) }   \arrow{e,t}{\phitilde }
 \node{ J^0( \tilde{C}_{\phi _W (s)}  )}    \\
\node{J^0(C) }\arrow{n,l} {\pi ^*}\arrow{e,t}{ \phi  } \node{ J^0(C)
}
\arrow{n,r}{\pi ^*}
\end{diagram}
\]
\label{lemurd2}
\end{lem}
\noindent
{\bf Proof.}  Note first that
\begin{equation}
\phitilde \pi ^* (L) = \pi ^* (M)  \;\; \mbox{for some fixed line
bundle} \;\;
M .
\label{equrd1}
\end{equation}
Indeed,
 $\phitilde $   is a global automorphism on $\unjac $ and so,
$\phitilde \,
\pi ^* (L)$ defines  a section of the map $H: \unjac \lra W^{\rm
reg}$.  Hence,
  by Proposition \ref{propsections},  it must have  the form $ \pi ^*
(M) $ for
some fixed line bundle
  $M$ on $C$.  By applying the norm map  to (\ref{equrd1}), we get
${\rm Nm}
  \phitilde \pi ^* (L)
  ={\rm Nm} \, \pi ^* (M) $.
 Corollary \ref{corurd1} implies that   $r \phi (L) =r M$.  Let
$n_r$ denote
the multiplication by $r$. By composing both sides of  (\ref{equrd1})
by $n_r$
and by using
  the last relation we get  that $n_r \phitilde \pi ^* (L) = n_r \pi
^*\phi
(L) $.
 Since the map $n_r $ is onto, this  completes the proof.
 \begin{flushright}
   $\Box$
   \end{flushright}

\begin{lem}  For any $E \in \urd $ and $M \in J^0(C)$ we have
$  \Phi (E \otimes M)= \Phi (E) \otimes \phi (M)$.
\label{lemurd3}
\end{lem}
\noindent
{\bf Proof.}  It suffices to prove the relation for all $E$ in a
Zariski open
$V$ of $\urd $. Take $\spec $ a  smooth spectral curve  and let $V$
be the
image of $ {\cal X}(r,d)
  \cap J^{\tilde{d}}(\spec )$ on $\urd$. According to  Corollary
\ref{cor17},
this is a Zariski open  of $\urd $.  Then $E=\pi _* \tilde {L} $. We
have

 \[ \begin{array}{rlll}
  \Phi (E \otimes M)   & = &  \Phi (\pi _* \tilde {L} \otimes M)
\\
     & =     &\Phi (\pi _* (\tilde {L} \otimes \pi ^* M))        &
\mbox{by
the projection formula},  \\
       & =       &    \pi _* \phitilded  (\tilde{L} \otimes \pi ^* M)
&
\mbox{by diagram (\ref {diagurd3})},
           \\
            & =    & \pi _*  (\phitilded(\tilde{L}) \otimes
\phitilde (\pi ^*
M)    )   &  \mbox {by the definition of the group map} \;\;
\phitilde ,  \\
           &  =    & \pi _*  (\phitilded(\tilde{L}) \otimes  \pi ^*
\phi ( M)
  )        & \mbox{by
 Lemma \ref{lemurd2}},  \\
  &  =   &    \pi _* \phitilded (\tilde{L})  \otimes \phi (M)   &
\mbox{by the
projection
formula},\\
& = & \Phi \pi _* (\tilde {L})  \otimes \phi (M)     \\
   &   =  & \Phi (E) \otimes \phi (M).
\end{array} \]
  \begin{flushright}
   $\Box$
   \end{flushright}

   \bigskip

\noindent
{\bf Proof of  Theorem \ref{theor3}.}   Let  $\Phi :\urd \lra \urd$
be a given
automorphism
 and let $\phi _d$ the  induced automorphism on $J^d(C)$ as in Lemma
\ref{lemurd1}. Given a vector bundle $E \in \urd $, we can find a
vector bundle
$E_{L_0} \in \su $ and a line bundle $\eta  \in J^0(C)$ such that $E=
E_{L_0}
\otimes \eta $. Note that $  \eta ^{\otimes r} = det E \otimes  L_0
^{-1}$. The
above  decomposition is unique up to a choice of an $r$-torsion
point.  Let
$\xi _1$ be a line bundle such that $ \xi _1 ^{ \otimes r}=
\phi _d(L_0)  \otimes  L_0 ^{-1} $.  Then,   $ T_{ \xi _1 ^{ -1} }
\circ  \Phi
$ induces an automorphism of $\su $. By the main Theorem
\ref{theor2}, we have
two cases.\\
\\
{\bf Case 1:}   $ T_{ \xi _1 ^{ -1} } \circ  \Phi (E_{L_0}) = \sigma
^* E_{L_0}
 \otimes \mu $, where $\sigma $ is an automorphism of the curve $C$
and $\mu $
a line bundle which satisfies  $ \mu ^{ \otimes  r} = L_0  \otimes
\sigma ^*
L_0 ^{-1}$.
 We choose now  $ \xi = \xi _1  \otimes\mu $  and so, $\xi ^{\otimes
r} =\phi
_d(L_0)
 \otimes    \sigma ^* L_0 ^{-1} $.   We have
 \[   \begin{array}{rlll}
  \Phi (E)      & = &  \Phi (E_{L_0} \otimes \eta )      \\
     & =     & \Phi ( E_{L_0}) \otimes \phi (\eta )       &
   \mbox{by  Lemma \ref{lemurd3}},   \\
       & =       &  T_{\xi _1} (T_{ \xi _1^{-1}  } \Phi (E_{L_0}))
\otimes
\phi (\eta )    &              \\
            & =    &      \sigma ^* E_{L_0}  \otimes \xi  \otimes
\phi  (\eta )
            &   \\
           &  =    &     \sigma ^* (E \otimes \eta  ^{-1})  \otimes
\xi
\otimes \phi  (\eta )   & \mbox{by the definition of} \;\;
       E_{L_0} , \\
  &  =   &   \sigma ^* E \otimes \xi \otimes \phi (\eta ) \otimes
\sigma
^*(\eta ^{-1} )  .
\end{array}     \]
Therefore,
\[
 \Phi (E) =   \sigma ^* E \otimes \xi \otimes \phi  (\eta ) \otimes
\sigma
^*(\eta  ^{-1})  \;\; \mbox {where} \;\;
   \eta ^{\otimes r} = {\rm det} E \otimes L_0 ^{-1}\;\; \mbox{and}
\;\;  \xi
^{\otimes r} = \phi _d(L_0)  \otimes  \sigma ^* L_0 ^{-1} .
\]
  Note that, due to Lemma \ref{lemurd3}, the map $\phi ^{-1} \circ
\sigma ^*
$ is the identity on the set of $r$-torsion points in $J^0(C)$.\\
\\
{\bf Case 2:}   $ T_{ \xi _1^{ -1} } \circ  \Phi (E_{L_0}) = \sigma
^*
E_{L_0}^{\vee}   \otimes \nu $, where $\sigma $ is an automorphism of
the curve
$C$ and $\nu $ is a  line bundle which satisfies
$  \nu ^{\otimes r}=  L_0  \otimes \sigma ^* L_0 $.  In this case we
choose
$\xi = \xi _1 \otimes \nu $  and so,  $ \xi ^{\otimes r} = \phi _d
(L_0)
\otimes \sigma ^* L_0 $.
      We have
 \[ \begin{array}{rlll}
  \Phi (E)       & =     & \Phi ( E_{L_0}) \otimes \phi (\eta )  &
\mbox{by
Lemma \ref{lemurd3}},   \\
       & =       &  T_{\xi _1 } (T_{ \xi _1^{ -1}  } \Phi (E_{L_0}))
\otimes
\phi (\eta )    &              \\
            & =    &      \sigma ^* E_{L_0}^ {\vee}   \otimes \xi
\otimes \phi
 (\eta )             &   \\
           &  =    &     \sigma ^* (E^{\vee } \otimes \eta  )
\otimes \xi
\otimes \phi  (\eta )  & \mbox{by the definition
  of} \;\;     E_{L_0} ,\\
  &  =   &      \sigma  ^* E^{\vee}    \otimes \xi \otimes \phi (\eta
) \otimes
  \sigma ^*(\eta ) .
\end{array} \]
Therefore,
\[
 \Phi (E) =   \sigma ^* E ^{\vee} \otimes \xi \otimes \phi (\eta )
\otimes
\sigma ^*(\eta )  \;\; \mbox {where} \;\;
    \eta ^{\otimes r}  = {\rm det} E \otimes L_0^{-1}\;\; \mbox{and}
\;\;  \xi
^{\otimes r} = \phi _d(L_0)  \otimes \sigma ^*  L_0 .
\]
  Again,  due to Lemma \ref{lemurd3}, the map $\phi \circ  \sigma ^*
$ has to
be the identity  map on the set of $r$-torsion points in $J^0(C)$.

\appendix
\section{Appendix: A proof of the Torelli Theorem for the moduli
space of
vector bundles}
Using the technics of the paper we sketch in the following the proof
of Theorem
\ref{Torelli} - the Torelli Theorem for vector bundles.
\begin{rem}
{\rm In the case $(r,d)=1$ the theorem has been proven in \cite{mn},
\cite{nr}
and \cite{ty}. The way they prove it is by showing that the
intermediate
Jacobian of  the moduli space is canonically isomorphic to the
principally
polarized Jacobian of the curve; the theorem then follows from the
usual
Torelli for Jacobian of curves. In the case $(r,d) \neq 1$, the
moduli space of
vector bundles is singular and the construction of the intermediate
Jacobian
does not apply. In the special case $r=2$ and trivial determinant,
Balaji has
proven a similar result for a desingularization of the space
constructed by
Seshadri, see \cite {ba}. Here we prove the theorem for any  $r, d$.}
\end{rem}
\noindent
{\bf Proof of Theorem \ref{Torelli}.}  The notation we use is in
analogy with
the
 one used in the paper. The
 subscripts $_1$, $_2$ in the notation refer to the curves $C_1$,
$C_2$
respectively.
   Let $\Phi : {\cal S}{\cal U}_{C_1}(r, L_1) \lra
 {\cal S}{\cal U}_{C_2}(r, L_2)$ be the given isomorphism. After
lifting to the
cotangent bundle we obtain  in a similar way  as in  the paper the
following
commutative diagram:
 \[
  \begin{diagram}
  \node{{\rm Prym}(\tilde{\cal C}_1,C_1)} \arrow{e,t}{\phitilde }
\arrow{s,l}{H_1}
\node{{\rm Prym}(\tilde{\cal C}_2,C_2)}
  \arrow{s,r}{H_2} \\
   \node{B_1^{\rm reg}} \arrow{e,b}{\phi } \node{B_2^{\rm reg}}
\end{diagram}
\]
where $\phitilde $ is the induced isomorphism,  $H_1$, $H_2$ are the
Hitchin
maps and the map $\phi $ is a linear isomorphism. Consider the
following
diagram
\[
  \begin{diagram}
  \node{S_2} \arrow{e,t}{i_2 }
\node{{\Bbb P}( \overline{B}_2^{ \vee})    }
  \arrow{s,r}{\overline{\phi}^{\, *}} \\
   \node{S_1} \arrow{e,b}{i_1 } \node{{\Bbb P}(\overline{B}_1^{\vee})
}
\end{diagram}
\]
where the maps $i_1, i_2 $ are those defined by the linear systems
$\overline{B}_1, \overline{B}_2 $ respectively and the isomorphism
$\overline{\phi} ^{\, *}$ is the one induced by the extension
$\overline{\phi}$
of  $\phi $.
 Using  arguments similar to those in Sections \ref{picfibprod} and
\ref{formphi}, we deduce that the existence of an isomorphism
$\phitilde $ in
the above diagram, implies that  $\overline{\phi} ^{\, *}$ has to map
$i_2(S_2)$ isomorphically to $i_1(S_1)$. Thus the curve $C_2$ maps to
$C_1$ and
since both are curves of the same genus $g \geq 3$, the map is an
isomorphism.

 \begin{flushright}
   $\Box$
   \end{flushright}

\end{document}